\documentclass[smallextended]{svjour3}       % onecolumn (second format)
\smartqed  % flush right qed marks, e.g. at end of proof
\usepackage{xcolor}

\definecolor{matlab_blue}{rgb}{        0 ,   0.4470 ,   0.7410}
\definecolor{matlab_red}{rgb}{1.0,    0.0,    0.0}
\definecolor{matlab_yellow}{rgb}{0.9290,    0.6940 ,   0.1250}
\definecolor{matlab_green}{rgb}{0.4660,    0.6740,    0.1880}
\usepackage{graphicx}
\usepackage{subfigure}

\usepackage{amssymb,amsmath,sansmath,lmodern}
\usepackage{color}

\usepackage{accents}
\usepackage{stmaryrd}
\usepackage{tabu}

\usepackage[numbers]{natbib}

\newcommand{\average}[1]{\left\{\!\!\left\{#1\right\}\!\!\right\}}   % average at element boundary

\newcommand{\aver}[1]{\ensuremath{\{\!\{#1\}\!\}}}
\newcommand{\jump}[1]{\ensuremath{\left\llbracket #1 \right\rrbracket}}

\newcommand\R{\mathbb R}

\newcommand\Rey{\mathrm{Re}}
\newcommand\Fr{\mathrm{Fr}}
\newcommand*\diff{\mathop{}\!\mathrm{d}}

\DeclareMathAccent{\svec}{\mathord}{letters}{126}
\newcommand\acclrvec[1]{\accentset{\,\leftrightarrow}{#1}}
\newcommand\stvec[1]{\mathbf #1}				
\newcommand\stvecg[1]{\boldsymbol #1}				
\newcommand\ssvec[1]{\acclrvec{\stvec{#1}}}	
\newcommand\ssvecg[1]{\acclrvec{\stvecg{#1}}}	
\newcommand\cssvec[1]{\acclrvec{\tilde{\stvec{#1}}}} 
\newcommand\cssvecg[1]{\acclrvec{\tilde{\stvecg{ #1}}}} 
\newcommand{\smat}[1]{\underline{\stvec{#1}}}

\newcommand{\snabla}{\svec{\nabla}}
\newcommand\interiorfaces{{\mathrm{interior}\atop\mathrm{faces}}}
\newcommand\boundaryfaces{{\mathrm{boundary}\atop\mathrm{faces}}}

\usepackage{xargs}    
\usepackage[colorinlistoftodos,prependcaption,textsize=tiny]{todonotes}
\newcommandx{\unsure}[2][1=]{\todo[linecolor=blue,backgroundcolor=blue!25,bordercolor=blue,#1]{#2}}
\newcommandx{\change}[2][1=]{\todo[linecolor=red,backgroundcolor=red!25,bordercolor=red,#1]{#2}}
\newcommandx{\info}[2][1=]{\todo[linecolor=OliveGreen,backgroundcolor=OliveGreen!25,bordercolor=OliveGreen,#1]{#2}}
\newcommandx{\improvement}[2][1=]{\todo[linecolor=Plum,backgroundcolor=Plum!25,bordercolor=Plum,#1]{#2}}
\newcommandx{\thiswillnotshow}[2][1=]{\todo[disable,#1]{#2}}

\usepackage{tikz}
\usetikzlibrary{plotmarks}
\newcommand\marksymbol[2]{\tikz[#2,scale=1.2]\pgfuseplotmark{#1};}

\begin{document}

\title{Entropy--stable discontinuous Galerkin approximation with summation--by--parts property for the incompressible
Navier--Stokes equations with variable density and artificial compressibility}
%\subtitle{Do you have a subtitle?\\ If so, write it here}

\titlerunning{ES DG approximation with SBP for the inc.
	NSE with artif. compressibility}        % if too long for running head

\author{Juan Manzanero        \and Gonzalo Rubio \and David A. Kopriva \and Esteban Ferrer \and 
Eusebio Valero 
}

%\authorrunning{Short form of author list} % if too long for running head

\institute{Juan Manzanero (\email{juan.manzanero@upm.es}), Gonzalo Rubio, Esteban Ferrer, Eusebio Valero  \at
	ETSIAE-UPM - School of Aeronautics, Universidad Polit\'ecnica de Madrid. Plaza Cardenal Cisneros 3, E-28040 Madrid, Spain. //
	Center for Computational Simulation, Universidad Polit\'ecnica de Madrid, Campus de Montegancedo, Boadilla del Monte, 28660, Madrid, Spain. \\
              David A. Kopriva \at Department of Mathematics, Florida State University and Computational Science Research Center, San Diego State University. 
}

\date{Received: date / Accepted: date}
% The correct dates will be entered by the editor

\maketitle

\begin{abstract}
We present a provably stable discontinuous Galerkin spectral element method for the incompressible Navier--Stokes equations with artificial compressibility and variable density. Stability proofs, which include boundary conditions, that follow a continuous entropy analysis are provided. We define a mathematical entropy function that 
combines the traditional kinetic energy and an additional energy term for the artificial compressiblity, and derive its associated entropy 
conservation law. The latter allows us to construct a provably stable split--form nodal Discontinuous 
Galerkin (DG) approximation that satisfies the summation--by--parts simultaneous--approximation--term (SBP--SAT) 
property.  The scheme and the stability proof are presented for general curvilinear three--dimensional hexahedral meshes.
We use the exact Riemann solver and the Bassi--Rebay 1 (BR1) scheme at the 
inter--element boundaries for inviscid and viscous fluxes respectively, and an explicit low storage Runge--Kutta RK3 scheme 
to integrate in time. We assess the accuracy and robustness of the method by solving the Kovasznay flow, the inviscid Taylor--Green vortex, and 
the Rayleigh--Taylor instability.
\keywords{Navier--Stokes \and Computational fluid dynamics \and High-Order methods \and Discontinuous Galerkin \and Spectral element.}
% \PACS{PACS code1 \and PACS code2 \and more}
% \subclass{MSC code1 \and MSC code2 \and more}
\end{abstract}

\section{Introduction}\label{sec:Introduction}

In this work, we provide an entropy--stable framework to solve the incompressible 
Navier--Stokes Equations (NSE) with varying density. The method presented herein 
can be used as a fluid flow engine for interface--tracking multiphase flow models  
(e.g. VOF \cite{1981:Hirt,1980:Nichols}, level--set \cite{1994:Sussman,1995:Adalsteinsson}, 
phase--field models \cite{1999:Jacqmin,2003:Badalassi}). 

Amongst the different incompressible NSE models, we use the artificial 
compressibility method \cite{1997:Shen}, which converts the elliptic problem into a 
hyperbolic system of equations, at the expense of a non--divergence free velocity 
field. However, it allows one to avoid constructing an approximation that satisfies the inf--sup condition \cite{2012:Ferrer,2013:Karniadakis,2014:Ferrer,2017:Ferrer}. 
Artificial compressibility is commonly combined with dual timestepping, which solves an inner pseudo--time step loop until velocity divergence errors are lower 
than a selected threshold, then performs the physical time marching \cite{2016:Cox}. 
In this work, we only address the spatial discretization. However, the method presented herein 
can be complemented with the pseudo--time step to control divergence 
errors. Nonetheless, we have found that solving the incompressible NSE with 
artificial compressibility can obtain satisfactory results even in transient 
simulations.

In this paper we present a nodal Discontinuous Galerkin (DG) spectral element method (DGSEM) 
for the incompressible Navier--Stokes equations with artificial compressibility. In particular, 
this work uses the Gauss--Lobatto version of the DGSEM, which makes it possible to 
obtain entropy stable schemes using the summation--by--parts simultaneous--approximation--term 
(SBP--SAT) property and two--point entropy conserving fluxes. Moreover, it handles arbitrary three 
 dimensional curvilinear hexahedral meshes while maintaining high order, spectral accuracy and 
 entropy stability.
 
 We present a novel entropy analysis for the incompressible NSE with artificial compressibility and variable
 density, where the traditional kinetic energy is complemented with an 
 artificial compressibility energy that forms the mathematical entropy. The 
 entropy conservation law is then mimicked semi--discretely, i.e. only considering 
 spatial discretization errors. The approximation uses a split--form DG 
 \cite{2016:gassner,2017:Gassner}, with the exact Riemann solver 
 \cite{2017:Bassi}, and the Bassi--Rebay 1 (BR1) \cite{1997:Bassi} to compute 
 inter--element boundary fluxes. We complete the analysis with a stability study of solid wall 
 boundary conditions. As a result, the numerical scheme is 
 entropy stable and parameter--free. 
 
 The rest of this paper is organized as follows: In Sec. \ref{sec:GoverningEquations} 
 we introduce the incompressible NSE with variable density and artificial 
 compressibility, and we perform the continuous entropy analysis in Sec. 
 \ref{subsec:GoverningEquations:Stability}. In Sec. \ref{sec:DGSEM} we describe
 the split--form DG scheme, and in Sec. \ref{sec:DiscreteEntropyAnalysis} we 
 study its entropy stability. Lastly, we perform numerical experiments in Sec. 
 \ref{sec:NumericalExperiments}. A convergence study using manufactured 
 solutions in Sec. \ref{subsec:num:conv}, the solution of the Kovasznay flow problem 
 in Sec. \ref{subsec:num:Kovasznay}, the inviscid three--dimensional Taylor--Green 
 vortex problem in Sec. \ref{subsec:num:TGV}, and the Rayleigh--Taylor 
 instability in Sec. \ref{subsec:num:RTI}.

\section{Governing equations. Continuous entropy analysis}\label{sec:GoverningEquations}
Given velocity $\svec{u}(\svec{x},t)=(u_1,u_2,u_3)=(u,v,w)$, pressure $p\left(\svec{x},t\right)$, and density $\rho(\svec{x},t)$ fields, the 
incompressible Navier--Stokes Equations (NSE) consist of the momentum 
equation,
\begin{equation}
\left(\rho \svec{u}\right)_t + \nabla\cdot\left(\rho\svec{u}\svec{u}\right) = 
-\nabla p + 
\nabla\cdot\left(\frac{1}{\Rey}\left(\nabla\svec{u}^{T}+\nabla\svec{u}\right)\right)+\frac{1}{\Fr^2}\rho\svec{e}_{g},
\label{eq:governing:iNS-momentum}
\end{equation}
with $\Rey$ and $\Fr$ being the Reynolds and Froude numbers, respectively, and $\svec{e}_{g}$ a unit vector in the gravity
direction. The artificial compressibility method \cite{1997:Shen} adds an equation for the time evolution of the pressure,
\begin{equation}
p_t + \frac{1}{M_0^2}\nabla\cdot\svec{u}=0,
\label{eq:governing:iNS-artificialcompressibility}
\end{equation}
where $M_0$ is the artificial compressibility model Mach number. 
Eqs. \eqref{eq:governing:iNS-momentum} and \eqref{eq:governing:iNS-artificialcompressibility} can be 
augmented with a transport equation for the density, which we allow to vary spatially,
\begin{equation}
  \rho_t + \nabla\cdot\left(\rho\svec{u}\right)=0.
  \label{eq:governing:iNS-densitytransport}
\end{equation}

Gathering
\eqref{eq:governing:iNS-momentum}--\eqref{eq:governing:iNS-densitytransport}, we 
regard the incompressible NSE with artificial compressibility as a hyperbolic 
system,

\begin{equation}
  \stvec{q}_t + \sum_{i=1}^{3}\frac{\partial\stvec{f}_{e,i}(\stvec{q})}{\partial x_i} = 
\sum_{i=1}^{3}\frac{\partial\stvec{f}_{v,i}(\stvec{q},\nabla\stvec{q})}{\partial x_i}+\stvec{s}\left(\stvec{q}\right),
  \label{eq:governing:iNS-system-sum}
\end{equation}
with conservative variables $\stvec{q}=\left(\rho,\rho\svec{u},p\right)$, inviscid 
fluxes $\stvec{f}_{e,i}(\stvec{q})$,

\begin{equation}
\stvec{f}_{e,1}=\stvec{f}_{e} = \left(\begin{array}{c}\rho u \\ \rho u^2+p \\ \rho u v \\ \rho u w \\ \frac{1}{M_0^2}u \end{array}\right),~~
\stvec{f}_{e,2}=\stvec{g}_{e} = \left(\begin{array}{c}\rho v \\ \rho uv \\ \rho v^2+p \\ \rho v w \\ \frac{1}{M_0^2}v \end{array}\right),~~
\stvec{f}_{e,3}=\stvec{h}_{e} = \left(\begin{array}{c}\rho w \\ \rho uw \\ \rho vw \\ \rho w^2+p \\ \frac{1}{M_0^2}w \end{array}\right),
   \label{eq:governing:iNS-inviscid-flux}
\end{equation}
viscous fluxes $\stvec{f}_{v,i}(\stvec{q},\nabla\stvec{q})$,

\begin{equation}
\stvec{f}_{v,1}=\stvec{f}_{v} = \left(\begin{array}{c}0 \\ \tau_{11} \\ \tau_{21} \\ \tau_{31} \\ 0 \end{array}\right),~~
\stvec{f}_{v,2}=\stvec{g}_{v} = \left(\begin{array}{c}0 \\ \tau_{12} \\ \tau_{22} \\ \tau_{32} \\ 0 \end{array}\right),~~
\stvec{f}_{v,3}=\stvec{h}_{v} = \left(\begin{array}{c}0 \\ \tau_{13} \\ \tau_{23} \\ \tau_{33} \\ 0 \end{array}\right),
   \label{eq:governing:iNS-viscous-flux}
\end{equation}
and source term $\stvec{s}\left(\stvec{q}\right)=\left(0,\frac{1}{\Fr^2}\rho\svec{e}_{g},0\right)$. In \eqref{eq:governing:iNS-viscous-flux}, we used the viscous tensor $\tens{\tau}$
\begin{equation}
\tau_{ij}=\frac{2}{Re}\tens{S}_{ij},
\label{eq:governing:stress-tensor}
\end{equation}
where
\begin{equation}
\tens{S}_{ij}=\frac{1}{2}\left(\frac{\partial u_i}{\partial x_j}+\frac{\partial u_j}{\partial x_i}\right)
\label{eq:governing:strain-tensor}
\end{equation}
are the components of the strain tensor, $\tens{S}$.

In this work we adopt the notation in \cite{2017:Gassner} to distinguish vectors with different nature. With an arrow on top of the variable, we define space vectors (e.g. $\svec{x}=\left(x,y,z\right)\in\R^{3}$). Vectors in bold are state vectors (e.g. $\stvec{q}=\left(\rho,\rho\svec{u},p\right)\in\R^{5}$). Lastly, we define block vectors as the result of stacking three state vectors spatial coordinates (e.g. fluxes),
\begin{equation}
\ssvec{f}_{e} = \left(\begin{array}{ccc}\stvec{f}_{e,1} \\ \stvec{f}_{e,1} \\ \stvec{f}_{e,1}\end{array}\right)= \left(\begin{array}{ccc}\stvec{f}_e \\ \stvec{g}_e \\ \stvec{h}_e\end{array}\right),~~\ssvec{f}_{v} = \left(\begin{array}{ccc}\stvec{f}_{v,1} \\ \stvec{f}_{v,1} \\ \stvec{f}_{v,1}\end{array}\right)= \left(\begin{array}{ccc}\stvec{f}_v \\ \stvec{g}_v \\ \stvec{h}_v\end{array}\right).
\end{equation}
We can then define the product of two block vectors,
\begin{equation}
\ssvec{f}\cdot\ssvec{g} = \sum_{i=1}^{3}\stvec{f}_{i}^T\stvec{g}_{i},
\label{eq:notation:block+block}
\end{equation}
the product of a space vector with a block vector,

\begin{equation}
\svec{g}\cdot\ssvec{f} = \sum_{i=1}^{3}g_i \stvec{f}_{i},
\label{eq:notation:space+block}
\end{equation}
and the product of a space vector with a state vector,

\begin{equation}
\svec{g}\stvec{f} = \left(\begin{array}{c}g_1\stvec{f} \\ g_2\stvec{f} \\ g_3\stvec{f}\end{array}\right),
\label{eq:notation:space+state}
\end{equation}
which results in a block vector. The operators \eqref{eq:notation:space+block} and \eqref{eq:notation:space+state} allow us to define the divergence and gradient operators,
\begin{equation}
\svec{\nabla}_{x}\cdot\ssvec{f} = \sum_{i=1}^{3}\frac{\partial\stvec{f}_{i}}{\partial x_i},~~\svec{\nabla}_{x}\stvec{q} = \left(\begin{array}{c}\stvec{q}_{x} \\ \stvec{q}_{y} \\ 
\stvec{q}_z\end{array}\right),
\end{equation}
so that we can write \eqref{eq:governing:iNS-system-sum} in the compact form,
\begin{equation}
\stvec{q}_{t} + \svec{\nabla}_{x}\cdot\ssvec{f}_{e}(\stvec{q}) = \svec{\nabla}_{x}\cdot\ssvec{f}_{v}\left(\stvec{q},\svec{\nabla}_{x}\stvec{q}\right) + \stvec{s}\left(\stvec{q}\right).
\label{eq:governing:iNS-system}
\end{equation}

Finally, we use state matrices (i.e. $5\times 5$ matrices) which we write using an underline $\smat{B}$, and we can combine state matrices to construct a block matrix,
\begin{equation}
\mathcal B = \left(\begin{array}{ccc}\smat{B}_{11} & \smat{B}_{12} & \smat{B}_{13} \\
\smat{B}_{21} & \smat{B}_{22} & \smat{B}_{23}\\
\smat{B}_{31} & \smat{B}_{32} & \smat{B}_{33}\end{array}\right),
\label{eq:notation:block-matrix}
\end{equation}
which we can directly multiply to a block vector to obtain another block vector. For instance, if we want to perform a matrix multiplication in space (e.g. a rotation),
\begin{equation}
\svec{g} = \tens{M}\svec{f},
\end{equation}
for each of the variables in the state vector, we construct the block matrix version of $\tens{M}$,
\begin{equation}
\mathcal M = \left(\begin{array}{ccc}M_{11}\smat{I}_{5} & M_{12}\smat{I}_{5} & M_{13}\smat{I}_{5}\\
M_{21}\smat{I}_{5} & M_{22}\smat{I}_{5} & M_{23}\smat{I}_{5}\\
M_{31}\smat{I}_{5} & M_{32}\smat{I}_{5} & M_{33}\smat{I}_{5}\end{array}\right),
\label{eq:notation:block-matrix-from-space-matrix}
\end{equation}
so that we can compactly write
\begin{equation}
\ssvec{g} = \mathcal M \ssvec{f}.
\end{equation}
In \eqref{eq:notation:block-matrix-from-space-matrix}, $\smat{I}_{5}$ is the $5\times 5$ identity matrix.

\subsection{Entropy analysis}\label{subsec:GoverningEquations:Stability}

The stability analysis rests on the existence of a scalar mathematical entropy $\mathcal E\left(\stvec{q}\right)$ that satisfies three properties:

\begin{enumerate}
\item It is a concave function,

\begin{equation}
\mathcal E\left(\stvec{q}\right) \geqslant 0,~~\forall \stvec{q},
\end{equation}
that satisfies
\begin{equation}
\stvec{k}^T\frac{\partial^2 \mathcal E}{\partial \stvec{q}^2}\stvec{k} \geqslant 0,~~\forall \stvec{k}.
\end{equation}

\item It defines a set of entropy variables $\stvec{w}$,

\begin{equation}
\stvec{w} = \frac{\partial \mathcal E\left(\stvec{q}\right)}{\partial\stvec{q}},
\end{equation}
which contract the inviscid part of the original equation system   \eqref{eq:governing:iNS-system} as
\begin{equation}
\stvec{w}^T\left(\stvec{q}_{t} +\svec{\nabla}_{x}\cdot\ssvec{f}_{e}\right) = \mathcal E_{t} + \svec{\nabla}_{x}\cdot\svec{f}_e^{\mathcal E},
\label{eq:stability-cont:entropy-contraction}
\end{equation}
where $\svec{f}_{e}^{\mathcal E}=\left(f_e^{\mathcal E},g_e^{\mathcal E},h_e^{\mathcal E}\right)$ is the entropy flux.

\item The viscous fluxes are always dissipative when multiplied by the entropy 
variables $\stvec{w}$,
\begin{equation}
\stvec{w}^{T}\left(\svec{\nabla}_{x}\cdot\ssvec{f}_{v}\right)= \svec{\nabla}_{x}\cdot\left(\stvec{w}^T\ssvec{f}_{v}\right)-\left(\svec{\nabla}_{x}\stvec{w}^T\right)\cdot\ssvec{f}_{v},
\end{equation}
with 
\begin{equation}
\left(\svec{\nabla}_{x}\stvec{w}\right)^T\cdot\ssvec{f}_{v} \geqslant 0.
\label{eq:stability-cont:viscous-def-pos}
\end{equation}
\end{enumerate}
Entropy stability is guaranteed if we show that the mathematical entropy satisfies the evolution equation,
\begin{equation}
\mathcal E_{t} +\svec{\nabla}_{x}\cdot\left(\svec{f}_{e}^{\mathcal E}-\stvec{w}^T\ssvec{f}_{v}\right) =-\left(\svec{\nabla}_{x}\stvec{w}\right)^T\cdot\ssvec{f}_{v}\leqslant 0
\label{eq:stability-cont:balance-derivative}
\end{equation}
in the absence of lower order source terms (e.g. gravity) in the original equations.
When \eqref{eq:stability-cont:balance-derivative} is integrated over the domain $\Omega$,
\begin{equation}
\frac{\diff}{\diff t}\int_{\Omega}\mathcal E\diff\svec{x} + \int_{\Omega}\svec{\nabla}_{x}\cdot\left(\svec{f}_{e}^{\mathcal E}-\stvec{w}^T\ssvec{f}_{v}\right)\diff \svec{x} = \frac{\diff \overline{\mathcal E}(t)}{\diff t} + \int_{\partial \Omega}\left(\svec{f}_{e}^{\mathcal E}-\stvec{w}^T\ssvec{f}_{v}\right)\cdot\svec{n}\diff S \leqslant 0
\label{eq:stability-cont:balance}
\end{equation}
shows that the total (mathematical) entropy $\overline{\mathcal E} = \int_{\Omega}\mathcal E\diff\svec{x}$, decreases 
(ignoring the effect of boundary conditions and low order terms). In \eqref{eq:stability-cont:balance}, $\svec{n}$ is the outward 
pointing normal vector to $\partial \Omega$, and $\diff S$ is 
the differential surface element. Eq.  \eqref{eq:stability-cont:balance} implies that the original system   \eqref{eq:governing:iNS-system} is well--posed in 
the sense that the entropy $\mathcal E\left(\stvec{q}\right)$  is bounded by the initial state, shown by integrating \eqref{eq:stability-cont:balance} in time,
\begin{equation}
\bar{\mathcal E}(t) \leqslant \bar{\mathcal E}(t_0) - \int_{t_0}^{t}\left( \int_{\partial \Omega}\left(\svec{f}_{e}^{\mathcal E}-\stvec{w}^T\ssvec{f}_{v}\right)\cdot\svec{n}\diff S\right)\diff t.
\label{eq:stability-cont:balance-time}
\end{equation}

The incompressible NSE with artificial compressibility \eqref{eq:governing:iNS-system}, with inviscid \eqref{eq:governing:iNS-inviscid-flux} and viscous
\eqref{eq:governing:iNS-viscous-flux} fluxes admits the following mathematical entropy
\begin{equation}
\mathcal E\left(\stvec{q}\right) = \frac{1}{2}\rho\left(u^2 + v^2 + w^2\right) + \frac{1}{2}M_0^2p^2 = \mathcal K + \mathcal E_{AC},
\label{eq:stability-cont:entropy-def}
\end{equation}
which is the sum of the kinetic energy,
\begin{equation}
\mathcal K = \frac{1}{2}\rho\left(u^2 + v^2 + w^2\right),
\end{equation}
and an additional energy term due to artificial compressibility effects,
\begin{equation}
\mathcal E_{AC} = \frac{1}{2}M_0^2 p^2.
\end{equation}
Note that the artificial compressibility effects vanish as $M_0$ tends to zero (incompressible limit).

The entropy variables for the entropy \eqref{eq:stability-cont:entropy-def} are
\begin{equation}
\stvec{w}=\frac{\partial \mathcal E}{\partial\stvec{q}} = \left(-\frac{1}{2}v_{tot}^2, u, v, w, M_0^2 
p\right),
\end{equation}
where $v_{tot}= \left(u^2 + v^2 + w^2\right)^{1/2}$ is the total velocity. 

We now show that the entropy \eqref{eq:stability-cont:entropy-def} satisfies the three properties we enumerated previously.

\subsubsection{Convex function and positive semi--definite Hessian}
By construction, we see that the entropy \eqref{eq:stability-cont:entropy-def} is positive if the density remains positive, i.e.
\begin{equation}
\mathcal E\left(\stvec{q}\right)\geqslant 0~~\mathrm{if}~~\rho\left(\svec{x};t\right) \geqslant 0.
\end{equation}
Furthermore, the Hessian matrix of the entropy with respect to conservative variables is
\begin{equation}
\frac{\partial^2\mathcal E}{\partial\stvec{q}^2} = \frac{\partial \stvec{w}}{\partial\stvec{q}} = \left(\begin{array}{ccccc}v_{tot}^2/\rho & -u/\rho & -v/\rho & -w/\rho & 0 \\
-u/\rho & 1/\rho & 0 & 0 & 0\\
-v/\rho & 0 & 1/\rho & 0 & 0\\
-w/\rho & 0 & 0 & 1/\rho & 0\\
0 & 0 & 0 & 0 & M_0^2\end{array}\right),
\end{equation}
which is also positive semi--definite if the density $\rho(\svec{x};t)$ remains 
positive, since its eigenvalues are $\lambda_{1,2}=\frac{1}{\rho}$, $\lambda_{3}=\frac{1+v_{tot}^2}{\rho}$, $\lambda_{4}=M_0^2$, and 
$\lambda_{5}=0$.

\subsubsection{Contraction of the inviscid fluxes }
We now show that the entropy \eqref{eq:stability-cont:entropy-def} also satisfies the property \eqref{eq:stability-cont:entropy-contraction}.
% by finding a suitable entropy flux $\svec{f}_{e}^{\mathcal E}$. 
First, the time derivative,
\begin{equation}
\begin{split}
\stvec{w}^T\stvec{q}_t =& -\frac{1}{2}\left(u^2 + v^2 + w^2\right)\rho_t+u\left(\rho u\right)_t + v\left(\rho v\right)_t + w\left(\rho w\right)_t + M_0^2 p\left(p\right)_t\\
=&\left(-\frac{1}{2}u^2\rho_t + u\left(\rho u\right)_t\right)+\left(-\frac{1}{2}v^2\rho_t + v\left(\rho v\right)_t\right)\\
&+\left(-\frac{1}{2}w^2\rho_t  + w\left(\rho w\right)_t\right)+\left(\frac{1}{2}M_0^2p^2\right)_{t}.
\end{split}
\end{equation}
To obtain the kinetic energy, we manipulate each of the three velocity components $u_j$ contributions as
\begin{equation}
\begin{split}
-\frac{1}{2}u_j^2\rho_t + u_j\left(\rho u_j\right)_t &= -\frac{1}{2}u_j^2\rho_t + \frac{1}{2}u_j\left(\rho u_j\right)_t + \frac{1}{2}u_j\left(\rho u_j\right)_t \\
&= -\frac{1}{2}u_j^2\rho_t + \frac{1}{2}\rho_t u_j^2 +\frac{1}{2}\rho u_j u_{j,t}+\left(\frac{1}{2}\rho u_j^2\right)_t - \frac{1}{2}\rho u_j u_{j,t} \\
&= \left(\frac{1}{2}\rho u_j^2\right)_t.
\end{split}
\end{equation}
Hence, we have constructed the time contribution of the inviscid part contraction,
\begin{equation}
\stvec{w}^T\stvec{q}_t = \left(\frac{1}{2}\rho\left(u^2 + v^2 + w^2\right) + \frac{1}{2}M_0^2 p^2\right)_t = \mathcal E_{t}.
\label{eq:stability-cont:iNSAC-entropy-contraction}
\end{equation}

Next, the contraction of the inviscid fluxes rests on the existence of an entropy flux, $\svec{f}_{e}^{\mathcal E}=\left(f_e^{\mathcal E},g_e^{\mathcal E}, h_e^{\mathcal E}\right)$, that implies
\begin{equation}
\stvec{w}^T\left(\svec{\nabla}_{x}\cdot\ssvec{f}_{e}\right) = \stvec{w}^T\left(\stvec{f}_{e,x}+\stvec{g}_{e,y} + \stvec{h}_{e,z}\right) = \svec{\nabla}_{x}\cdot\svec{f}_{e}^{\mathcal E} = f_{e,x}^{\mathcal E}+g_{e,y}^{\mathcal E} + h_{e,z}^{\mathcal E}.
\end{equation}
Therefore, we must show that
\begin{equation}
\stvec{w}^{T}\stvec{f}_{e,x} = f_{e,x}^{\mathcal E},~~\stvec{w}^{T}\stvec{g}_{e,x} = g_{e,x}^{\mathcal E},~~\stvec{w}^{T}\stvec{h}_{e,x} = h_{e,x}^{\mathcal E}.
\end{equation}
 For the $x$--component,
\begin{equation}
\begin{split}
\stvec{w}^T \stvec{f}_{e,x} &= -\frac{1}{2}v_{tot}^2\left(\rho u\right)_x + u\left(\rho u^2 + p\right)_x + v\left(\rho u v\right)_x + w\left(\rho u w\right)_x + M_0^2 p \frac{1}{M_0^2}u_x \\
&= -\frac{1}{2}\left(u^2+v^2+w^2\right)\left(\rho u\right)_x + u\left(\rho u^2\right)_x + v\left(\rho u v\right)_x + w\left(\rho u w\right)_x +  p u_x + up_x,
\end{split}
\end{equation}
which consists of the sum of one term per velocity component $u_j$, 
\begin{equation}
\begin{split}
&-\frac{1}{2}u_j^2\left(\rho u\right)_x + u_j\left(\rho u u_j\right)_x = -\frac{1}{2}u_j^2\left(\rho u\right)_x + \frac{1}{2}u_j\left(\rho u u_j\right)_x+\frac{1}{2}u_j\left(\rho u u_j\right)_x  \\
&=-\frac{1}{2}u_j^2\left(\rho u\right)_x+ \frac{1}{2}u_j\left(\rho u u_j\right)_x+\left(\frac{1}{2}\rho u_j^2 u\right)_x - \frac{1}{2}\rho u u_j u_{j,x} = \left(\frac{1}{2}\rho u_j^2 
u\right)_x,
\end{split}
\end{equation}
plus the pressure work $\left(pu\right)_x$. Analogously, for $y$-- and $z$--components,
\begin{equation}
\begin{split}
\stvec{w}^T \stvec{f}_{e,x} &= \left(\left(\frac{1}{2}\rho\left(u^2+v^2+w^2\right) + p\right)u\right)_x=f_{e,x}^{\mathcal E}, \\
\stvec{w}^T \stvec{g}_{e,y} &= \left(\left(\frac{1}{2}\rho\left(u^2+v^2+w^2\right) + p\right)v\right)_y =g_{e,y}^{\mathcal E},\\
\stvec{w}^T \stvec{h}_{e,z} &= \left(\left(\frac{1}{2}\rho\left(u^2+v^2+w^2\right) + p\right)w\right)_z =h_{e,z}^{\mathcal E}.
\end{split}
\end{equation}
We therefore obtain the entropy flux,
\begin{equation}
\svec{f}_{e}^{\mathcal E} = \left(\frac{1}{2}\rho\left(u^2+v^2+w^2\right)+p\right)\svec{u},
\label{eq:stability-cont:entropy-flux}
\end{equation}
which drives the entropy evolution equation \eqref{eq:stability-cont:entropy-contraction}. We note that despite the fact that the entropy \eqref{eq:stability-cont:entropy-def} is particular to the incompressible Navier--Stokes equations with artificial compressibility, the entropy flux obtained is the same as in the incompressible NSE without artificial compressibility.

\subsubsection{Positive semi--definiteness of the viscous fluxes }

Finally,  we show that viscous fluxes and entropy variables satisfy \eqref{eq:stability-cont:viscous-def-pos}. To do so, we re--write the viscous fluxes as a function of the entropy vector $\stvec{w}$ (instead of the conservative vector $\stvec{q}$), which in general (both for compressible and incompressible formulations) can be linearly spanned in the gradient of the entropy variables,
\begin{equation}
\stvec{f}^{i}_{v}\left(\stvec{w},\svec{\nabla}_{x}\stvec{w}\right) = \sum_{j=1}^{3}\smat{B}^{\mathcal E}_{ij}\left(\stvec{w}\right)\frac{\partial\stvec{w}}{\partial x_j},
\end{equation}
or, using the block matrix representation introduced in \eqref{eq:notation:block-matrix},
\begin{equation}
\ssvec{f}_{v}\left(\stvec{w},\svec{\nabla}_{x}\stvec{w}\right) = \mathcal B^{\mathcal E}(\stvec{w})\svec{\nabla}_{x}\stvec{w}.
\label{eq:stability-cont:viscous-flux-Y=AX}
\end{equation}
In general, the coefficients $\smat{B}^{\mathcal E}_{ij}(\stvec{w})$ are non--linear and depend on the entropy variables. In particular to the incompressible NSE, the matrices $\smat{B}_{ij}^{\mathcal E}$ are constant (i.e. they do not depend on $\stvec{w}$),
\begin{equation}
\begin{split}
\smat{B}^{\mathcal E}_{11} = \frac{1}{\Rey}\left(\begin{array}{ccccc} 0 & 0 & 0 & 0 & 0 \\ 0 & 2 & 0 & 0 & 0 \\ 0 & 0 & 1 & 0 & 0 \\ 0 & 0 & 0 & 1 & 0 \\ 0 & 0 & 0 & 0 & 0 \end{array}\right),~~
\smat{B}^{\mathcal E}_{12} = \frac{1}{\Rey}\left(\begin{array}{ccccc} 0 & 0 & 0 & 0 & 0 \\ 0 & 0 & 0 & 0 & 0 \\ 0 & 1 & 0 & 0 & 0 \\ 0 & 0 & 0 & 0 & 0 \\ 0 & 0 & 0 & 0 & 0 \end{array}\right), ~~
\smat{B}^{\mathcal E}_{13} = \frac{1}{\Rey}\left(\begin{array}{ccccc} 0 & 0 & 0 & 0 & 0 \\ 0 & 0 & 0 & 0 & 0 \\ 0 & 0 & 0 & 0 & 0 \\ 0 & 1 & 0 & 0 & 0 \\ 0 & 0 & 0 & 0 & 0 \end{array}\right), \\
\smat{B}^{\mathcal E}_{21} = \frac{1}{\Rey}\left(\begin{array}{ccccc} 0 & 0 & 0 & 0 & 0 \\ 0 & 0 & 1 & 0 & 0 \\ 0 & 0 & 0 & 0 & 0 \\ 0 & 0 & 0 & 0 & 0 \\ 0 & 0 & 0 & 0 & 0 \end{array}\right),~~
\smat{B}^{\mathcal E}_{22} = \frac{1}{\Rey}\left(\begin{array}{ccccc} 0 & 0 & 0 & 0 & 0 \\ 0 & 1 & 0 & 0 & 0 \\ 0 & 0 & 2 & 0 & 0 \\ 0 & 0 & 0 & 1 & 0 \\ 0 & 0 & 0 & 0 & 0 \end{array}\right), ~~
\smat{B}^{\mathcal E}_{23} = \frac{1}{\Rey}\left(\begin{array}{ccccc} 0 & 0 & 0 & 0 & 0 \\ 0 & 0 & 0 & 0 & 0 \\ 0 & 0 & 0 & 0 & 0 \\ 0 & 0 & 1 & 0 & 0 \\ 0 & 0 & 0 & 0 & 0 \end{array}\right),\\
\smat{B}^{\mathcal E}_{31} = \frac{1}{\Rey}\left(\begin{array}{ccccc} 0 & 0 & 0 & 0 & 0 \\ 0 & 0 & 0 & 1 & 0 \\ 0 & 0 & 0 & 0 & 0 \\ 0 & 0 & 0 & 0 & 0 \\ 0 & 0 & 0 & 0 & 0 \end{array}\right),~~
\smat{B}^{\mathcal E}_{32} = \frac{1}{\Rey}\left(\begin{array}{ccccc} 0 & 0 & 0 & 0 & 0 \\ 0 & 0 & 0 & 0 & 0 \\ 0 & 0 & 0 & 1 & 0 \\ 0 & 0 & 0 & 0 & 0 \\ 0 & 0 & 0 & 0 & 0 \end{array}\right), ~~
\smat{B}^{\mathcal E}_{33} = \frac{1}{\Rey}\left(\begin{array}{ccccc} 0 & 0 & 0 & 0 & 0 \\ 0 & 1 & 0 & 0 & 0 \\ 0 & 0 & 1 & 0 & 0 \\ 0 & 0 & 0 & 2 & 0 \\ 0 & 0 & 0 & 0 & 0 \end{array}\right).
\end{split}
\label{eq:stability-cont:viscous-matrices}
\end{equation}

For entropy stability, then, the set of matrices $\smat{B}^{\mathcal E}_{ij}$ must satisfy:
\begin{enumerate}
\item Symmetry,
\begin{equation}
\smat{B}^{\mathcal E}_{ij} = \left(\smat{B}^{\mathcal E}_{ji}\right)^T,
\label{eq:entropy-stability:viscous-matrices-symmetry}
\end{equation}
and
\item Positive semi--definiteness with respect to the gradient of the entropy variables, 

\begin{equation}
\left(\svec{\nabla}_{x}\stvec{w}\right)^T\cdot\ssvec{f}_{v} = \sum_{i,j=1}^{3}\left(\frac{\partial \stvec{w}}{\partial x_i}\right)^T\smat{B}^{\mathcal E}_{ij}\left(\frac{\partial\stvec{w}}{\partial x_j}\right) \geqslant  0.
\label{eq:entropy-stability:viscous-matrices-pos-def}
\end{equation}
\end{enumerate}
For the matrices given in \eqref{eq:stability-cont:viscous-matrices}, both the first property \eqref{eq:entropy-stability:viscous-matrices-symmetry} and second property \eqref{eq:entropy-stability:viscous-matrices-pos-def} are satisfied. The first is immediate, and the second is,

\begin{equation}
\begin{split}
\sum_{i,j=1}^{3}\left(\frac{\partial \stvec{w}}{\partial x_i}\right)^T\smat{B}^{\mathcal E}_{ij}\left(\frac{\partial\stvec{w}}{\partial x_j}\right) &= \frac{1}{\Rey}\left(2u_x^2 + v_x^2 + w_x^2 +u_y^2+2v_y^2+w_y^2+u_z^2+v_z^2+2w_z^2+ 2u_yv_x + 2u_zw_x+2v_zw_y\right)\\
&=\frac{1}{\Rey}\left(2u_x^2 +2v_y^2+2w_z^2+\left(u_y+v_x\right)^2+ \left(u_z+w_x\right)^2+\left(v_z+w_y\right)^2\right)\geqslant 0,
\end{split}
\label{eq:stability-cont:viscous-work-expression}
\end{equation}
thus proving \eqref{eq:stability-cont:viscous-def-pos}. The last identity \eqref{eq:stability-cont:viscous-work-expression} is commonly represented as

\begin{equation}
\left(\svec{\nabla}_{x}\stvec{w}\right)^T\cdot\ssvec{f}_{v}= \frac{2}{\Rey}\tens{S}:\tens{S} \geqslant 0,
\label{eq:stability-cont:viscous-strain-tensor}
\end{equation}
where $\tens{S}$ is the strain tensor \eqref{eq:governing:strain-tensor}.

Hence, we can write the entropy equation \eqref{eq:stability-cont:balance-derivative} in its particular version for the incompressible NSE with artificial compressibility, but ignoring low order terms, as
\begin{equation}
\begin{split}
\left(\frac{1}{2}\rho v_{tot}^2 + \frac{1}{2}M_0^2p^2\right)_{t} + \svec{\nabla}_{x}\cdot\left(\left(\frac{1}{2}\rho v_{tot}^2+p\right)\svec{u}-\svec{u}\cdot\tens{\tau}\right) =-\frac{2}{\Rey}\tens{S}:\tens{S}\leqslant 
0,
\end{split}
\label{eq:stability-cont:balance-final}
\end{equation}
or, integrated in the domain $\Omega$,
\begin{equation}
\begin{split}
\frac{\diff}{\diff t}\int_{\Omega}\left(\frac{1}{2}\rho v_{tot}^2 + \frac{1}{2}M_0^2p^2\right)\diff\svec{x} + \int_{\partial\Omega}\left(\left(\frac{1}{2}\rho v_{tot}^2+p\right)\svec{u}-\svec{u}\cdot\tens{\tau}\right)\cdot\svec{n}\diff S =-\int_{\Omega}\frac{2}{\Rey}\tens{S}:\tens{S}\diff\svec{x}\leqslant 
0.
\end{split}
\label{eq:stability-cont:balance-final-omega}
\end{equation}
To the best of our knowledge, this is the first time that this entropy analysis 
was performed for this set of equations.

The boundary integral is studied for free-- and no--slip wall boundary 
conditions. In this continuous analysis is zero in both of them, since in free--slip 
walls we set $\svec{u}\cdot\svec{n}=0$ and $\svec{\tau}\cdot\svec{n}=0$, and in 
no--slip walls is $\svec{u}=0$. Therefore,

\begin{equation}
\begin{split}
\frac{\diff}{\diff t}\int_{\Omega}\left(\frac{1}{2}\rho v_{tot}^2 + \frac{1}{2}M_0^2p^2\right)\diff\svec{x}  =-\int_{\Omega}\frac{2}{\Rey}\tens{S}:\tens{S}\diff\svec{x}\leqslant 
0,
\end{split}
\label{eq:stability-cont:balance-final-omega-bc}
\end{equation}
so that the entropy is conserved when the flow is inviscid, and dissipated 
otherwise.

\section{Discontinuous Galerkin approximation}\label{sec:DGSEM}

In the this section, we construct an entropy--stable DG scheme that satisfies a semi--discrete version of the bound \eqref{eq:stability-cont:balance-final} using split--forms and the SBP--SAT property.

We detail the construction of the nodal Discontinuous Galerkin Spectral Element Method (DGSEM).
 From all the variants, we restrict ourselves to the tensor product DGSEM with 
 Gauss--Lobatto (GL) points, since it satisfies the Summation--By--Parts
  Simultaneous--Approximation--Term (SBP--SAT) property \cite{2014:Carpenter}. 
  This property is used to prove the scheme's stability without relying on exact integration.

The reference domain $\Omega$ is tessellated with non--overlapping hexahedral elements, $e$, which are geometrically transformed to a reference element $E=[-1,1]^3$. This transformation is performed using a (polynomial) transfinite mapping $\svec{X}$ that relates physical ($\svec{x}=\left(x^1,x^2,x^3\right)=\left(x,y,z\right)=x\hat{x}+y\hat{y}+z\hat{z}$), and local reference ($\svec{\xi}=\left(\xi^1,\xi^2,\xi^3\right)=\left(\xi,\eta,\zeta\right)=\xi\hat{\xi}+\eta\hat{\eta}+\zeta\hat{\zeta}$) coordinates through
\begin{equation}
\svec{x}=\svec{X}\left(\svec{\xi}\right)=\svec{X}\left(\xi,\eta,\zeta\right).
\label{eq:Mapping}
\end{equation}
The space vectors $\hat{x}_{i}$ and $\hat{\xi}_{i}$ are unit vectors in the three Cartesian directions of physical and reference coordinates, respectively. 

The transformation \eqref{eq:Mapping} defines three covariant basis vectors,
\begin{equation}
\svec{a}_{i} = \frac{\partial\svec{X}}{\partial \xi^{i}},~~i=1,2,3,
\end{equation}
and three contravariant basis vectors,
\begin{equation}
\svec{a}^{i} = \svec{\nabla}_{x}\xi^{i} = \frac{1}{J}\left(\svec{a}_{j}\times\svec{a}_k\right),~~(i,j,k)\text{ cyclic},
\end{equation}
where 
\begin{equation}
J = \svec{a}_i\cdot\left(\svec{a}_{j}\times\svec{a}_{k}\right)
\end{equation}
is the Jacobian of the mapping $\svec{X}$. The contravariant coordinate vectors satisfy the metric identities \cite{2006:Kopriva},
\begin{equation}
\sum_{i=1}^{3}\frac{\partial \left(Ja_n^{i}\right)}{\partial\xi^{i}} = 0,~~n=1,2,3,
\label{eq:metrics:metric-id-cont}
\end{equation}
where $a_{n}^{i}$ is the $n$--th Cartesian component of the contravariant vector $\svec{a}^{i}$.

We use the volume weighted contravariant basis $J\svec{a}^{i}$ to transform differential operators from physical ($\svec{\nabla}_{x}$) to reference ($\svec{\nabla}_{\xi}$) space. The divergence of a vector is
\begin{equation}
\snabla_x\cdot\svec{f} = \frac{1}{J}\snabla_{\xi}\cdot\left(\tens{M}^T\svec{f}\right),
\label{eq:metrics:div-transf-1eqn}
\end{equation}
where $\tens{M}=\left(J\svec{a}^{\xi},J\svec{a}^{\eta},J\svec{a}^{\zeta}\right)$ is the metrics matrix. We use \eqref{eq:notation:block-matrix-from-space-matrix} to write the divergence of an entire block vector compactly. Thus, we define the metrics block matrix $\mathcal M$,
\begin{equation}
\mathcal M = \left(\begin{array}{ccc}Ja^{1}_{1}\smat{I}_{5} & Ja^{2}_{1}\smat{I}_{5} & Ja^{3}_{1}\smat{I}_{5}\\
Ja^{1}_{2}\smat{I}_{5} & Ja^{2}_{2}\smat{I}_{5} & Ja^{3}_{2}\smat{I}_{5}\\
Ja^{1}_{3}\smat{I}_{5} & Ja^{2}_{3}\smat{I}_{5} & Ja^{3}_{3}\smat{I}_{5}\end{array}\right),
\label{eq:metrics:metrics-matrix}
\end{equation}
which allows us to write \eqref{eq:metrics:div-transf-1eqn} for all the state variables,
\begin{equation}
\snabla_x\cdot\ssvec{f} = \frac{1}{J}\snabla_{\xi}\cdot\left(\mathcal M^T\ssvec{f}\right)=\frac{1}{J}\snabla_{\xi}\cdot\cssvec{f}  ,
\label{eq:metrics:div-transf}
\end{equation}
with $\cssvec{f}$ being the block vector of the contravariant fluxes,
\begin{equation}
\cssvec{f}=\mathcal M^{T}\ssvec{f},~~\tilde{\stvec{f}}^{i} = J\svec{a}^{i}\cdot\ssvec{f}.
\label{eq:metrics:contravariant-flux}
\end{equation}
The gradient of a scalar is
\begin{equation}
\snabla_{x} u = \frac{1}{J}\tens{M}\snabla_{\xi}u.
\label{eq:metrics:gradient-transf-1eqn}
\end{equation}
which we can also extend to all state variables using \eqref{eq:metrics:metrics-matrix},
\begin{equation}
\snabla_x \stvec{u} = \frac{1}{J}{\mathcal M}\snabla_{\xi}\stvec{u}.
\label{eq:metrics:gradient-transf}
\end{equation}

To transform the incompressible NSE \eqref{eq:governing:iNS-system} into reference space, we first write them as a first order system, defining the auxiliary variable $\ssvec{g}=\snabla_x\stvec{w}$ so that
\begin{equation}
\begin{split}
\stvec{q}_{t} +\snabla_{x}\cdot\ssvec{f}_{e}(\stvec{q}) &= \snabla_{x}\cdot\ssvec{f}_{v}\left(\ssvec{g}\right) + \stvec{s}(\stvec{q}),\\
\ssvec{g}&= \snabla_x\stvec{w}.
\end{split}
\end{equation}
Recall that the incompressible NSE viscous fluxes depend only on the gradient of the entropy variables, $\ssvec{g}$. 

Next, we transform the operators using \eqref{eq:metrics:div-transf} and \eqref{eq:metrics:gradient-transf},
\begin{subequations}\label{eq:dg:system-1st-order-computational}
\begin{align}
J\stvec{q}_{t} +\snabla_{\xi}\cdot\cssvec{f}_{e}(\stvec{q}) &= \snabla_{\xi}\cdot\cssvec{f}_{v}\left(\ssvec{g}\right) + J\stvec{s}(\stvec{q}),
\label{eq:dg:system-1st-order-computational-q}\\
J\ssvec{g}&= {\mathcal M}\snabla_{\xi}\stvec{w}
\label{eq:dg:system-1st-order-computational-g}
\end{align}
\end{subequations}
to get the final form of the equations to be approximated.

The DG approximation is obtained from weak forms of the equations \eqref{eq:dg:system-1st-order-computational}. We first define the inner product in the reference element, $E$, for state and block vectors
\begin{equation}
\left\langle \stvec{f},\stvec{g}\right\rangle_{E} = \int_{E}\stvec{f}^{T}\stvec{g}\diff E,~~\left\langle\ssvec{f},\ssvec{g}\right\rangle_{E} = \int_{E}\ssvec{f}\cdot\ssvec{g}\diff E.
\label{eq:dg:weak-forms}
\end{equation}
We construct two weak forms by multiplying \eqref{eq:dg:system-1st-order-computational-q} by a test function $\stvecg{\phi}$, and \eqref{eq:dg:system-1st-order-computational-g} by a block vector test function $\ssvecg{\psi}$, then integrating over the reference element $E$, and finally integrating by parts to get
\begin{equation}
\begin{split}
\left\langle J\stvec{q}_{t},\stvecg{\phi}\right\rangle_{E} &+\int_{\partial E}\stvecg{\phi}^{T}\left(\cssvec{f}_{e}-\cssvec{f}_{v}\right)\cdot\hat{n}\diff S_{\xi}-\left\langle\cssvec{f}_{e},\snabla_{\xi}\stvecg{\phi}\right\rangle_{E} =- \left\langle\cssvec{f}_{v},\snabla_{\xi}\stvecg{\phi}\right\rangle_{E} + \left\langle J\stvec{s},\stvecg{\phi}\right\rangle_{E},\\
\left\langle J\ssvec{g},\ssvecg{\psi}\right\rangle_{E}&= \int_{\partial E}\stvec{w}^{T}\left(\cssvecg{\psi}\cdot\hat{n}\right)\diff S_{\xi} - \left\langle \stvec{w},\snabla_{\xi}\cdot\cssvecg{\psi}\right\rangle_{E}.
\end{split}
\label{eq:dg:weak-form-1}
\end{equation}
The quantities $\hat{n}$ and $\diff S_{\xi}$ are the unit outward pointing normal and surface differential at the faces of $E$, respectively. The contravariant test function $\cssvecg{\psi}$ follows the definition \eqref{eq:metrics:contravariant-flux}. Finally, surface integrals extend to all six faces of an element,
\begin{equation}
\int_{\partial E} \svec{\tilde{f}}\cdot\hat{n}\diff S_{\xi} = \int_{[-1,1]^{2}}\!\!\!\!\!\!\tilde{f}^{1}\diff \eta \diff \zeta\biggr|_{\xi=-1}^{\xi=1}\!\!+\int_{[-1,1]^{2}}\!\!\!\!\!\!\tilde{f}^{2}\diff \xi \diff \zeta\biggr|_{\eta=-1}^{\eta=1}\!\!+\int_{[-1,1]^{2}}\!\!\!\!\!\!\tilde{f}^{3}\diff \xi\diff \eta \biggr|_{\zeta=-1}^{\zeta=1}.
\label{eq:dg:surface-integral-def-cont}
\end{equation}

We can write surface integrals in either physical or reference space. The relation to the physical surface integration variables is
\begin{equation}
\diff S^i = \left| J\stvec{a}^i\right|\diff\xi^{j}\diff\xi^{k} = \mathcal J_f^i \diff S_\xi^i,
\end{equation}
where we have defined the face Jacobian $\mathcal J_{f}^i = \left|\mathcal J \stvec{a}^i\right|$. We can relate the surface flux in both reference element, $\tilde{\boldsymbol{f}}\cdot\hat{\boldsymbol{n}}$, and physical, $\stvec{f}\cdot\stvec{n}$, variables through
\begin{equation}
\tilde{\stvec{f}}\cdot\hat{\stvec{n}}^{i}\diff S_\xi = \left(\boldsymbol{\mathcal M}^T\stvec{f}\right)\cdot\hat{\stvec{n}}^{i}\diff S_\xi = \stvec{f}\cdot\left(\boldsymbol{\mathcal M}\hat{\stvec{n}}^{i}\right)\diff S_\xi = \stvec{f}\cdot\stvec{n}\left|J\stvec{a}^i\right|\diff S_\xi = \stvec{f}\cdot\stvec{n}^{i}\diff S.
\label{eq:dg:local-physical-fluxes-relationship}
\end{equation}
Therefore, the surface integrals can be represented both in physical and reference 
spaces,
\begin{equation}
  \int_{\partial E}\tilde{\stvec{f}}\cdot\hat{\stvec{n}}\diff S_\xi = \int_{\partial 
  e}\stvec{f}\cdot\stvec{n}\diff S,
  \label{eq:dg:surface-integrals-relation}
\end{equation}
and we will use one or the other depending on whether we are studying an 
isolated element (reference space) or the whole combination of elements in 
the mesh (physical space).

\subsection{Polynomial approximation in the DGSEM}\label{subsec:DGSEM:Approximation}

We now construct the discrete version of \eqref{eq:dg:weak-form-1}. The approximation of the state vector inside each element $E$ is an order $N$ polynomial,
\begin{equation}
\stvec{q}\approx \mathbb I^{N}(\stvec{q})=\stvec{Q}\left(\svec{\xi}\right) = \sum_{i,j,k=0}^{N}\stvec{Q}_{ijk}(t)l_i(\xi)l_j(\eta)l_k(\zeta)\in\mathbb P^{N}.
\label{eq:dg:interpolation}
\end{equation}
where $\mathbb P^{N}$ is the space of polynomials of degree less than or equal to $N$ on $[-1,1]^{3}$ and $\mathbb I^{N}$ is the polynomial interpolation operator. The state values $\stvec{Q}_{ijk}(t)=\stvec{Q}(\xi_i,\eta_j,\zeta_k;t)$ are the nodal degrees of freedom (time dependent coefficients) at the tensor product of each of the Gauss--Lobatto points $\{\xi_i\}_{i=0}^{N}$, where we write Lagrange polynomials $l_i(\xi)$,
\begin{equation}
l_i(\xi) = \prod_{\substack{j=0 \\ j\neq i}}^{N}\frac{\xi-\xi_j}{\xi_i-\xi_j}.
\end{equation}

The geometry and metric terms are also approximated with order N polynomials. The transfinite mapping is approximated using $\mathcal X=\mathbb I^{N}\left(X\right)$, but special attention must be paid to its derivatives (i.e. the contravariant basis) since $\mathcal J\svec{a}^{i}\neq\mathbb I^{N}\left(\svec{a}_{j}\times\svec{a}_k\right)$. For the discrete version of the metric identities \eqref{eq:metrics:metric-id-cont} to hold,
\begin{equation}
\sum_{i=1}^{3}\frac{\partial \mathbb I^{N}\left(\mathcal J a_n^i\right)}{\partial \xi^i} = 0,~~n=1,2,3,
\end{equation}
we compute metric terms using the curl form \cite{2006:Kopriva},
\begin{equation}
\mathcal J a_{n}^{i} = -\hat{x}_i \cdot\snabla_{\xi}\times\left(\mathbb I^{N}\left(\mathcal X_{l}\snabla_{\xi}\mathcal X_{m}\right)\right),~~~i=1,2,3,~~n=1,2,3,~~(n,m,l)\text{ cyclic}.
\label{eq:dg:contravariant-basis}
\end{equation}
If we compute $\mathcal J \svec{a}^{i}$ using \eqref{eq:dg:contravariant-basis}, we ensure discrete free--stream preservation, which is crucial to avoid grid induced motions and importantly also to guarantee entropy stability.

Next, we approximate integrals that arise in the weak formulation using Gauss quadratures. Let $\{w_i\}_{i=0}^{N}$ be the quadrature weights associated to Gauss--Lobatto nodes $\{\xi_{i}\}_{i=0}^{N}$. Then, in one dimension,
\begin{equation}
\int_{-1}^{1}f(\xi)\diff \xi \approx \int_{N}f(\xi)\diff \xi \equiv \sum_{m=0}^{N}w_m f(\xi_{m})=\sum_{m=0}^{N}w_mF_{m}.
\end{equation}
For Gauss--Lobatto points, the approximation is exact if $f\left(\xi\right)\in\mathbb P^{2N-1}$. The extension to three dimensions has three nested quadratures in each of the three reference space dimensions,
\begin{equation}
\left\langle f,g\right\rangle_{E} \approx\left\langle f,g\right\rangle_{E,N} = 
\sum_{m,n,l=0}^{N}w_{mnl}F_{mnl}G_{mnl},~~w_{mnl}=w_mw_nw_l,
\label{eq:dg:discrete-inner}
\end{equation}
with a similar definition for block vectors.
The inner product is computed exactly if $f(\xi)g(\xi)\in\mathbb P^{2N-1}$. The approximation of surface integrals is performed similarly, replacing exact integrals by Gauss quadratures in \eqref{eq:dg:surface-integral-def-cont},
\begin{equation}
\begin{split}
\int_{\partial E} \svec{\tilde{f}}\cdot\hat{n}\diff S_{\xi} &\approx\int_{\partial E,N} \svec{\tilde{f}}\cdot\hat{n}\diff S_{\xi} = \int_{N}\tilde{f}^{1}\diff \eta \diff \zeta\biggr|_{\xi=-1}^{\xi=1}\!\!+\int_{N}\tilde{f}^{2}\diff \xi \diff \zeta\biggr|_{\eta=-1}^{\eta=1}\!\!+\int_{N}\tilde{f}^{3}\diff \xi\diff \eta 
\biggr|_{\zeta=-1}^{\zeta=1}\\
&\equiv\sum_{j,k=0}^{N}w_{jk}\left(\tilde{{F}}^{1}_{Njk}-\tilde{{F}}^{1}_{0jk}\right)+\sum_{i,k=0}^{N}w_{ik}\left(\tilde{{F}}^{2}_{iNk}-\tilde{{F}}^{2}_{i0k}\right)+\sum_{i,j=0}^{N}w_{ij}\left(\tilde{{F}}^{3}_{ijN}-\tilde{{F}}^{3}_{ij0}\right).
\end{split}
\label{eq:dg:surface-integral-def-disc}
\end{equation}

Gauss--Lobatto points are used to construct entropy--stable schemes using split--forms. Since boundary points are included, there is no need to perform an interpolation of the volume polynomials in \eqref{eq:dg:surface-integral-def-disc}  to the boundaries, which is known as the Simultaneous--Approximation--Term (SAT) property. The SAT property, alongside the exactness of the quadrature yields the Summation--By--Parts (SBP)%
\begin{equation}
\left\langle \svec{\nabla}_{\xi}\cdot\cssvec{F},\stvec{V}\right\rangle_{E,N} = \int_{\partial E,N}\left(\cssvec{F}\cdot\hat{n}\right)\stvec{V}\diff S_{\xi} - \left\langle \cssvec{F},\snabla_\xi \stvec{V}\right\rangle_{E,N}
\label{eq:dg:discreteGaussLaw}
\end{equation}
property \cite{2014:Carpenter}.

With the polynomial framework \eqref{eq:dg:interpolation}, \eqref{eq:dg:discrete-inner}, and \eqref{eq:dg:surface-integral-def-disc} in the continuous weak forms \eqref{eq:dg:weak-form-1}, we get the \textit{standard} DG version of the incompressible NSE,
\begin{equation}
\begin{split}
\left\langle \mathcal J\stvec{Q}_{t},\stvecg{\phi}\right\rangle_{E,N} &+\int_{\partial E,N}\stvecg{\phi}^{T}\left(\cssvec{F}_{e}-\cssvec{F}_{v}\right)\cdot\hat{n}\diff S_{\xi}-\left\langle\cssvec{F}_{e},\snabla_{\xi}\stvecg{\phi}\right\rangle_{E,N} =- \left\langle\cssvec{F}_{v},\snabla_{\xi}\stvecg{\phi}\right\rangle_{E,N} + \left\langle \mathcal J \stvec{S},\stvecg{\phi}\right\rangle_{E,N},\\
\left\langle \mathcal J\ssvec{G},\ssvecg{\psi}\right\rangle_{E,N}&= \int_{\partial E,N}\stvec{W}^{T}\left(\cssvecg{\psi}\cdot\hat{n}\right)\diff S_{\xi} - \left\langle \stvec{W},\snabla_{\xi}\cdot\cssvecg{\psi}\right\rangle_{E,N},
\end{split}
\label{eq:dg:discrete-weak-form}
\end{equation}
where the test functions $\stvecg{\phi}$ and $\ssvecg{\psi}$ are now restricted to polynomial spaces in $\mathbb P^{N}$. 

Inter--element coupling is enforced using \textit{numerical fluxes} at the element boundaries, 
\begin{equation}
 \ssvec{F}_{e}\approx{\ssvec{F}}_{e}^{\star}\left(\stvec{Q}_L,\stvec{Q}_R\right),~~ \ssvec{F}_{v}\approx{\ssvec{F}}_{v}^{\star}\left(\ssvec{G}_L,\ssvec{G}_R\right),~~ \stvec{W}\approx \stvec{W}^{\star}\left(\stvec{Q}_L,\stvec{Q}_R\right),
 \label{eq:dg:numerical-fluxes}
\end{equation}
where $\stvec{Q}_L$, $\ssvec{G}_L$, $\stvec{Q}_R$, and $\ssvec{G}_R$ are state vectors and gradients at the left and right sides of the face. Details on the {numerical flux} functions are given below in Sec. \ref{subsec:NumericalFluxes}. With \eqref{eq:dg:numerical-fluxes}, \eqref{eq:dg:discrete-weak-form} becomes,
\begin{subequations}\label{eq:dg:discrete-weak-form-numflux}
\begin{align}
\left\langle \mathcal J\stvec{Q}_{t},\stvecg{\phi}\right\rangle_{E,N} &+\int_{\partial E,N}\stvecg{\phi}^{T}\left(\tilde{\stvec{F}}_{e}^{\star}-\tilde{\stvec{F}}_{v}^{\star}\right)\cdot\hat{n}\diff S_{\xi}-\left\langle\cssvec{F}_{e},\snabla_{\xi}\stvecg{\phi}\right\rangle_{E,N} =- \left\langle\cssvec{F}_{v},\snabla_{\xi}\stvecg{\phi}\right\rangle_{E,N} + \left\langle \mathcal J \stvec{S},\stvecg{\phi}\right\rangle_{E,N},\label{eq:dg:discrete-weak-form-numflux:q}\\
\left\langle \mathcal J\ssvec{G},\ssvecg{\psi}\right\rangle_{E,N}&= \int_{\partial E,N}\stvec{W}^{\star,T}\left(\cssvecg{\psi}\cdot\hat{n}\right)\diff S_{\xi} - \left\langle \stvec{W},\snabla_{\xi}\cdot\cssvecg{\psi}\right\rangle_{E,N}.\label{eq:dg:discrete-weak-form-numflux:g}
\end{align}
\end{subequations}

To construct entropy stable schemes using the SBP--SAT property and two--point volume fluxes following \cite{2016:gassner,2017:Gassner}, we apply the discrete Gauss Law \eqref{eq:dg:discreteGaussLaw} to transform \eqref{eq:dg:discrete-weak-form-numflux:q} into a strong DG form,
\begin{equation}
	\begin{split}
	\left\langle \mathcal J\stvec{Q}_{t},\stvecg{\phi}\right\rangle_{E,N} &+\int_{\partial E,N}\stvecg{\phi}^{T}\left({\cssvec{F}}_{e}^{\star}-\cssvec{F}_{e}\right)\cdot\hat{n}\diff S_{\xi}+\left\langle\snabla_{\xi}\cdot\cssvec{F}_{e},\stvecg{\phi}\right\rangle_{E,N} \\
	&= \int_{\partial E,N}\stvecg{\phi}^{T}\left({\cssvec{F}}_{v}^{\star}-\cssvec{F}_{v}\right)\cdot\hat{n}\diff S_{\xi}+\left\langle\snabla_{\xi}\cdot\cssvec{F}_{v},\stvecg{\phi}\right\rangle_{E,N} + \left\langle \mathcal J \stvec{S},\stvecg{\phi}\right\rangle_{E,N},\\
	\left\langle \mathcal J\ssvec{G},\ssvecg{\psi}\right\rangle_{E,N}&= \int_{\partial E,N}\stvec{W}^{\star,T}\left(\cssvecg{\psi}\cdot\hat{n}\right)\diff S_{\xi} - \left\langle \stvec{W},\snabla_{\xi}\cdot\cssvecg{\psi}\right\rangle_{E,N}.
	\end{split}
	\label{eq:dg:discrete-strong-form}
\end{equation}
To compute the divergence of the viscous flux, we use the standard DG approach (i.e. direct differentiation of the polynomials),
\begin{equation}
\snabla_{\xi}\cdot\cssvec{F}_{v} = \sum_{m=0}^{N}\left(\tilde{\stvec{F}}^{1}_{mjk}l_m'(\xi)+\tilde{\stvec{F}}^{2}_{imk}l_m'(\eta)+\tilde{\stvec{F}}^{3}_{ijm}l_m'(\zeta)\right).
\label{eq:dg:standard-divergence}
\end{equation}
However to get a stable approximation, we compute the divergence of a two--point inviscid flux \cite{2003:Tadmor,2013:Fisher}
\begin{equation}
\begin{split}
\snabla_{\xi}\cdot\cssvec{F}_{e}\approx \mathbb {D}\left(\cssvec{F}_{e}\right)^{\#}(\xi,\eta,\zeta)\equiv 2\sum_{m=0}^{N}&\phantom{{}+{}}l'_m(\xi)\tilde{\stvec{F}}^{1,\#}_{e}\left(\xi,\eta,\zeta;\xi_m,\eta,\zeta\right)\\
&+l'_m(\eta)\tilde{\stvec{F}}^{2,\#}_{e}\left(\xi,\eta,\zeta;\xi,\eta_m,\zeta\right)\\
&+l'_m(\zeta)\tilde{\stvec{F}}^{3,\#}_{e}\left(\xi,\eta,\zeta;\xi,\eta,\zeta_m\right).
\end{split}
\end{equation}
The contravariant numerical volume flux is constructed from a flux function 
$\ssvec{F}^{\#}_{e}(\cdot,\cdot)$ called two--point flux that is a function of two states
 \cite{2003:Tadmor,2016:gassner}
\begin{equation}
\tilde{\stvec{F}}^{l,\#}_{e}\left(\xi,\eta,\zeta;\alpha,\beta,\gamma\right) \equiv \ssvec{F}^{\#}_{e}\left(\stvec Q(\xi,\eta,\zeta),\stvec Q(\alpha,\beta,\gamma)\right)\cdot\frac{1}{2}\left(\mathcal J\svec{a}^{l}(\xi,\eta,\zeta)+\mathcal J\svec{a}^{l}(\alpha,\beta,\gamma)\right),~~l=1,2,3.
\end{equation}
At a node, we write this contravariant flux as
\begin{equation}
\left(\tilde{\stvec{F}}^{l,\#}_{e}\right)_{(im)jk} \equiv \ssvec{F}^{\#}_{e}\left(\stvec Q_{ijk},\stvec Q_{mjk}\right)\cdot\frac{1}{2}\left(\mathcal J\svec{a}^{l}_{ijk}+\mathcal J\svec{a}^{l}_{mjk}\right),~~l=1,2,3.
\end{equation} It can be simplified if we define the two--point average
\begin{equation}
\aver{u}_{\left(im\right)jk}=\frac{u_{ijk}+u_{mjk}}{2}
\end{equation}
leading to the definition of the divergence of the two--point flux
\begin{equation}
\begin{split}
\mathbb D\left(\cssvec{F}_{e}\right)^{\#}_{ijk} = 2\sum_{m=0}^{N}&\phantom{{}+{}}l'_{m}(\xi_{i})\ssvec{F}_{e}^{\#}\left(Q_{ijk},Q_{mjk}\right)\cdot\aver{\mathcal J\svec{a}^{1}}_{(im)jk}\\
&+l'_m(\eta_j)\ssvec{F}_{e}^{\#}\left(Q_{ijk},Q_{imk}\right)\cdot\aver{\mathcal J\svec{a}^{2}}_{i(jm)k} \\
&+l'_m(\zeta_k)\ssvec{F}_{e}^{\#}\left(Q_{ijk},Q_{ijm}\right)\cdot\aver{\mathcal 
J\svec{a}^{3}}_{ij(km)}.
\end{split}
\label{eq:dg:split-form-op}
\end{equation}

Two--point fluxes are also defined in terms of averages of the two states. 
Leaving off the subscripts for the node locations, we construct two--point entropy conserving 
versions, $\ssvec{F}^{ec}_{e}$, of the fluxes $\ssvec F_{e}^{\#}$ for the incompressible NSE 
using the general formula
\begin{equation}
\stvec{F}^{ec}_{e} = \left(\begin{array}{c}\tilde{\rho u} \\
\tilde{\rho u}\aver{u}+\aver{p} \\
\tilde{\rho u}\aver{v} \\
\tilde{\rho u}\aver{w} \\
\frac{1}{M_0^2}\aver{u} 
\end{array}\right),~~\stvec{G}^{ec}_{e} = \left(\begin{array}{c}\tilde{\rho v} \\
\tilde{\rho v}\aver{u} \\
\tilde{\rho v}\aver{v} +\aver{p}\\
\tilde{\rho v}\aver{w} \\
\frac{1}{M_0^2}\aver{v} 
\end{array}\right),~~\stvec{H}^{ec}_{e} = \left(\begin{array}{c}\tilde{\rho w} \\
\tilde{\rho w}\aver{u} \\
\tilde{\rho w}\aver{v} \\
\tilde{\rho w}\aver{w}+\aver{p} \\
\frac{1}{M_0^2}\aver{w} 
\end{array}\right).
\label{eq:dg:two-point-flux-iNS}
\end{equation}

We consider two options for the approximation of the momentum $\tilde{\rho u_i}, i = 1,2,3$. In the first choice, we take the average of the product, while in the second we take the product of the averages,
\begin{equation}
\tilde{\rho u_i}^{(1)}=\aver{\rho u_i}, \text{or}~\tilde{\rho 
u_i}^{(2)}=\aver{\rho}\aver{u_i}.
\label{eq:stability:momentum-skewsymmetric-options}
\end{equation}

Using the two--point entropy conserving fluxes, the final version of the approximation is written as
\begin{equation}
	\begin{split}
	\left\langle \mathcal J\stvec{Q}_{t},\stvecg{\phi}\right\rangle_{E,N} &+\int_{\partial e,N}\stvecg{\phi}^{T}\left({\ssvec{F}}_{e}^{\star}-\ssvec{F}_{e}\right)\cdot\svec{n}\diff S+\left\langle\mathbb D\left(\cssvec{F}_{e}\right)^{ec},\stvecg{\phi}\right\rangle_{E,N} \\
	&= \int_{\partial e,N}\stvecg{\phi}^{T}\left({\ssvec{F}}_{v}^{\star}-\ssvec{F}_{v}\right)\cdot\svec{n}\diff S+\left\langle\snabla_{\xi}\cdot\cssvec{F}_{v},\stvecg{\phi}\right\rangle_{E,N} + \left\langle \mathcal J \stvec{S},\stvecg{\phi}\right\rangle_{E,N},\\
	\left\langle \mathcal J\ssvec{G},\ssvecg{\psi}\right\rangle_{E,N}&= \int_{\partial e,N}\stvec{W}^{\star,T}\left(\ssvecg{\psi}\cdot\svec{n}\right)\diff S - \left\langle \stvec{W},\snabla_{\xi}\cdot\cssvecg{\psi}\right\rangle_{E,N},
	\end{split}
	\label{eq:dg:scheme-final}
\end{equation}
where we wrote the surface integrals in physical variables using \eqref{eq:dg:surface-integrals-relation}. The evolution equation that one would implement for the coefficients $\stvec{Q}_{ijk}$, and the equation for the gradients $\ssvec{G}_{ijk}$ are obtained by replacing the test functions $\stvecg{\phi}$ and $\ssvecg{\psi}$ by Lagrange basis functions 
$l_i(\xi)l_j(\eta)l_k(\zeta)$,
\begin{equation}
	\begin{split}
\mathcal J_{ijk}\frac{\mathrm{d}\stvec{Q_{ijk}}}{\mathrm{d}t} &+\left(\frac{\delta_{im}}{w_i}\left({\tilde{\stvec{F}}}_{e}^{1,\star}-\tilde{\stvec{F}}_{e}^{1}\right)_{ijk}\!\!\!\!+\frac{\delta_{jm}}{w_j}\left({\tilde{\stvec{F}}}_{e}^{2,\star}-\tilde{\stvec{F}}_{e}^{2}\right)_{ijk}\!\!\!\!+\frac{\delta_{km}}{w_k}\left({\tilde{\stvec{F}}}_{e}^{3,\star}-\tilde{\stvec{F}}_{e}^{3}\right)_{ijk}\right)_{m=0}^{m=N}\\
  &+\mathbb D\left(\cssvec{F}_{e}\right)^{ec}_{ijk}= 
\left(\svec{\nabla}_{\xi}\cdot\cssvec{F}_{v}\right)_{ijk}\\
	&+\left(\frac{\delta_{im}}{w_i}\left({\tilde{\stvec{F}}}_{v}^{1,\star}-\tilde{\stvec{F}}_{v}^{1}\right)_{ijk}\!\!\!\!+\frac{\delta_{jm}}{w_j}\left({\tilde{\stvec{F}}}_{v}^{2,\star}-\tilde{\stvec{F}}_{v}^{2}\right)_{ijk}\!\!\!\!+\frac{\delta_{km}}{w_k}\left({\tilde{\stvec{F}}}_{v}^{3,\star}-\tilde{\stvec{F}}_{v}^{3}\right)_{ijk}\right)_{m=0}^{m=N}\\
	&+ \mathcal J_{ijk} \stvec{S}_{ijk},\\
	\mathcal J_{ijk}\ssvec{G}_{ijk}&=\stvec{W}_{ijk}^{\star}\left(\frac{\delta_{im}}{w_i}\mathcal J\svec{a}_{ijk}^{1}+\frac{\delta_{jm}}{w_j}J\svec{a}_{ijk}^{2}+\frac{\delta_{km}}{w_k}J\svec{a}_{ijk}^{3}\right)_{m=0}^{m=N}\\
	&-\sum_{m=0}^{N}w_m\left(\frac{l_i'(\xi_m)}{w_i}\stvec{W}_{mjk}\mathcal J\svec{a}^{1}_{ijk}+\frac{l_j'(\xi_m)}{w_j}\stvec{W}_{imk} \mathcal J\svec{a}^{2}_{ijk}+\frac{l_k'(\xi_m)}{w_k}\stvec{W}_{ijm}\mathcal 
	J\svec{a}^{3}_{ijk}\right).
	\end{split}
	\label{eq:dg:scheme-nodal}
\end{equation}

Finally, we note that the two--point flux form of the equation is algebraically equivalent to the 
approximation of a split form of the original equations that averages conservative and 
nonconservative forms of the original PDEs. See Appendix \ref{appendix:SplitFormAppendix}. 
Eq. \eqref{eq:dg:scheme-nodal} is integrated in time using a third order low--storage explicit Runge--Kutta 
RK3 scheme \cite{1980:Williamson}.

\subsection{Numerical fluxes}\label{subsec:NumericalFluxes}

The approximation \eqref{eq:dg:scheme-final} is completed with the addition of numerical fluxes ${\stvec{F}}_{e}^{\star}$, ${\stvec{F}}_{v}^{\star}$, and $\stvec{W}^{\star}$, and boundary conditions, the first of which we describe in this section.

For the inviscid fluxes we use the exact Riemann solver derived in \cite{2017:Bassi}. Using the rotational invariance of the flux \cite{2009:Toro}, we write the normal flux as
\begin{equation}
\ssvec{F}_{e}\cdot\svec{n} = \smat{T}^{T}\stvec{F}_{e}\left(\smat{T}\stvec{Q}\right)=\smat{T}^{T}\stvec{F}_{e}\left(\stvec{Q}_{n}\right),~~\smat{T} = \left(\begin{array}{ccccc}1 & 0 & 0 & 0 & 0 \\
0 & n_{x} & n_{y} & n_{z} & 0 \\
0 & t_{1,x} & t_{1,y} & t_{1,z} & 0 \\
0 & t_{2,x} & t_{2,y} & t_{2,z} & 0 \\
0 & 0 & 0 & 0 & 1
\end{array}\right),
\label{eq:Riemann:rotational-invariance}
\end{equation}
where $\smat{T}$ is a rotation matrix that only affects velocities, $\svec{n}=\left(n_x,n_y,n_z\right)$ is the normal unit vector to the face, and $\svec{t}_{1}$ and $\svec{t}_{2}$ are two tangent unit vectors to the face. Recall that $\stvec{F}_{e}$ is the $x$--component of the inviscid flux $\ssvec{F}_{e}=\left(\stvec{F}_{e},\stvec{G}_{e},\stvec{H}_{e}\right)$ \eqref{eq:governing:iNS-inviscid-flux}. When the rotation matrix $\smat{T}$ multiplies the state vector $\stvec{Q}$, we obtain the rotated state vector $\stvec{Q}_{n}$,
\begin{equation}
\stvec{Q}_n = \smat{T}\stvec{Q}=\left(\rho, \rho U_{n}, \rho V_{t1}, \rho V_{t2}, P\right),
\label{eq:Riemann:rotated-state}
\end{equation}
where $U_n = \svec{U}\cdot\svec{n}$ the normal velocity, and $V_{ti}=\svec{U}\cdot\svec{t}_{i}$ ($i=1,2$) are the two tangent velocities. Taken in context, $\rho$  refers to the polynomial approximation of the density. Note that the reference system rotation does not affect the total speed 
\begin{equation}
V_{tot}^2=U^2+V^2+W^2=U_n^2+V_{t1}^2+V_{t2}^2.
\end{equation}

The rotational invariance allows us to transform the 3D Riemann problem into a one dimensional problem for the rotated state vector,
\begin{equation}
\frac{\partial \stvec{Q}_{n}}{\partial t} + \frac{\partial \stvec{F}_{e}\left(\stvec{Q}_{n}\right)}{\partial x} = 0,~~\stvec{Q}_{n}(x,0) = \left\{\begin{array}{ccc}\left(\rho_L,\rho_L U_{nL}, \rho_L V_{t1L}, \rho_L V_{t2L}, P_L \right) & \text{ if }& x \leqslant 0, \\\\
\left(\rho_R,\rho_R U_{nR}, \rho_R V_{t1R}, \rho_R V_{t2R}, P_R \right) & \text{ if }& x > 0,     \end{array} \right.
\end{equation}
whose exact solution is (the details can be found in \cite{2017:Bassi}),

\begin{equation}
\stvec{Q}_{n}^{\star} = \left(\begin{array}{c}\rho^{\star} \\ \rho^{\star}U_n^{\star} \\
\rho^{\star}V_{t1}^{\star} \\
\rho^{\star}V_{t2}^{\star} \\
P^{\star}
\end{array}\right),~~ \stvec{F}_{e}\left(\stvec{Q}^{\star}_{n}\right) = \left(\begin{array}{c}
\rho^{\star}U_n^{\star} \\
\rho^{\star}\left(U_n^{\star}\right)^2 + P^{\star} \\
\rho^{\star}U_n^{\star}V_{t1}^{\star}\\
\rho^{\star}U_n^{\star}V_{t2}^{\star}\\
\frac{1}{M_0^2}\rho^{\star}U_n^{\star}\end{array}\right),
\label{eq:Riemann:star-fluxes}
\end{equation}
with star region solution,
\begin{equation}
\begin{split}
U_n^{\star} &= \frac{P_{L}-P_{R} + \rho_L U_{nL}\lambda_{L}^{+}-\rho_{R}U_{nR}\lambda_{R}^{-}}{\rho_L\lambda_L^{+}-\rho_R \lambda_R^{-}}, ~~ P^{\star} = P_L + \rho_L \lambda^{+}_L\left(U_{nL} - U_n^\star\right),\\
\rho^{\star} &= \left\{\begin{array}{ccc}\rho_{L}^{\star} & \text{ if } & U_n^\star \geqslant 0 \\\rho_{R}^{\star} & \text{ if } & U_n^\star < 0 
\end{array}\right.,~~ \rho_L^{\star} = \frac{\rho_L \lambda_L^{+}}{U_n^{\star}-\lambda_L^-},~~\rho_R^{\star} = \frac{\rho_R \lambda_R^{-}}{U_n^{\star}-\lambda_R^+},~~V_{ti}^{\star}=\left\{\begin{array}{ccc}V_{tiL} & \text{ if } & U_n^\star \geqslant 0 \\V_{tiR} & \text{ if } & U_n^\star < 0 
\end{array}\right. ,
\end{split}
\label{eq:Riemann:star-region-solution}
\end{equation}
and eigenvalues,
\begin{equation}
\lambda_L^{\pm} = \frac{U_{nL} \pm a_L}{2},~~\lambda_R^{\pm} = \frac{U_{nR} \pm a_R}{2},~~ a = \sqrt{U_n^2 + \frac{4}{M_0^2\rho}}.
\end{equation}
Note that in this model, eigenvalues with the positive superscript, $\lambda^{+}_{LR}$, are 
always positive, and eigenvalues with negative superscript , $\lambda^{-}_{LR}$, are 
always negative (i.e. the flow is always subsonic). Following \eqref{eq:Riemann:rotational-invariance}, we multiply 
$\stvec{F}_{e}\left(\stvec{Q}_{n}^{\star}\right)$ in \eqref{eq:Riemann:star-fluxes} by the transposed rotation 
matrix $\smat{T}^{T}$ to obtain the numerical flux $\ssvec{F}_{e}^{\star}\cdot\svec{n}$.

For viscous fluxes we use the Bassi--Rebay 1 (BR1) scheme, which uses the average between adjacent elements in both entropy variables and fluxes,
\begin{equation}
\stvec{W}^{\star} = \aver{\stvec{W}},~~{\ssvec{F}}_{v}^{\star} = \aver{{\ssvec{F}}_{v}}.%= \mathcal B^{\mathcal E}\aver{\ssvec{G}}\cdot\hat{n}
\label{eq:Riemann:BR1}
\end{equation}
\subsection{Boundary conditions}\label{sec:BoundaryConditions}

The full approximation is completed with the addition of boundary conditions. 
Here we show how to impose free-- and no--slip wall boundary conditions.

\subsubsection{Inviscid flux}

The inviscid flux is responsible for canceling the normal velocity $u_n=0$, and 
does not control tangential velocities nor viscous stresses. Thus, the 
procedure for both free-- and no--slip wall boundary conditions is identical.

We consider two ways to enforce the wall boundary condition through the numerical flux $\ssvec{F}_{e}^{\star}$. 
The first directly computes the wall numerical flux. The second sets an artificial reflection state to be used as the external state for the exact Riemann solver.

 If we do not want to use a Riemann solver at the boundaries, we prescribe the zero normal 
velocity through the numerical flux ($U^\star_n=0$), and we take the external pressure from the interior ($P^\star=P$), giving
\begin{equation}
  \ssvec{F}_{e}^{\star}\cdot\svec{n} = \left( \begin{array}{c}\rho U_{n}\\ \rho U_{n}\svec{U}+P\svec{n} 
  \\\frac{1}{M_0^2}U_{n}\end{array}\right)^{\star}=\left( \begin{array}{c}0 \\ P\svec{n} 
  \\0\end{array}\right).
\end{equation}
Alternatively we can use the exact Riemann solver at the physical boundaries, and construct an external state $\stvec{Q}_{n}^{e}$ (rotated using \eqref{eq:Riemann:rotated-state}) taking the interior state $\stvec{Q}^{i}_{n}$ but changing sign of the normal velocity,
\begin{equation}
\stvec{Q}^{i}_n = \left(\begin{array}{c}\rho \\ \rho U_{n}\\ \rho V_{t1}\\ \rho V_{t2}\\ P\end{array}\right),~~\stvec{Q}^{e}_n = \left(\begin{array}{c}\rho \\ -\rho U_{n}\\ \rho V_{t1}\\ \rho V_{t2}\\ P\end{array}\right).
\label{eq:wallbc:riemann-states}
\end{equation}

\subsubsection{Viscous flux}\label{subsubsec:DGSEM:BC-wall-viscous}

For the viscous fluxes, we have to give appropriate values for both the entropy variables 
$\stvec{W}^{\star}$ and the viscous numerical fluxes $\ssvec{F}_{v}^{\star}$ at the 
boundaries. However, only the way the velocities in $\stvec{W}^{\star}$ are specified is important, since the viscous 
fluxes (the stress tensor \eqref{eq:governing:stress-tensor}) are independent of the pressure gradient. 

For the free--slip wall, we use a Neumann boundary condition: we take the entropy variables 
from the interior $\stvec{W}^{\star}=\stvec{W}$ and set viscous fluxes to zero 
$\ssvec{F}_v^{\star}\cdot\svec{n}=\stvec{0}$.

For the no--slip wall, we use a Dirichlet 
boundary condition, thus we use zeroed entropy variables $\stvec{W}^{\star}=\left(0,0,0,0,P\right)$ 
(recall that the pressure is not relevant), and take the viscous fluxes from 
the interior $\ssvec{F}_{v}^{\star} = \ssvec{F}_{v}$.

\section{Stability analysis}\label{sec:DiscreteEntropyAnalysis}

In this section we show that the approximation satisfies the discrete version of the continuous entropy bound \eqref{eq:stability-cont:balance-time}. To do so, we follow \cite{2017:Gassner} and replace $\stvecg{\phi}=\stvec{W}$ and $\ssvecg{\psi}=\ssvec{F}_{v}$ in the discrete system \eqref{eq:dg:scheme-final}, giving
\begin{equation}
\begin{split}
\left\langle \mathcal J\stvec{Q}_{t},\stvec{W}\right\rangle_{E,N} &+\int_{\partial e,N}\stvec{W}^{T}\left({\ssvec{F}}_{e}^{\star}-\ssvec{F}_{e}\right)\cdot\svec{n}\diff S+\left\langle\mathbb D\left(\cssvec{F}_{e}\right)^{ec},\stvec{W}\right\rangle_{E,N} \\
&= \int_{\partial e,N}\stvec{W}^{T}\left({\ssvec{F}}_{v}^{\star}-\ssvec{F}_{v}\right)\cdot\svec{n}\diff S+\left\langle\snabla_{\xi}\cdot\cssvec{F}_{v},\stvec{W}\right\rangle_{E,N},\\
\left\langle \mathcal J\ssvec{G},\ssvec{F}_{v}\right\rangle_{E,N}&= \int_{\partial e,N}\stvec{W}^{\star,T}\left(\ssvec{F}_{v}\cdot\svec{n}\right)\diff S - \left\langle \stvec{W},\snabla_{\xi}\cdot\cssvec{F}_{v}\right\rangle_{E,N}.
\end{split}
\label{eq:stability:scheme-with-testfcns}
\end{equation}
In this analysis we do not consider the source term, as in the continuous case. We can combine both equations of \eqref{eq:stability:scheme-with-testfcns} into a single equation since both share the last quadrature on the right hand side, so
\begin{equation}
\begin{split}
\left\langle \mathcal J\stvec{Q}_{t},\stvec{W}\right\rangle_{E,N} &+\int_{\partial e,N}\stvec{W}^{T}\left({\ssvec{F}}_{e}^{\star}-\ssvec{F}_{e}\right)\cdot\svec{n}\diff S+\left\langle\mathbb D\left(\cssvec{F}_{e}\right)^{ec},\stvec{W}\right\rangle_{E,N} \\
&= \int_{\partial e,N}\left(\stvec{W}^{T}\left({\ssvec{F}}_{v}^{\star}-\ssvec{F}_{v}\right)+\stvec{W}^{\star,T}\ssvec{F}_{v}\right)\cdot\svec{n}\diff S-\left\langle \mathcal J\ssvec{G},\ssvec{F}_{v}\right\rangle_{E,N}.
\end{split}
\label{eq:stability:one-eqn}
\end{equation}

The first term in \eqref{eq:stability:one-eqn} is the discrete integral (quadrature) of the entropy time derivative. Since there is no discrete approximation in time in this analysis (i.e., we consider only semi--discrete stability), the chain rule--based contraction \eqref{eq:stability-cont:iNSAC-entropy-contraction} holds and
\begin{equation}
\left\langle \mathcal J \stvec{Q}_{t},\stvec{W}\right\rangle_{E,N} = \sum_{ijk=0}^{N}w_{ijk}\mathcal J_{ijk}\stvec{W}_{ijk}^{T}\frac{\diff \stvec{Q}_{ijk}}{\diff t}=\sum_{ijk=0}^{N}w_{ijk}\mathcal J_{ijk}\frac{\diff \mathcal E_{ijk}}{\diff t} =\frac{\diff}{\diff t} \left\langle \mathcal J \mathcal E, 1\right\rangle_{E,N}.
\end{equation}

Next, we work on the inviscid flux contribution. In \cite{2017:Gassner}, the authors proved that for arbitrary states $\stvec{Q}_{L}$ and $\stvec{Q}_{R}$, if the two--point entropy conserving flux satisfies 
Tadmor's jump condition \cite{2003:Tadmor,2013:Fisher,2016:gassner},
\begin{equation}
\jump{\stvec{W}^{T}} \stvec{F}^{ec,l}\left(\stvec{Q}_{L},\stvec{Q}_{R}\right)- \jump{\stvec{W}^{T}\stvec{F}^{l}}+\jump{{F}^{\mathcal E,l}} = 0,~~l=1,2,3,
\label{eq:stability:two-point-entropy-contraction}
\end{equation}
with 
\begin{equation}
\jump{u} = u_R-u_L,
\label{eq:stability:jump-def}
\end{equation}
being the jump operator, then the inner product of the split--form divergence with the entropy variables satisfies,
\begin{equation}
\left\langle\mathbb D\left(\cssvec{F}_{e}\right)^{ec},\stvec{W}\right\rangle_{E,N}  = \int_{\partial E,N}\svec{\tilde{F}}^{\mathcal E}\cdot\hat{n}\diff S_{\xi}=\int_{\partial e,N}\svec{{F}}^{\mathcal E}\cdot\svec{n}\diff S.
\label{eq:stability:split-form-to-boundary}
\end{equation}
To show that the two--point flux \eqref{eq:dg:two-point-flux-iNS} satisfies \eqref{eq:stability:two-point-entropy-contraction}, we replace the 
entropy conserving flux \eqref{eq:dg:two-point-flux-iNS} and the entropy 
variables, and see that
\begin{equation}
  \begin{split}
\jump{\stvec{W}^{T}}  \stvec{F}^{ec,l} -\jump{\stvec{W}^{T}\stvec{F}^{l}}=& -\frac{1}{2}\jump{V_{tot}^2}\tilde{\rho 
  U_{l}}+\tilde{\rho U_l}\left(\aver{U}\jump{U}+\aver{V}\jump{V}+\aver{W}\jump{W}\right) + \aver{P}\jump{U_l} 
  \\
  &+\aver{U_l}\jump{P}+\frac{1}{2}\jump{\rho U_l V_{tot}^2 }-\jump{\rho \left(U^2+V^2+W^2\right)U_l 
  }-\jump{U_l p}-\jump{U_l P} \\
  =&-\frac{1}{2}\tilde{\rho U_l}\jump{U^2+V^2+W^2} + \frac{1}{2}\tilde{\rho 
  u_l}\jump{U^2+V^2+W^2}+\jump{P U_{l}}\\
  &-\jump{\left(\frac{1}{2}\rho \left(U^2+V^2+W^2\right)+P\right)U_{l}}-\jump{P 
  U_{l}}\\
  =&-\jump{\left(\frac{1}{2}\rho \left(U^2+V^2+W^2\right)+P\right)U_{l}} = -\jump{F^{\mathcal 
  E,l}},
  \end{split}
  \label{eq:stability:entropy-flux-contraction-proof}
\end{equation}
where $F^{\mathcal E,l}$ is the entropy flux 
\eqref{eq:stability-cont:entropy-flux}. In the process, we used the two arithmetic properties of the 
average and jump operators,
\begin{equation}
\jump{ab}=\aver{a}\jump{b}+\jump{a}\aver{b},~~\aver{a}\jump{a}=\frac{1}{2}\jump{a^2}.
\label{eq:stability:average-jump-relation}
\end{equation}
Note that \eqref{eq:stability:entropy-flux-contraction-proof} holds 
independently of which momentum approximation $\tilde{\rho u_l}$
\eqref{eq:stability:momentum-skewsymmetric-options} is used, so from the point of view of stability, either approximation in \eqref{eq:stability:momentum-skewsymmetric-options} is acceptable.

Thus, using \eqref{eq:stability:split-form-to-boundary}, \eqref{eq:stability:one-eqn} becomes,
\begin{equation}
\begin{split}
\frac{\diff}{\diff t} \left\langle \mathcal J \mathcal E, 1\right\rangle_{E,N}&+\int_{\partial e,N}\left(\stvec{W}^{T}\left({\stvec{F}}_{e}^{\star}-\ssvec{F}_{e}\right)+\svec{{F}}^{\mathcal E}\right)\cdot\svec{n}\diff S \\
&= \int_{\partial e,N}\left(\stvec{W}^{T}\left({\stvec{F}}_{v}^{\star}-\ssvec{F}_{v}\right)+\stvec{W}^{\star,T}\ssvec{F}_{v}\right)\cdot\svec{n}\diff S-\left\langle \mathcal J\ssvec{G},\mathcal B^{\mathcal E}\ssvec{G}\right\rangle_{E,N},
\end{split}
\label{eq:stability:stability-with-twopoint}
\end{equation}
where we have written the viscous flux as a function of the entropy variable gradient following \eqref{eq:stability-cont:viscous-flux-Y=AX}, and we 
wrote the surface integrals in physical coordinates using \eqref{eq:dg:surface-integrals-relation}. The former allows us to bound the viscous flux volume contribution using the viscous positive semi--definiteness property \eqref{eq:entropy-stability:viscous-matrices-pos-def},
\begin{equation}
\left\langle \mathcal J\ssvec{G},\mathcal B^{\mathcal E}\ssvec{G}\right\rangle_{E,N} \geqslant \min_{E,N}\left(\mathcal J\right)\left\langle\ssvec{G},\mathcal B^{\mathcal E}\ssvec{G}\right\rangle_{E,N} \geqslant 0,
\label{eq:stability:viscous-semidefinite-positive}
\end{equation}
since all nodal values of the Jacobian, $\mathcal J_{ijk}$, are strictly positive in an admissible quality mesh. Now, \eqref{eq:stability:viscous-semidefinite-positive} is the discrete version of \eqref{eq:stability-cont:viscous-strain-tensor}, therefore, it can be written as
\begin{equation}
\left\langle\mathcal J\ssvec{G},\mathcal B^{\mathcal E}\ssvec{G}\right\rangle_{E,N}=\frac{2}{\Rey}\left\langle\mathcal J \mathcal S,\mathcal S\right\rangle_{E,N}\geqslant 0,
\end{equation}
where $\mathcal S$ is the strain tensor computed from the approximated entropy variable gradient $\mathcal S=\tens{S}\left(\ssvec{G}\right)$. 

Now that the volume terms have been discussed and shown to match the continuous versions, we address element boundary terms. Stability of the boundary terms only makes sense when all elements are considered. Therefore, we sum \eqref{eq:stability:stability-with-twopoint} over all the elements, to get the time rate of change of the total entropy,
\begin{equation}
\frac{\diff\bar{\mathcal E}}{\diff t} + \text{IBT} + \text{PBT} = -\sum_{e}\frac{2}{\Rey}\left\langle\mathcal J \mathcal S,\mathcal 
S\right\rangle_{E,N}\leqslant 0,
\label{eq:stability:entropy-bound-all-el}
\end{equation}
where
\begin{equation}
\bar{\mathcal E} = \sum_{e}\left\langle \mathcal J \mathcal E, 1\right\rangle_{E,N}
\end{equation}
is the total entropy, IBT is the contribution of the interior faces to the surface integral,
\begin{equation}
  \begin{split}
\text{IBT} =&-\sum_{\interiorfaces}\int_{N}\left(\jump{\stvec{W}^{T}}\ssvec{F}^{\star}_{e}+\jump{\svec{F}^{\mathcal E}}-\jump{\stvec{W}^{T}\ssvec{F}_{e}}\right)\cdot\svec{n}_L\diff 
S \\
&+ \sum_{\interiorfaces}\int_{N}\left(\jump{\stvec{W}^T}\ssvec{F}^{\star}_{v}+\stvec{W}^{\star,T}\jump{\ssvec{F}_{v}}-\jump{\stvec{W}^{T}\ssvec{F}_{v}}\right)\cdot\svec{n}_{L}\diff 
S \\
=&\mathrm{IBT}_{e} + \mathrm{IBT}_{v},
 \end{split}
\label{eq:stability:IBT}
\end{equation}
and PBT is the physical boundary contribution,
\begin{equation}
  \begin{split}
    \text{PBT} =&\phantom{{}+{}}\sum_{\boundaryfaces}\int_{N}\left(\stvec{W}^{T}\left(\ssvec{F}^{\star}_{e}-\ssvec{F}_{e}\right)+\svec{F}^{\mathcal 
    E}\right)\cdot\svec{n}\diff S \\
    &-\sum_{\boundaryfaces}\int_{N}\left(\stvec{W}^{T}\left(\ssvec{F}^{\star}_{v}-\ssvec{F}_{v}\right)+\stvec{W}^{\star,T}\ssvec{F}_{v}\right)\cdot\svec{n}\diff S 
    \\
    =&\mathrm{PBT}_{e} - \mathrm{PBT}_{v}.
  \end{split}
  \label{eq:stability:PBT}
\end{equation}

In \eqref{eq:stability:IBT}, we chose to write the contributions of the elements on the left and the right  of the surface integrals taking the left face normal $\svec{n}_{L}$ as a reference. Thus, the right side face terms are added with opposing sign since $\svec{n}_R = -\svec{n}_L$. As a result, we get the jumps in the variables as defined in \eqref{eq:stability:jump-def}. Note that the numerical fluxes $\ssvec{F}_{e}^{\star}$ and $\ssvec{F}_{v}^{\star}$ factor out of the jump operators since they are shared by both left and right states. 
%In 
%Sec. \ref{subsec:IBT} and \ref{subsec:PBT-wall}, we show that both interior and boundary operators are stable,
%%
%\begin{equation}
%  \text{IBT} \geqslant 0,~~\text{PBT}\geqslant 0,
%\end{equation}
%%
%thereby obtaining a discrete entropy bound 
%%
%\begin{equation}
%  \frac{\diff \bar{\mathcal E}}{\diff t} \leqslant -\sum_{e}\frac{2}{\Rey}\left\langle\mathcal J \mathcal S,\mathcal 
%  S\right\rangle_{E,N},
%\end{equation}
%%
%which is consistent with the continuous bound 
%\eqref{eq:stability-cont:balance-final}.

\subsection{Stability of the Interior Boundary Terms (IBT)}\label{subsec:IBT}

Following \eqref{eq:stability:entropy-bound-all-el}, interior boundary terms \eqref{eq:stability:IBT} are stable if $\text{IBT}\geqslant 0$.  
For inviscid fluxes we consider two possibilities: using the two--point entropy conserving flux 
\eqref{eq:dg:two-point-flux-iNS}, or the exact Riemann solver 
\eqref{eq:Riemann:star-region-solution}.

\subsubsection{Inviscid fluxes: two--point entropy conserving numerical flux}

If we use the two--point entropy conserving flux \eqref{eq:dg:two-point-flux-iNS} as the numerical flux $\ssvec{F}^{\star}_{e} 
=\ssvec{F}^{ec}_e$, we obtain $\mathrm{IBT}_{e}=0$ as a result of 
Tadmor's jump condition \eqref{eq:stability:two-point-entropy-contraction},
\begin{equation}
\text{IBT}_{e}=-\sum_{\interiorfaces}\int_{N}\left(\jump{\stvec{W}^{T}}\ssvec{F}^{ec}_{e}+\jump{\svec{F}^{\mathcal E}}-\jump{\stvec{W}^{T}\ssvec{F}_{e}}\right)\cdot\svec{n}_{L}\diff 
S = 0.
\label{eq:stability:IBT-e-central}
\end{equation}

\subsubsection{Inviscid fluxes: exact Riemann solver}
It is not trivial to show the stability of the inviscid contribution to IBT using the exact Riemann solver to compute $\ssvec{F}_{e}^{\star}$. We move that analysis to Appendix \ref{appendix:IBT-exactRiemann}, where we prove that
\begin{equation}
\text{IBT}_{e}=-\sum_{\interiorfaces}\int_{N}\left(\jump{\stvec{W}^{T}}\ssvec{F}^{\star}_{e}+\jump{\svec{F}^{\mathcal E}}-\jump{\stvec{W}^{T}\ssvec{F}_{e}}\right)\cdot\svec{n}_{L}\diff 
S\geqslant 0.
\label{eq:stability:IBT-e}
\end{equation}

\subsubsection{Viscous fluxes: BR1 method}

For the viscous fluxes we write viscous contribution in IBT by replacing the numerical values using the BR1 scheme \eqref{eq:Riemann:BR1},
\begin{equation}
\text{IBT}_{v}= \sum_{\interiorfaces}\int_{N}\left(\jump{\stvec{W}^T}\aver{\ssvec{F}_{v}}+\aver{\stvec{W}^{T}}\jump{\ssvec{F}_{v}}-\jump{\stvec{W}^{T}\ssvec{F}_{v}}\right)\cdot\svec{n}_L\diff 
S=0,
\end{equation}
since the algebraic identity \eqref{eq:stability:average-jump-relation} holds.

Thus, we conclude that interior boundary terms  are stable, $\text{IBT}\geqslant 0$. 
Moreover, we get an entropy conserving scheme using the two--point entropy flux \eqref{eq:dg:two-point-flux-iNS}
as the numerical flux, $\mathrm{IBT}=0$, and an entropy stable scheme if we use 
the exact Riemann solver, $\mathrm{IBT}\geqslant 0$.

\subsection{Physical Boundary Terms (PBT): wall boundary conditions}\label{subsec:PBT-wall}

Like at interior boundaries, boundary condition prescriptions are stable if $\text{PBT}\geqslant 0$. In this section we 
follow \cite{2019:Hinderlang} and show that the two approaches described in Sec. \ref{sec:BoundaryConditions}  stably impose both free-- and no--slip wall boundary conditions.

\subsubsection{Inviscid fluxes}
The stability condition for the inviscid boundary flux is,
\begin{equation}
  \text{PBT}_{e}=\sum_{\boundaryfaces}\int_{N}\left(\stvec{W}^{T}\left(\ssvec{F}^{\star}_{e}-\ssvec{F}_{e}\right)\cdot\svec{n}+\svec{F}^{\mathcal 
    E}\cdot\svec{n}\right)\diff S  \geqslant 0.
 \label{eq:PBTeIntegralCondition}
\end{equation}
A sufficient condition for  \eqref{eq:PBTeIntegralCondition} to hold is that the argument
\begin{equation}
\Delta_{e}=  \stvec{W}^{T}\left(\ssvec{F}^{\star}_{e}-\ssvec{F}_{e}\right)\cdot\svec{n}+\svec{F}^{\mathcal 
    E}\cdot\svec{n}\geqslant 0.
    \label{eq:wallbc:euler-inequality}
\end{equation}
We presented two choices by which to enforce the wall boundary condition through the numerical flux $\ssvec{F}_{e}^{\star}$, and examine their stability in turn.
\begin{enumerate}
\item If we do not to use a Riemann solver at the boundaries, we write \eqref{eq:wallbc:euler-inequality} replacing the 
fluxes,
\begin{equation}
 \Delta_{e}= \left(-\frac{1}{2}V_{tot}^{2},\svec{u},M_0^2 P\right)\left[\left(\begin{array}{c}0 \\ P\svec{n} \\ 0\end{array}\right)-\left(\begin{array}{c}\rho U_{n} \\ \rho U_{n}\svec{U}+P\svec{n} \\ \frac{1}{M_0^2}U_n\end{array}\right)\right] 
  + \left(\frac{1}{2}\rho V_{tot}^2 + P\right)U_n=0.
\end{equation}
Therefore, enforcing the wall boundary condition through direct prescription of the numerical flux is neutrally stable $\text{PBT}=0$.

\item If we use the exact Riemann solver at the physical boundaries and we construct the external state \eqref{eq:wallbc:riemann-states}, we construct the star region solution \eqref{eq:Riemann:star-region-solution} for the particular states  $\stvec{Q}_{n}^{e}$ ,
\begin{equation}
U^{\star} = \frac{P-P +\rho \lambda^+U_n- \rho\left(-\lambda^{+}\right)\left(-U_n\right)}{\rho \lambda^+-\rho\left(-\lambda^{+}\right)}=0 ,~~P^{\star} = P + \rho \lambda^{+}U_n,
\end{equation}
where,
\begin{equation}
\lambda^{+}_{L} = \frac{U_n+a}{2} = \lambda^{+}\geq 0, ~~ \lambda^{-}_{R} = \frac{-U_n-a}{2} = -\lambda^{+}<0.
\end{equation}
Therefore,
\begin{equation}
\stvec{F}_{e}\left(\stvec{Q}_n^{\star}\right) = \left(\begin{array}{c}\rho^{\star}U_{n}^{\star} \\\rho^\star \left(U_n^\star\right)^2+ P^\star \\ \rho^\star U_n^\star V_{t1}^{\star} \\ \rho^\star V_n^\star v_{t2}^{\star} \\ \frac{1}{M_0^2}U_n^\star \end{array}\right)=\left(\begin{array}{c}0 \\  P + \rho \lambda^{+}U_n \\ 0 \\ 0 \\ 0 \end{array}\right).
\end{equation}
Now we write \eqref{eq:wallbc:euler-inequality} using the rotational invariance \eqref{eq:Riemann:rotational-invariance}, $\ssvec{F}_{e}^{\star}\cdot\svec{n} = \smat{T}^{T}\stvec{F}_{e}\left(\stvec{Q}_{n}^{\star}\right)$, and $\ssvec{F}_{e}\cdot\svec{n}=\smat{T}^{T}\stvec{F}_{e}\left(\stvec{Q}_{n}\right)$,
\begin{equation}
\begin{split}
\Delta_{e}&=  \stvec{W}^{T}\smat{T}^{T}\left(\stvec{F}_{e}\left(\stvec{Q}_{n}^{\star}\right)-\stvec{F}_{e}\left(\stvec{Q}_{n}^{i}\right)\right)+\svec{F}^{\mathcal 
	E}\cdot\svec{n} \\
&=\stvec{W}_{n}^{T}\left(\stvec{F}_{e}\left(\stvec{Q}_{n}^{\star}\right)-\stvec{F}_{e}\left(\stvec{Q}_{n}^{i}\right)\right)+\svec{F}^{\mathcal 
	E}\cdot\svec{n} \\
&=\left(-\frac{1}{2}V_{tot}^2,U_n,V_{t1},V_{t2},P\right)\left(\begin{array}{c}-\rho U_n \\  P + \rho \lambda^{+}U_n-\rho U_n^2 - P \\ -\rho U_n V_{t1} \\ -\rho U_n V_{t2} \\ -\frac{1}{M_0^2}U_n \end{array}\right)+\frac{1}{2}\rho V_{tot}^2 U_n + PU_n \\
&=\frac{1}{2}\rho V_{tot}^2U_n + \rho \lambda^{+}U_n^2 - \rho U_n\left(U_n^2 + V_{t1}^2 + V_{t2}^{2}\right) - PU_n + \frac{1}{2}\rho V_{tot}^2U_n + P U_n\\
&=\rho \lambda^{+}U_n^2 \geqslant 0,
\end{split}
\end{equation}
which is always dissipative since $\lambda^{+}\geqslant 0$. The dissipation introduced at the boundaries increases with the square of the velocity normal to the wall (which vanishes when the boundary condition is exactly satisfied, $U_n=0$).
\end{enumerate}

We conclude that wall boundary conditions can be stably enforced either by direct imposition through the numerical flux, or using the exact Riemann solver constructing an external state. The former is neutrally stable and does not add any dissipation, the latter introduces numerical dissipation that vanishes as the normal velocity converges (weakly).

\subsubsection{Viscous fluxes}

The difference between free-- and no--slip boundary condition rests on the 
entropy variables and viscous numerical flux implentation. We will see the effect of the 
choices given in Sec. \ref{subsubsec:DGSEM:BC-wall-viscous} on viscous fluxes 
stability,

\begin{equation}
\begin{split}
    \text{PBT}_{v} =&-\sum_{\boundaryfaces}\int_{N}\left(\stvec{W}^{T}\ssvec{F}^{\star}_{v}+\stvec{W}^{\star,T}\ssvec{F}_{v}-\stvec{W}^{T}\ssvec{F}_{v}\right)\cdot\svec{n}\diff 
    S.
  \end{split}
  \end{equation}

For the free--slip wall, the choice $\stvec{W}^{\star}=\stvec{W}$ and $\ssvec{F}_v^{\star}\cdot\svec{n}=\stvec{0}$ 
implies that the viscous contribution is zero,
\begin{equation}
\begin{split}
    \text{PBT}_{v} =&-\sum_{\boundaryfaces}\int_{N}\left(\stvec{W}^{T}\ssvec{F}^{\star}_{v}+\stvec{W}^{\star,T}\ssvec{F}_{v}-\stvec{W}^{T}\ssvec{F}_{v}\right)\cdot\svec{n}\diff S \\
    =&-\sum_{\boundaryfaces}\int_{N}\left(\stvec{W}^{T}\ssvec{F}_{v}-\stvec{W}^{T}\ssvec{F}_{v}\right)\cdot\svec{n}\diff S 
    =0.
    \end{split}
\end{equation}
The same result holds for the no--slip wall choice 
$\stvec{W}^{\star}=(0,0,0,0,P)$ and 
$\ssvec{F}_v^{\star}\cdot\svec{n}=\ssvec{F}_v\cdot\svec{n}$.

\subsection{Final remarks}\label{subsec:Stability:Conclusions}

We prove that the split--form DG scheme with the two--point entropy conserving 
flux \eqref{eq:dg:two-point-flux-iNS} satisfies the discrete entropy law 
\eqref{eq:stability:entropy-bound-all-el}. The boundary terms are identically 
zero if the two--point entropy conserving flux is used as the numerical flux at 
interior boundaries and if the wall boundary condition is enforced via direct 
prescription, and stable if we use the exact Riemann solver 
\eqref{eq:Riemann:star-region-solution}. The former yields an entropy conserving 
scheme (where all the dissipation is due to the physical viscosity), while the 
latter produces an entropy stable (dissipative) scheme, where the physical viscosity is 
complemented with numerical dissipation at element boundaries.

\section{Numerical experiments}\label{sec:NumericalExperiments}

The purpose of this section is to address the accuracy and robustness of the 
method. For the former, we perform a convergence analysis in 
three--dimensions using the manufactured solution method, and in two--dimensions 
using the Kovasznay flow problem \cite{1948:Kovasznay}. For the latter, we solve 
the three--dimensional inviscid Taylor--Green vortex, and the two--dimensional 
Rayleigh--Taylor instability \cite{2000:Guermond}.

As noted in Sec. \ref{subsec:DGSEM:Approximation}, in all the numerical experiments  we use a third order explicit low--storage 
Runge--Kutta scheme \cite{1980:Williamson}. 

\subsection{Convergence study}\label{subsec:num:conv}

In this section, we address the accuracy of the scheme by considering a 
manufactured solution with a divergence free velocity field,
\begin{equation}
  \begin{array}{ll}
    \rho_0 &= 1, \\
    u_0 &= \cos\left(\pi(x+y+z-2t)\right), \\
    v_0 &= -2\cos\left(\pi(x+y+z-2t)\right), \\
    w_0 &=\cos\left(\pi(x+y+z-2t)\right), \\
    p_0 &= 2\cos\left(\pi(x+y+z-2t)\right) -3\frac{1}{\Rey} \pi \sin\left(\pi(x+y+z-2t)\right),\\ 
  \end{array}
  \label{eq:num:conv:man-sol}
\end{equation}
with $\Rey=1000$, and the corresponding source term,
\begin{equation}
  \begin{array}{ll}
    q_\rho &= 0, \\
    q_{\rho u} &= 0 \\
    q_{\rho v} &= -6\pi\sin\left(\pi(x+y+z-2t)\right)-\frac{9}{\Rey}\pi^2\cos\left(\pi(x+y+z-2t)\right), \\
    q_{\rho w} &=0, \\
    q_p &= 4\pi\sin\left(\pi(x+y+z-2t)\right)+\frac{6}{\Rey}\pi^2\cos\left(\pi(x+y+z-2t)\right),\\ 
  \end{array}
  \label{eq:num:conv:source}
\end{equation}
solved in a fully periodic box with size $[-1,1]^3$. The artificial compressibility Mach number is $M_0^2=10^{-3}$. 
In these simulations, we use a $\text{CFL}$ type condition, with the local maximum wave speed and 
relative grid size $\Delta x/(N+1)$. We use $\text{CFL}=0.5$ for all 
computations, and the end time is set to $t_{F}=10$. After time integration, we 
compute the $L^2$ error as,

\begin{equation}
  \Vert u-u_0\Vert_{L^2} = \sqrt{\sum_{e}\sum_{i,j,k=0}^{N}w_iw_jw_k (
  u_{ijk}-u_{0,ijk})^2 \mathcal J_{ijk}}
\end{equation}

\begin{figure}[h]
  \centering
\subfigure[$\rho$ error]{\includegraphics[scale=0.2]{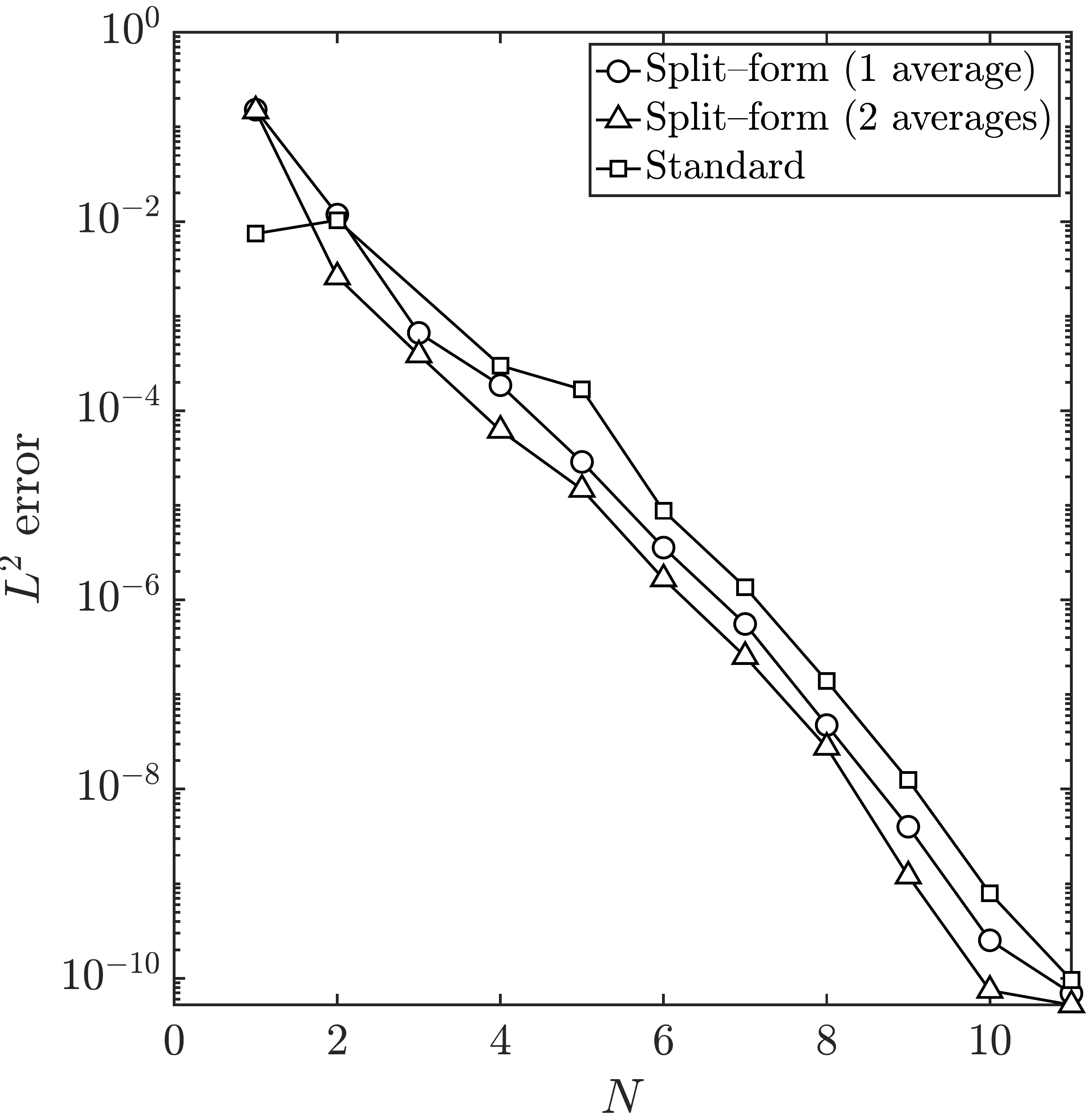}\label{fig:num:conv:p-Conv:rho}}
\subfigure[$\rho u$ and $\rho w$ errors]{\includegraphics[scale=0.2]{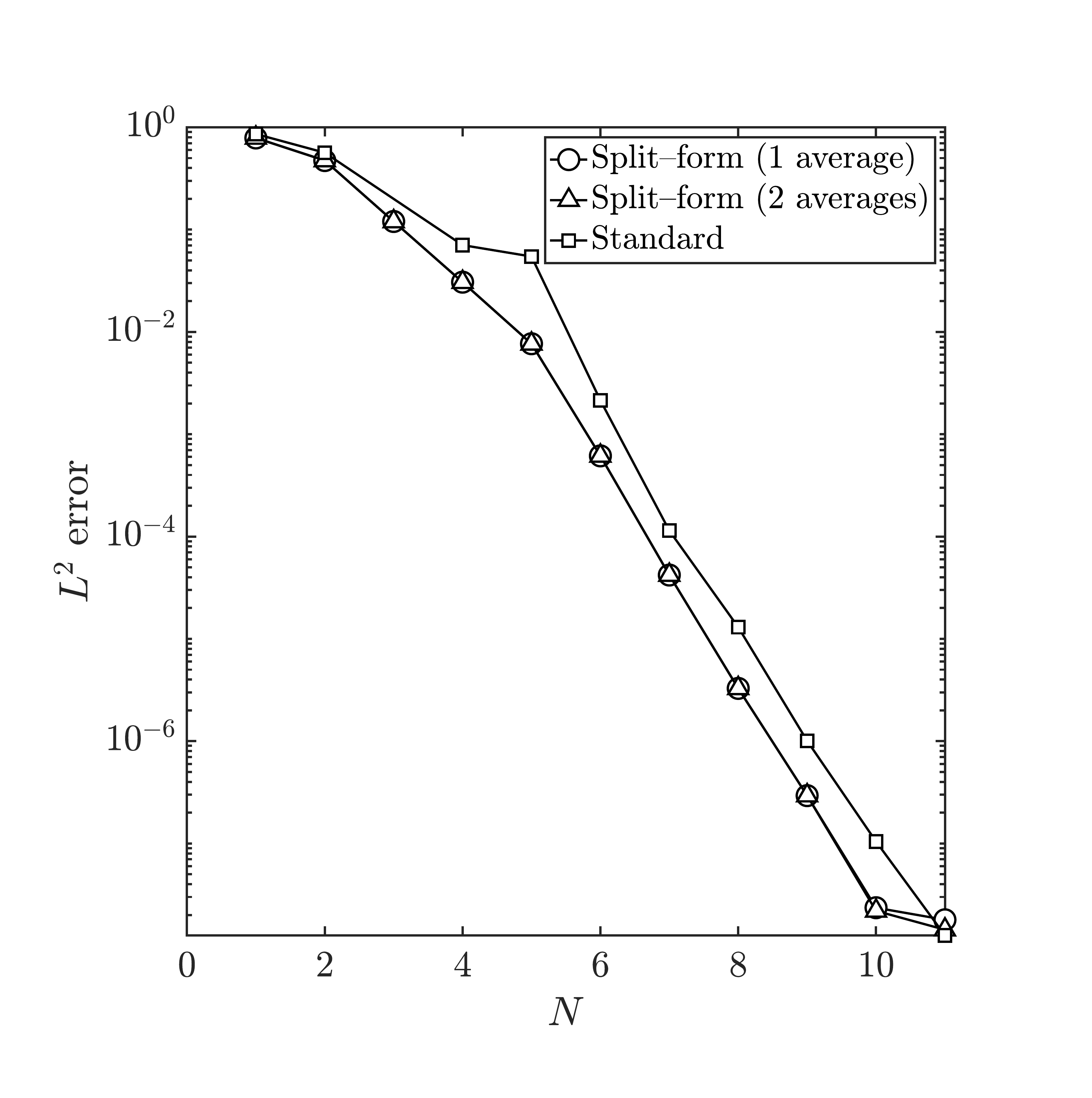}\label{fig:num:conv:p-Conv:rhouw}}
\subfigure[$\rho v$ error]{\includegraphics[scale=0.2]{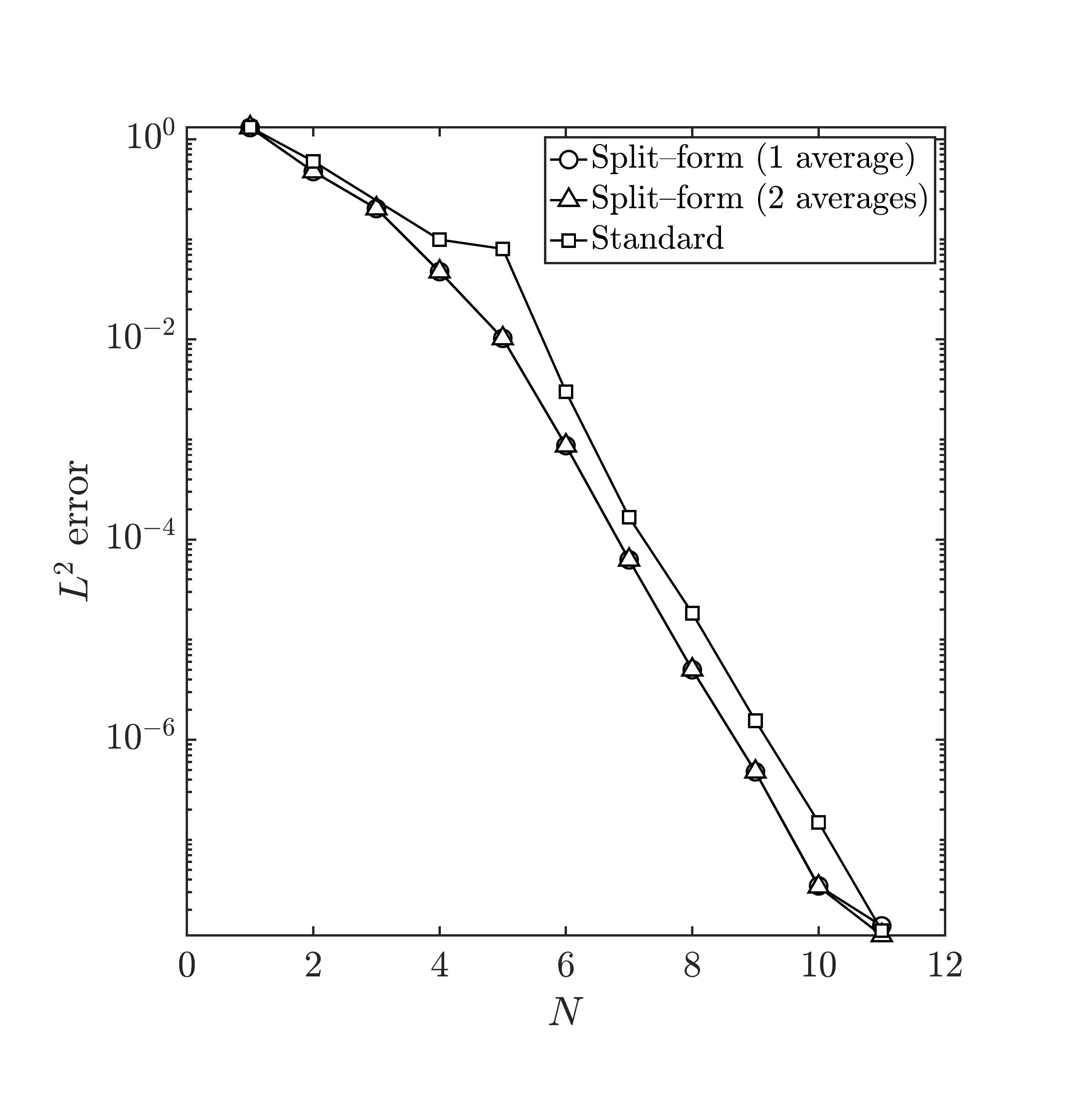}\label{fig:num:conv:p-Conv:rhov}}
\subfigure[$p$ error]{\includegraphics[scale=0.2]{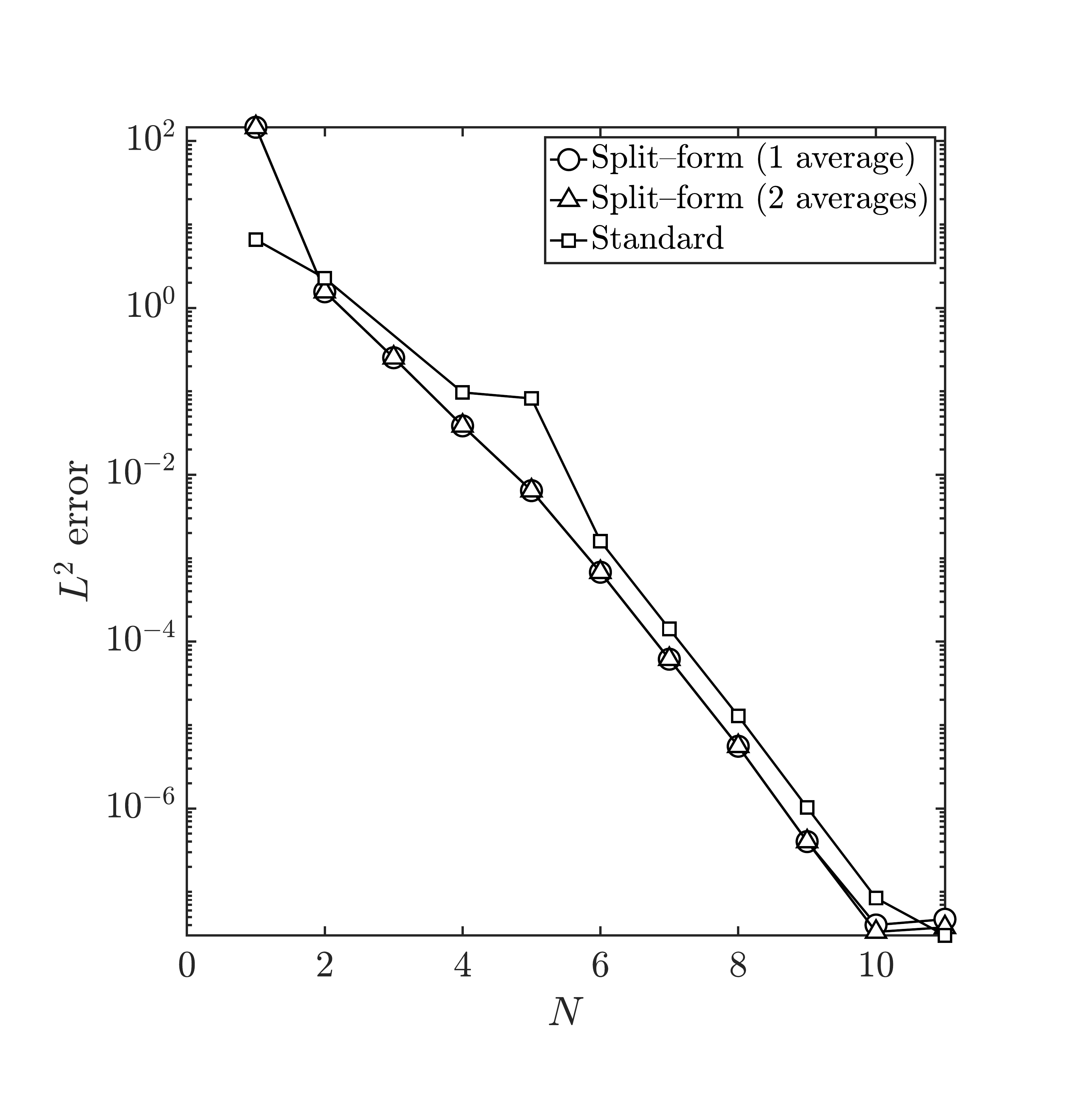}\label{fig:num:conv:p-Conv:p}}
  \caption{Polynomial order convergence study using the manufactured solution \eqref{eq:num:conv:man-sol} on a $4^3$ cartesian grid. We consider three schemes: split--form DG
  with one and two averages in momentum, and standard DG. All three schemes show spectral convergence behavior. The uneven rate of convergence is an effect of the BR1 scheme}
  \label{fig:num:conv:p-Conv}
\end{figure}

We first perform a polynomial order study on a regular a $4^3$ mesh, with the  polynomial order ranging from $N=1$ to $N=11$. We solve the incompressible NSE with the three schemes 
presented in this work: split--form with both one ($\marksymbol{o}{black}$) and two ($\marksymbol{triangle}{black}$) averages in momentum, and 
standard DG ($\marksymbol{square}{black}$). We show the $L^2$ errors obtained for the five state variables in Figs. \ref{fig:num:conv:p-Conv:rho}--\ref{fig:num:conv:p-Conv:p}.
(Because of the symmetry, the $\rho u$ and $\rho w$ $L^2$ errors are identical, both represented in Fig. \ref{fig:num:conv:p-Conv:rhouw}.)
We find that the $L^2$ errors are systematically smaller with the split--form, 
and the convergence is smoother than standard DG. Both remain stable in this 
smooth problem. The slight non--optimality found in the convergence rates is a result 
of the BR1 scheme even--odd behavior (note that the Reynolds number is moderate, $\Rey=1000$), which we have not 
found when solving the purely inviscid problem.

Next, we address the mesh convergence and consider different meshes: from $4^3$ to $16^3$, and polynomial orders: 
from $N=2$ to $N=5$. For this experiment, we use the split--form with one average, and maintain the same parameters used 
in the polynomial order study. The results are reported in Table \ref{tab:num:conv:h-conv}, 
with the $L^2$ errors obtained for each of the five state variables, and the 
estimated order of convergence, showing good agreement with the reference \cite{2017:Bassi}. 
As in \cite{2017:Bassi}, we note that the pressure convergence order is systematically 
smaller than the others.

{\begin{table}[h]
  \caption{Manufactured solution \eqref{eq:num:conv:man-sol} convergence analysis: we use $4^3$, $6^3$, $8^3$, $12^3$ and $16^3$ meshes, and $N=2,3,4$ and 5.
  We select the split--form DG with one average in momentum, and run all simulations until $t_F=10$ using the explicit RK3 scheme with $CFL=0.5$}
  \label{tab:num:conv:h-conv}
\resizebox{\textwidth}{!}{%  
\begin{tabular}{llllllllll}
\hline
      & Mesh   & $\rho$ error & order & $\rho u,w$ error & order & $\rho v$ error & order & $p$ error & order \\ \hline
N=2 & $4^3$ & 1.18E-2 & -- & 4.73E-1 & -- & 4.75E-1 & -- & 1.56E-0 & -- \\ 
 & $6^3$ & 3.29E-3 & 3.16 & 2.12E-1 & 1.98 & 2.16E-1 & 1.94 & 7.23E-1 & 1.89 \\ 
 & $8^3$ & 1.55E-3 & 2.61 & 9.98E-2 & 2.62 & 1.11E-1 & 2.32 & 4.26E-1 & 1.84 \\ 
 & $12^3$ & 5.44E-4 & 2.58 & 2.88E-2 & 3.07 & 3.32E-2 & 2.98 & 1.91E-1 & 1.98 \\ 
 & $16^3$ & 2.39E-4 & 2.86 & 1.24E-2 & 2.92 & 1.56E-2 & 2.63 & 9.90E-2 & 2.29 \\ 
N=3 & $4^3$ & 6.65E-4 & -- & 1.20E-1 & -- & 2.03E-1 & -- & 2.53E-1 & -- \\ 
 & $6^3$ & 1.58E-4 & 3.55 & 2.57E-2 & 3.80 & 4.72E-2 & 3.60 & 5.82E-2 & 3.62 \\ 
 & $8^3$ & 6.79E-5 & 2.92 & 9.59E-3 & 3.42 & 1.70E-2 & 3.54 & 2.51E-2 & 2.92 \\ 
 & $12^3$ & 1.45E-5 & 3.80 & 2.24E-3 & 3.58 & 3.13E-3 & 4.18 & 6.79E-3 & 3.23 \\ 
 & $16^3$ & 4.44E-6 & 4.12 & 8.03E-4 & 3.57 & 1.13E-3 & 3.55 & 2.56E-3 & 3.39 \\ 
N=4 & $4^3$ & 1.86E-4 & -- & 3.06E-2 & -- & 4.75E-2 & -- & 3.85E-2 & -- \\ 
 & $6^3$ & 3.04E-5 & 4.47 & 3.00E-3 & 5.73 & 4.87E-3 & 5.62 & 6.78E-3 & 4.28 \\ 
 & $8^3$ & 8.02E-6 & 4.63 & 7.66E-4 & 4.74 & 1.20E-3 & 4.88 & 1.90E-3 & 4.43 \\ 
 & $12^3$ & 9.92E-7 & 5.16 & 8.34E-5 & 5.47 & 1.40E-4 & 5.29 & 2.76E-4 & 4.76 \\ 
 & $16^3$ & 2.02E-7 & 5.54 & 1.79E-5 & 5.35 & 3.07E-5 & 5.27 & 6.62E-5 & 4.96 \\ 
N=5 & $4^3$ & 2.88E-5 & -- & 7.67E-3 & -- & 1.02E-2 & -- & 6.44E-3 & -- \\ 
 & $6^3$ & 2.29E-6 & 6.24 & 4.59E-4 & 6.94 & 7.10E-4 & 6.58 & 5.42E-4 & 6.11 \\ 
 & $8^3$ & 4.59E-7 & 5.59 & 6.67E-5 & 6.70 & 1.06E-4 & 6.61 & 1.10E-4 & 5.54 \\ 
 & $12^3$ & 3.78E-8 & 6.15 & 5.57E-6 & 6.13 & 8.89E-6 & 6.11 & 1.11E-5 & 5.66 \\ 
 & $16^3$ & 6.14E-9 & 6.32 & 1.07E-6 & 5.73 & 1.77E-6 & 5.61 & 2.09E-6 & 5.81 \\ 
\hline
\end{tabular}}
\label{tab:num:conv:h-conv}
\end{table}
}

In conclusion, the scheme and its implementation show the expected 
convergence behavior.

\subsection{Kovasznay test case}\label{subsec:num:Kovasznay}

\begin{figure}[h]
  \centering
  \includegraphics[scale=0.2]{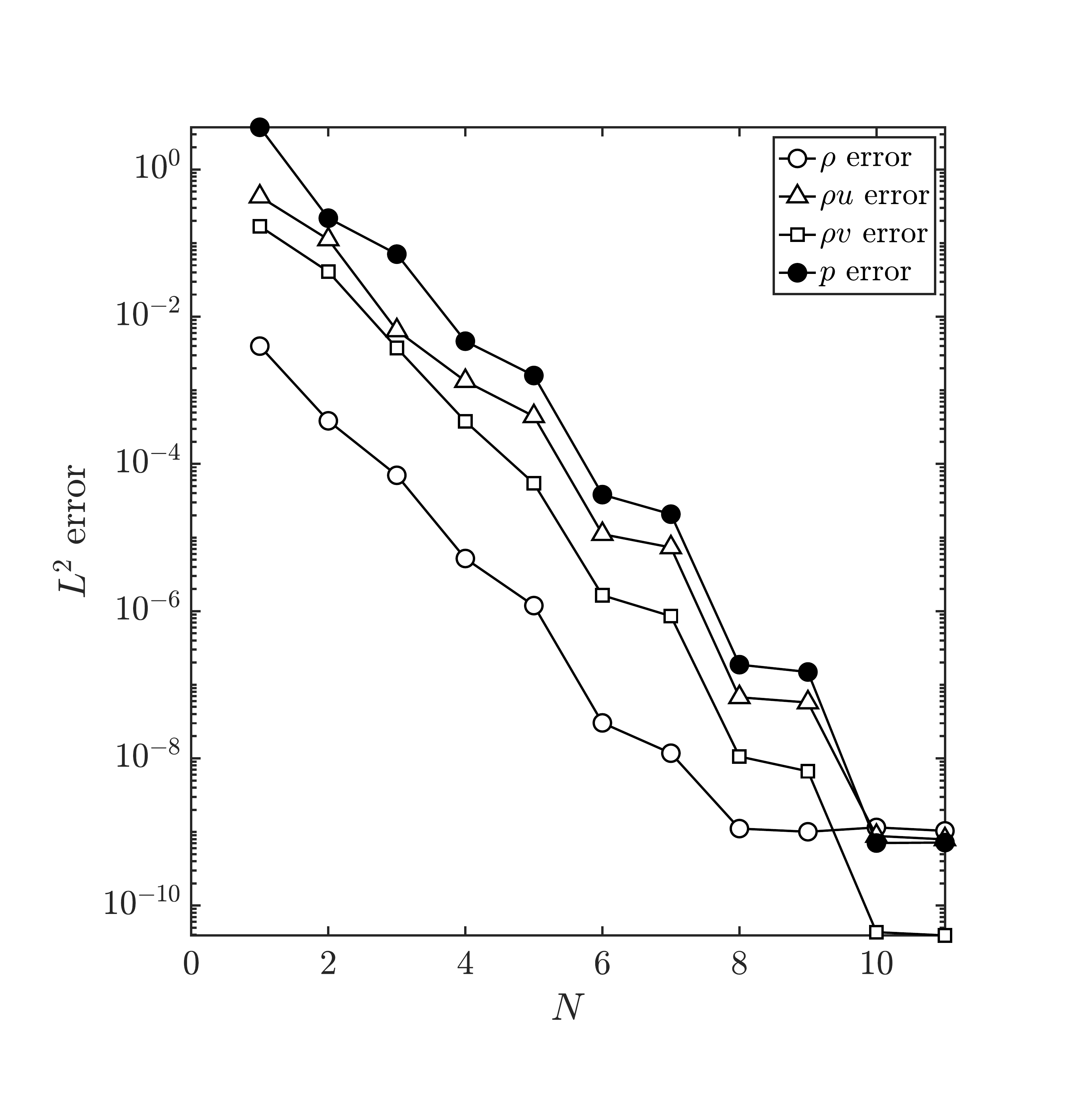}
  \caption{$L^2$ error of the Kovasznay test case solution with an $8^2$ mesh, for different polynomial orders. We use the split--form scheme with two averages,
  and the residual threshold in the time integration is $10^{-9}$}
  \label{fig:num:kov:p-Conv}
\end{figure}

We investigate the accuracy of the scheme on the Kovasznay two dimensional 
steady flow problem \cite{1948:Kovasznay}, with Reynolds number $\Rey=40$ on the domain 
$\Omega=[-0.5,1.5]\times[0,2]$. At the four boundaries (weakly through the exact Riemann solver), we apply the analytical solution,
\begin{equation}
  \begin{array}{l}
\displaystyle{  \rho_0 = 1, }\\
\displaystyle{  u_0 = 1-e^{\lambda x}\cos\left(2\pi y\right)},\\
\displaystyle{  v_0 = \frac{\lambda}{2\pi}e^{\lambda x}\sin\left(2\pi y\right)},\\
\displaystyle{  p_0 = \frac{1-e^{2\lambda x}}{2}},
  \end{array}
\end{equation}
where $\lambda$ is a parameter related to the Reynolds number 
$\lambda=\frac{\Rey}{2}-\sqrt{\frac{\Rey}{4}+4\pi^2}$. We use a uniform flow with uniform pressure as the 
initial condition. The mesh is an $8^2$ cartesian grid, and the polynomial order ranges from $N=2$ to $N=11$. 
We use the split--form scheme with two averages, $\mathrm{CFL}=0.75$, and we integrate in 
time for a residual threshold of $10^{-9}$. The $L^2$ error obtained is 
represented in Fig. \ref{fig:num:kov:p-Conv}. Although the solution converges to the residual threshold, the convergence rate is non--optimal. This is because the problem 
is viscous dominant, and the BR1 scheme suffers from even--odd behavior \cite{2017:Gassner}. As 
mentioned in \cite{2018:Manzanero,2019:Manzanero}, this can be solved by adding interface 
penalisation to the BR1 inter--element fluxes, which, for simplicity, we are not considering in 
this work. However, the even--odd effect is minimal when solving high Reynolds number 
flows.

\subsection{Inviscid Taylor--Green vortex}\label{subsec:num:TGV}

The high--order community has being using the three--dimensional inviscid Taylor--Green Vortex (TGV) as the 
reference problem to assess the 
robustness of the different methods to solve under--resolved transitional/turbulent flows. 
The configuration is a three--dimensional periodic box $[-1,1]^3$, with the initial condition,
\begin{equation}
  \begin{array}{l}
\displaystyle{  \rho_0 = 1, }\\
\displaystyle{  u_0 = \sin \pi x \cos \pi y \cos \pi z},\\
\displaystyle{  v_0 = -\cos \pi x \sin \pi y \cos \pi z},\\
\displaystyle{  w_0 = 0},\\
\displaystyle{  p_0 = \frac{1}{16}\left(\cos 2\pi x + \cos 2\pi y\right)\left(2 + \cos 
2\pi z\right)}.
  \end{array}
\end{equation}

Since the viscosity is zero, the mathematical entropy $\mathcal E$ 
should be constant according to the continuous bound \eqref{eq:stability-cont:balance-final}, when periodic boundary conditions are applied. The method is stable in the sense that the mathematical entropy $\mathcal E$ 
\eqref{eq:stability-cont:entropy-def} is bounded. 
%Since the mathematical entropy 
%is the sum of the kinetic energy $\mathcal K$ and the artificial compressibility 
%energy $\mathcal E_{AC}$, the kinetic energy rate is not monotonic, being these 
%errors proportional to the magnitude of the velocity divergence. 
%Taking the time 
%derivative in \eqref{eq:stability-cont:entropy-def},
%
%\begin{equation}
%\mathcal E_t = \mathcal K_t + M_0^2 p p_t = \mathcal K_t - p\svec{\nabla}\cdot\svec{u}  \leqslant 
%0.
%\end{equation}

The purpose of this problem is only to show that the method is entropy preserving in under--resolved
conditions. We construct a cartesian $8^3$ mesh, approximate the solution with $N=4$ order polynomials, and 
integrate in time until $t_F=20$ using the explicit RK3 scheme and 
$\mathrm{CFL}=0.75$.  Using the exact Riemann solver (ERS), we consider the split--form scheme (\rule[1pt]{10pt}{2pt}), the 
standard scheme ({\color{matlab_red}\rule[1pt]{10pt}{2pt}}), and the standard scheme with 
Gauss points ({\color{matlab_yellow}\rule[1pt]{10pt}{2pt}}). The first scheme 
should remain entropy stable ($\mathcal E_t\leqslant 0$), while the 
standard DG does not satisfy an entropy law. Additionally, we consider the 
split--form scheme using the two--point entropy conserving numerical flux (Central, {\color{matlab_blue}\rule[1pt]{10pt}{2pt}}), which is entropy 
conserving $\mathcal E_t=0$.

\begin{figure}[h]
  \centering
  \includegraphics[scale=0.3]{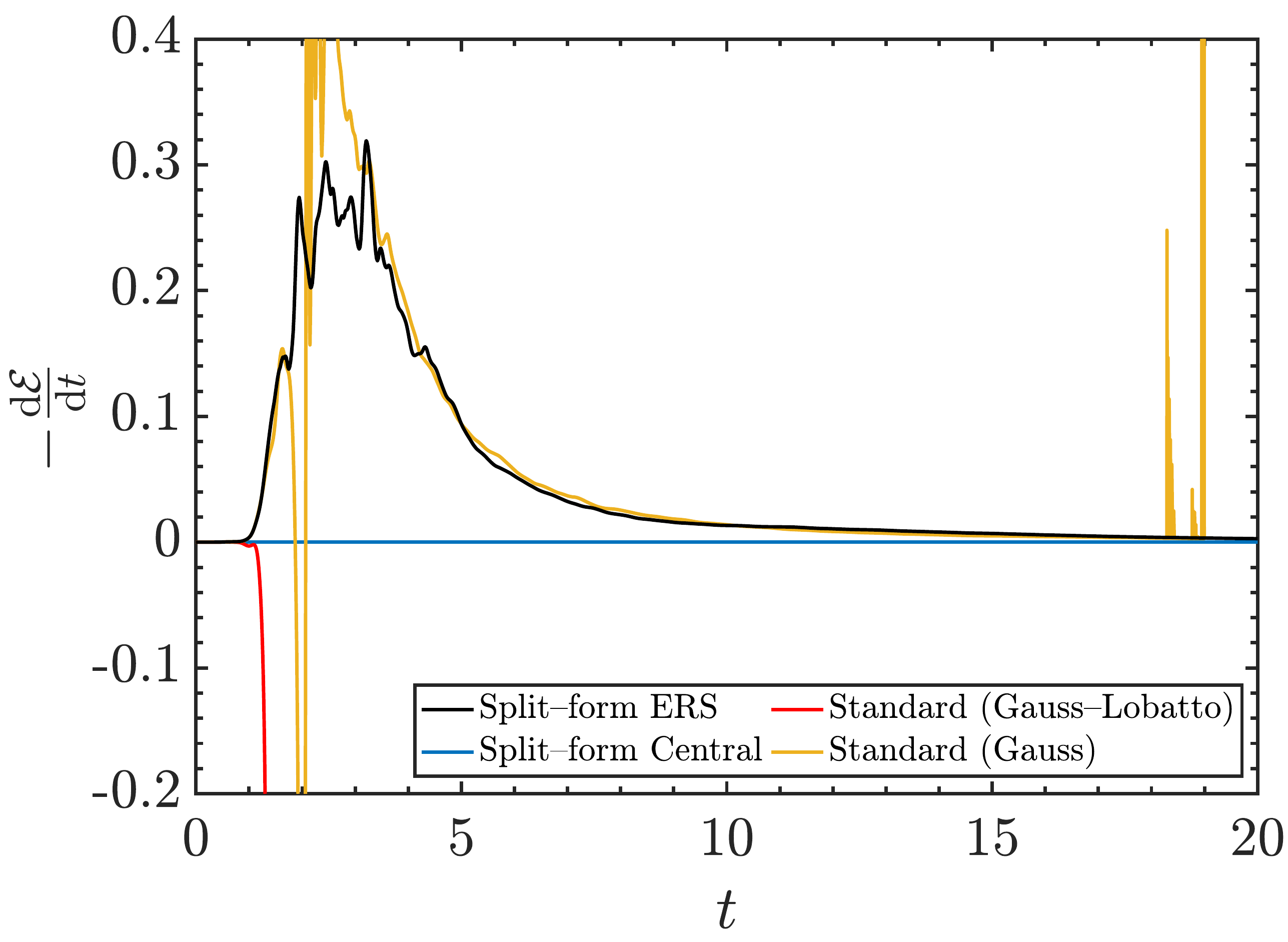}
  \caption{Evolution of the entropy time derivative for the four schemes considered.}
  \label{fig:num:TGV:entropy8P4}
\end{figure}

The evolution of the (negative) entropy time derivative is represented in Fig.~\ref{fig:num:TGV:entropy8P4} for the four 
schemes considered. On the one hand, the standard scheme is unstable and crashes for both 
Gauss ($t\approx 19.0$) and Gauss--Lobatto ($t\approx 1.4$). On the other hand, 
the split--form scheme is stable using both the ERS and the two--point entropy 
conserving numerical flux. while the former dissipates entropy at the inter--element 
boundaries, thus $\mathcal E_t \leqslant 0$, the latter is entropy conserving 
and $\mathcal E_{t}$ is machine precision order in all time instants. We are 
aware that to solve this problem, it is not enough to be entropy preserving, but to 
dissipate entropy at the appropriate rate. However, we confirm with these findings that using an entropy preserving 
scheme can be a suitable baseline scheme to which add additional artificial viscosity 
or LES models as subgrid--scale models 
\cite{2017:Flad,2017:Fernandez,2018:Manzanero-role}, which we do not cover here.

\subsection{Rayleigh--Taylor instability}\label{subsec:num:RTI}

To assess the robustness of the approximation in more challenging conditions, we solve 
the (surface tension free) two--dimensional Rayleigh--Taylor instability with Reynolds numbers $\Rey=1000$ and $\Rey=5000$,  
and compare our results to those in \cite{2017:Bassi} (where the authors used an incompressible solver with artificial 
viscosity). The initial condition is a heavy 
fluid ($\rho_{2}=3$) placed on top of a lighter one ($\rho_{1}=1$). We follow the same configuration used 
in \cite{2000:Guermond}, where the domain is a rectangle $\Omega = [0,0.5]\!\!\times\!\![0,4]$, 
which we discretize following \cite{2017:Bassi} using $16\!\!\times\!\!128$ equally--spaced elements 
with polynomial order $N=6$. The initial position of the interface is 
regularized following \cite{2000:Guermond},
\begin{equation}
  \rho = 2+\tanh\left(\frac{y-\eta(x)}{0.01}\right), ~~\eta(x) = 2-0.01\cos\left(2\pi 
  x\right).
\end{equation}
The boundary conditions are free--slip walls at left and right boundaries, and no--slip walls at bottom and top boundaries.
We set the parameter $M_0^2=2\cdot 10^{-4}$, the Froude number $\Fr=1$, we use the split--form scheme with
two averages, and a fixed timestep $\Delta t=1.5\cdot 
10^{-5}$ for the explicit RK3 scheme.

The approximation should remain stable even for 
under--resolved flows. Thus, we should not need to include artificial viscosity to stabilize it (although is desirable
to enhance the accuracy \cite{2018:Manzanero-role}). 
However, the analysis was performed assuming density positivity, which might not be the case if 
a sharp discontinuity in density is encountered. If the continuity equation is 
not regularized using a multiphase method (e.g. the Cahn--Hilliard equation 
\cite{1958:Cahn,2019:Manzanero}), the flow around the discontinuity suffers from Gibbs phenomena 
that eventually leads to negative density values, thus leading to a code crash. To tackle this problem without increasing the complexity of the 
underlying physics (i.e. without additional regularization of the density), we enforce a density limiter which is mass 
and momentum conserving. Leaving the state vector inaltered, we modify the 
fluxes and instead of dividing by the density to get the velocities (which can be close to zero and lead to divergence), we divide by the limited densities 
$\tilde{\rho}=f\left(\rho\right)$,
\begin{equation}
  \frac{\partial}{\partial t}\left\{\begin{array}{c}\rho \\ \rho u \\ \rho v \\ \rho w \\ p \end{array}\right\} 
  +\svec{\nabla} \cdot \left(\begin{array}{ccc}\rho u & \rho v & \rho w \\
  \frac{\left(\rho u\right)^2}{\tilde{\rho}} +p&   \frac{\rho u \cdot\rho v}{\tilde{\rho}} & \frac{\rho u \cdot\rho 
  w}{\tilde{\rho}}\\
  \frac{\rho u \cdot\rho v}{\tilde{\rho}}&   \frac{(\rho v)^2}{\tilde{\rho}} +p& \frac{\rho v \cdot\rho 
  w}{\tilde{\rho}}\\
 \frac{\rho u \cdot\rho 
  w}{\tilde{\rho}}&   \frac{\rho w \cdot\rho v}{\tilde{\rho}} & \frac{(\rho w)^2}{\tilde{\rho}} +p\\  
  \frac{\rho u}{M_0^2\tilde{\rho}} &   \frac{\rho v}{M_0^2\tilde{\rho}} &   \frac{\rho w}{M_0^2\tilde{\rho}}
  \end{array}\right) = \svec{\nabla}\cdot\ssvec{f}_{v}.
\end{equation}
The limited density $\tilde{\rho}$ can be any positive continuous function of the 
density. In this work, we use
\begin{equation}
  \tilde{\rho}=f(\rho)=\left\{\begin{array}{ccc}\rho & \text{if} & \rho \geqslant \rho_0\\
  \rho_0 & \text{if} & \rho < \rho_0 \end{array}\right..
\end{equation}
With this approach, the density acts as a working variable 
whose overshoots are controlled using the auxiliary value $\rho_0$  in troubled cells. In this 
work, we set $\rho_0=0.9$.

\begin{figure}[h]
  \centering
  \subfigure[$t=1.0$]{\includegraphics[height=10.0cm]{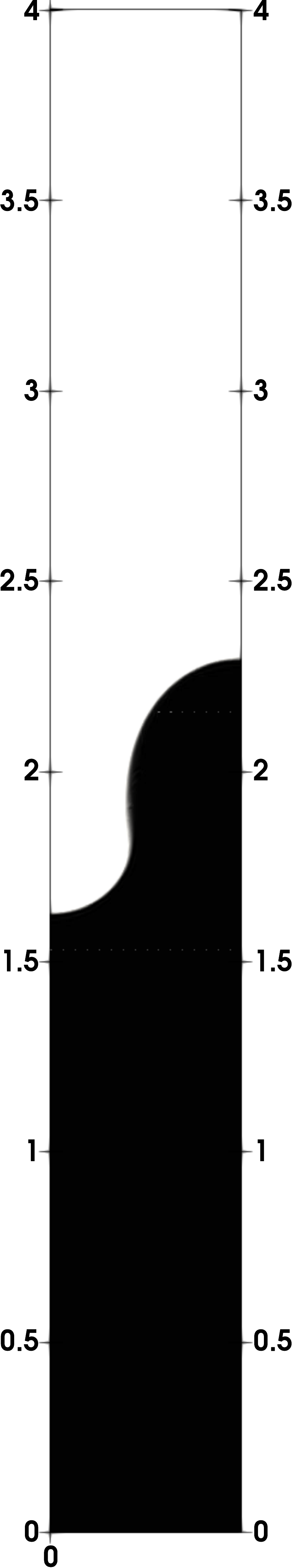}}  
  \subfigure[$t=1.5$]{\includegraphics[height=10.0cm]{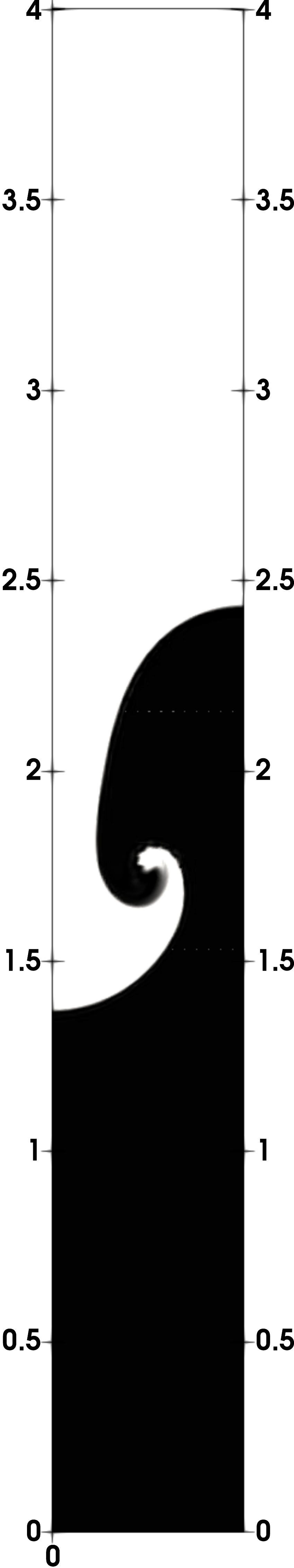}}  
  \subfigure[$t=1.75$]{\includegraphics[height=10.0cm]{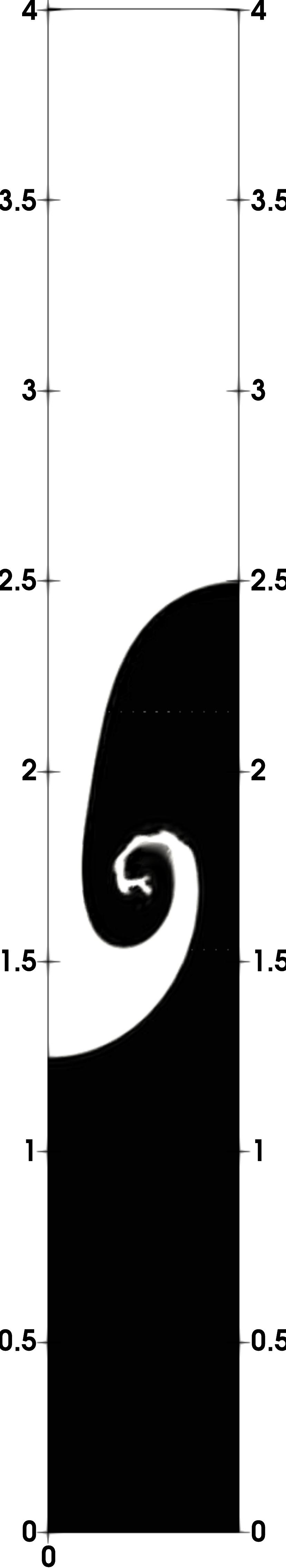}}  
  \subfigure[$t=2.0$]{\includegraphics[height=10.0cm]{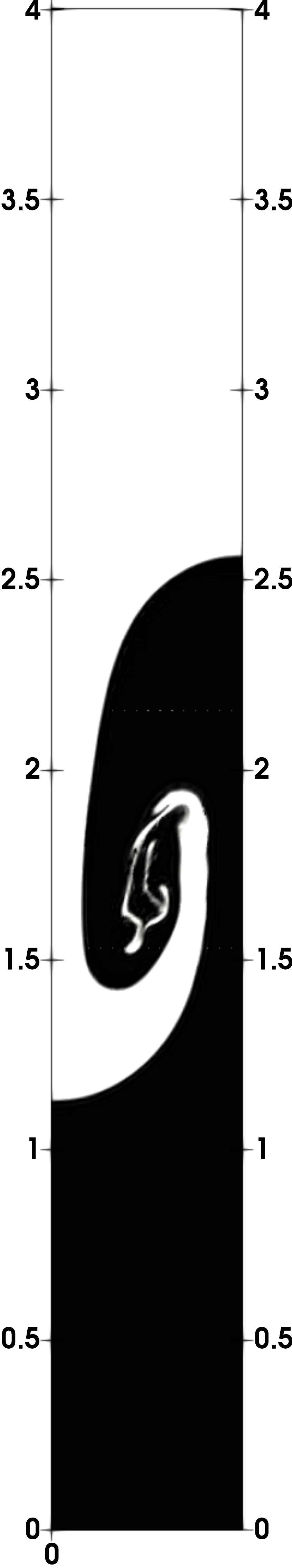}}  
  \subfigure[$t=2.25$]{\includegraphics[height=10.0cm]{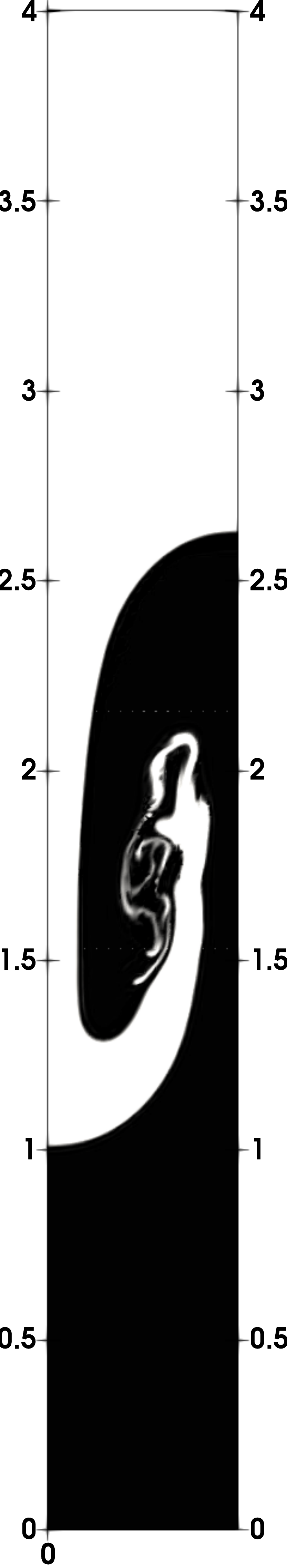}}  
  \subfigure[$t=2.5$]{\includegraphics[height=10.0cm]{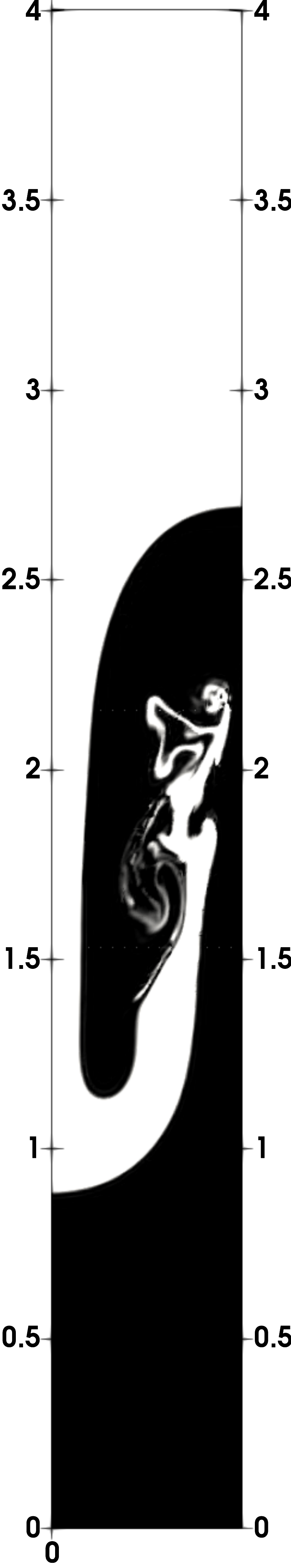}}  
  \caption{Rayleigh--Taylor instability with $\Rey=1000$: density contours. The 
  artificial compressibility Mach number is $M_0^2=2\cdot 10^{-4}$}
  \label{fig:num:RTI:density:Re1000}
\end{figure}

In Fig. \ref{fig:num:RTI:density:Re1000} we represent the density contours for $\Rey=1000$ at 
the same time snapshots shown in \cite{2017:Bassi} (in Trygvarsson's time scale 
$t=\sqrt{2}t_{Tryg}$). In this moderate Reynolds number simulation, we look the same 
as \cite{2017:Bassi}, although once it becomes under--resolved, 
the small structures patterns are different due to the artificial viscosity effect. 
Nonetheless, we find the differences subtle in this lower Reynolds number 
configuration.

\begin{figure}[h]
  \centering
  \subfigure[$t=1.0$]{\includegraphics[height=10.0cm]{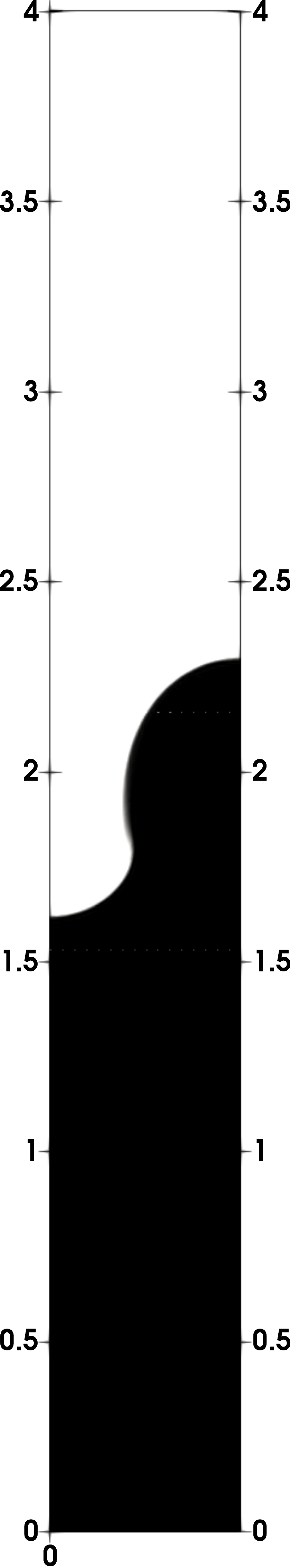}}  
  \subfigure[$t=1.5$]{\includegraphics[height=10.0cm]{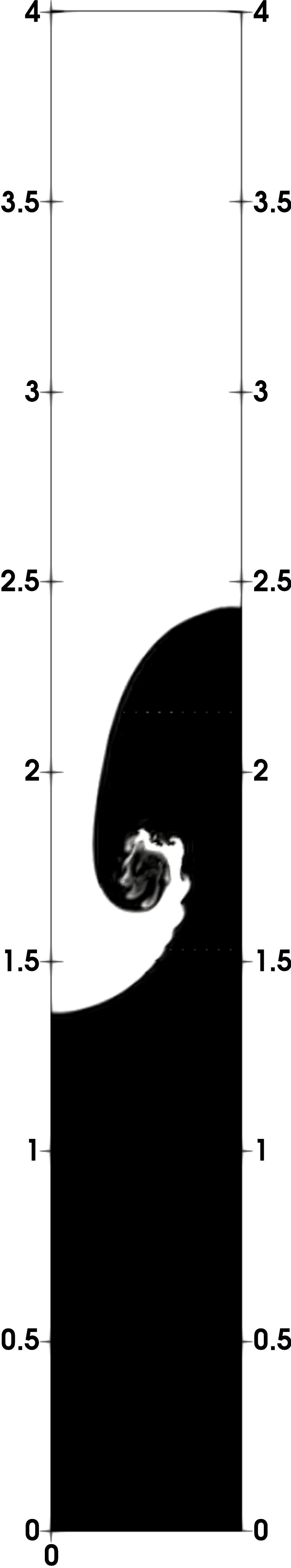}}  
  \subfigure[$t=1.75$]{\includegraphics[height=10.0cm]{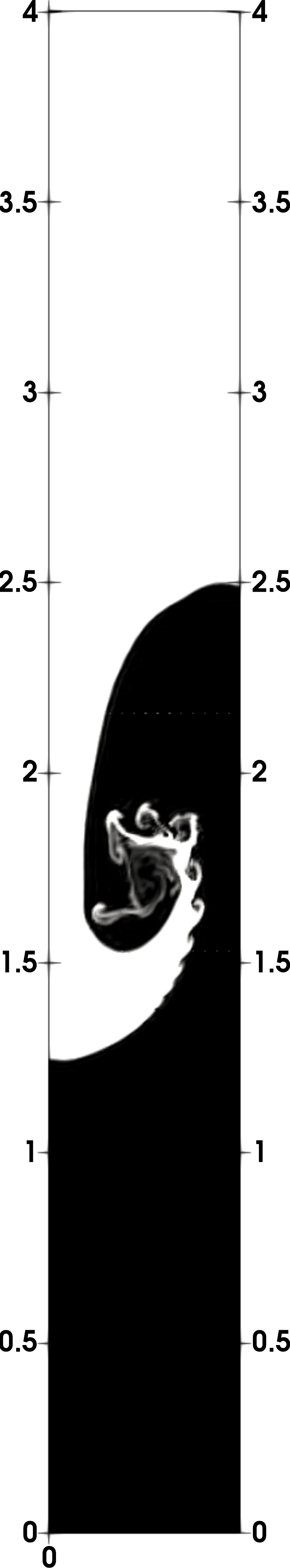}}  
  \subfigure[$t=2.0$]{\includegraphics[height=10.0cm]{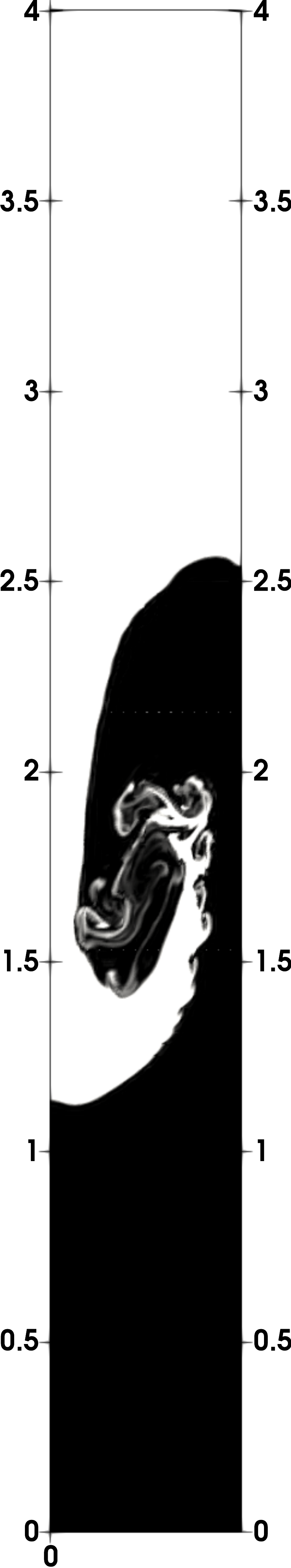}}  
  \subfigure[$t=2.25$]{\includegraphics[height=10.0cm]{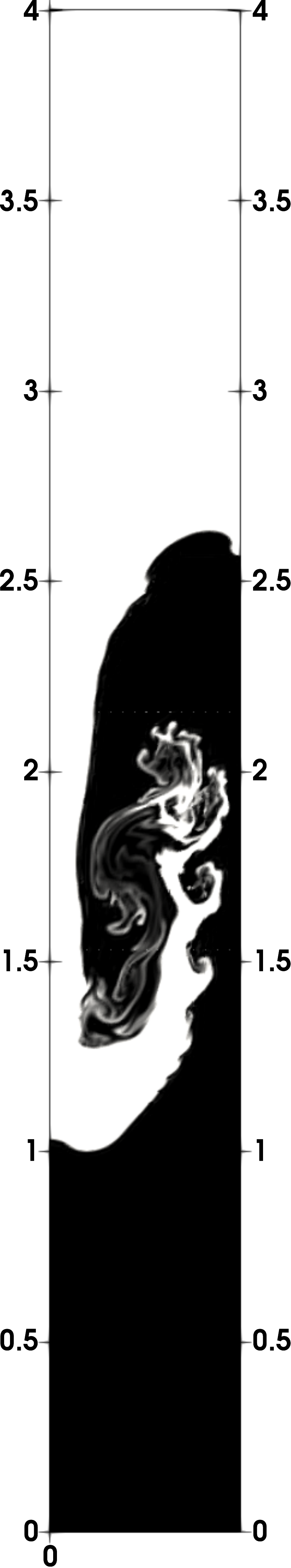}}  
  \subfigure[$t=2.5$]{\includegraphics[height=10.0cm]{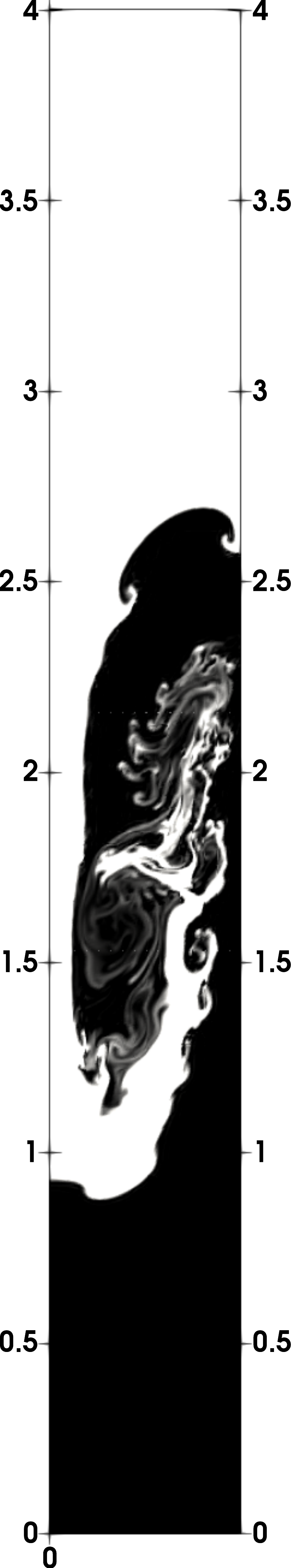}}  
  \caption{Rayleigh--Taylor instability with $\Rey=5000$: density contours. The 
  artificial compressibility Mach number is $M_0^2=2\cdot 10^{-4}$}
  \label{fig:num:RTI:density:Re5000}
\end{figure}

We represent the same density contours for $\Rey=5000$ in Fig. \ref{fig:num:RTI:density:Re5000}. 
The solution presented here looks the same as 
that presented in \cite{2017:Bassi}. However when the solution is 
under--resolved, $t>1.5$, the differences are not subtle. This is because we have 
not introduced any artificial dissipation, and thus we can represent smaller scales. By not introducing artificial dissipation, we confirm the robustness 
of the method.

\begin{figure}[h]
  \centering
  \includegraphics[scale=0.3]{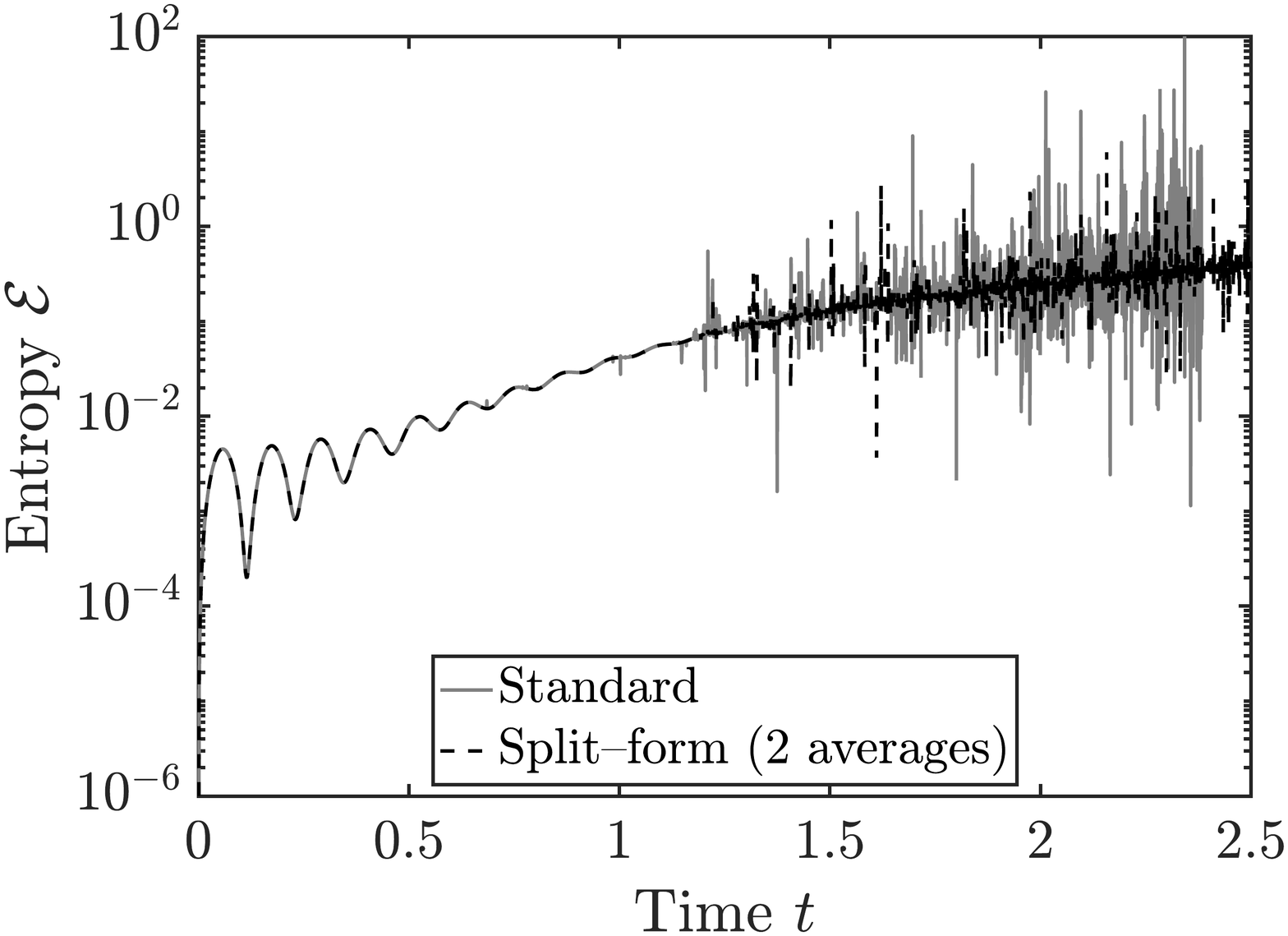}
  \caption{Rayleigh--Taylor instability with $\Rey=5000$: time evolution of the entropy \eqref{eq:stability-cont:entropy-def} for the split--form scheme 
  (with two averages in momentum) and the standard DG scheme. We find that not only the standard scheme
  is unstable and crashes, but also that the amplitude of the oscillations in the under--resolved stages are higher. Although entropy of the under-resolved split form solution is oscillatory, it is more robust}
  \label{fig:num:RTI:entropy}
\end{figure}

Finally, we solve the $\Rey=5000$ Rayleigh--Taylor instability also using the standard DG scheme, and 
compare the results in terms of the entropy evolution in Fig. 
\ref{fig:num:RTI:entropy}. First we note that the entropy increases due to gravitational force work. while the split--form is stable until the final time 
($t=2.5$), the standard scheme crashes close to $t=2.4$, when the flow gets 
severely under--resolved. Even in the laminar stage of the problem ($t<1.2$), 
the standard scheme starts to show minor noise. Not only is the standard scheme less 
robust, but the amplitude of the oscillations in the under--resolved stages 
is higher in the standard scheme, as a result of the uncontrolled creation of 
small structures without any stabilizing mechanism. 
%As a result, we find that using the 
%split--form scheme is the physically--consistent method to solve under--resolved 
%problems, since its robustness prevents to rest on artificial viscosity or other 
% techniques to be stable.

\section{Summary and conclusions}\label{sec:Conclusions}

We have developed an entropy stable DG approximation for the incompressible NSE 
with variable density and artificial compressibility. To do so, we first performed the continuous 
entropy analysis on a novel entropy function that includes both kinetic energy and artificial compressibility effects.
 We showed that the particular mathematical entropy for this set 
of equations is bounded in time, not including gravitational effects. Next, we constructed a DG scheme using Gauss--Lobatto points and the SBP--SAT 
property that mimics the continuous entropy analysis discretely. This was 
achieved using a two--point entropy conserving flux, with two options to 
discretize momentum, and the exact Riemann solver \cite{2017:Bassi} at both interior and physical boundaries. 
For viscous fluxes, we use the Bassi--Rebay 1 (BR1) scheme, obtaning as a result 
a parameter--free numerical discretization. The analysis was completed with the study of the solid wall 
boundary condition. 

 Lastly, we tested the numerical convergence of the scheme using manufactured solutions 
and the Kovasznay flow, and the robustness by solving the inviscid Taylor--Green Vortex (TGV) and the Rayleigh--Taylor
Instability (RTI) with $\Rey=1000$ and $\Rey=5000$. We show that the split--form 
scheme remains entropy stable (and entropy conserving if we use the two--point entropy conserving flux as the numerical flux) 
in severely under--resolved conditions, while the standard scheme can solve 
accurately the $\Rey=1000$ problem, but it is unstable and crashes for the inviscid TGV and RTI with $\Rey=5000$. 

We conclude that:

\begin{enumerate}
\item The DGSEM with two--point entropy conserving fluxes and wall boundary conditions is entropy stable with the exact Riemann solver.
\item Stability is reflected in robustness: the original scheme is unstable and crashes without the 
enhancements.
\item Numerical experiments show that the scheme is spectrally accurate, with exponential convergence for smooth flows.
\item Even when they were both stable, the entropy stable scheme was more 
accurate.
\item The splitting of the nonlinear terms in the two--point fluxes is not unique. Two splittings were studied. The one where the momentum is written as the product of two averages was slightly more accurate than using the average of the momentum. Both are entropy stable, however.
\end{enumerate}

\acknowledgement{The authors would like to thank Dr. Gustaaf Jacobs of the San Diego State University for his hospitality. This work was supported by a grant from the Simons Foundation ($\#426393$, David Kopriva). This project has received funding from the European Union’s Horizon 2020 research and innovation programme under grant agreement No 785549 (FireExtintion: H2020-CS2-CFP06-2017-01).
The authors acknowledge the computer resources and technical assistance provided by the Centro de Supercomputaci\'on y Visualizaci\'on de Madrid (CeSViMa).
}

\appendix

\section{Stability analysis of the exact Riemann solver}\label{appendix:IBT-exactRiemann}

Interior boundary terms stability rests on the positivity of \eqref{eq:stability:IBT-e}, which is satisfied if
\begin{equation}
\Delta_{e} = -\jump{\stvec{W}^{T}}\ssvec{F}_{e}^{\star}\cdot\svec{n}_L - \jump{\svec{F}^{\mathcal E}}\cdot\svec{n}_L + \jump{\stvec{W}^{T}\ssvec{F}_{e}}\cdot\svec{n}_L \geqslant 0,
\label{eq:IBT-stability:delta-e}
\end{equation}
at all face nodes. We use the rotational invariance \eqref{eq:Riemann:rotational-invariance} to transform \eqref{eq:IBT-stability:delta-e} to the face oriented system,
\begin{equation}
\begin{split}
\Delta_{e} =& -\jump{\stvec{W}^{T}}\smat{T}^{T}\stvec{F}_{e}\left(\stvec{Q}_n^{\star}\right) - \jump{{F}^{\mathcal E}\left(\stvec{Q}_{n}\right)} + \jump{\stvec{W}^{T}\smat{T}^{T}\stvec{F}_{e}\left(\stvec{Q}_{n}\right)} \\
&-\jump{\stvec{W_n}^{T}}\stvec{F}_{e}\left(\stvec{Q}_n^{\star}\right) - \jump{{F}^{\mathcal E}\left(\stvec{Q}_{n}\right)} + \jump{\stvec{W}^{T}_n\stvec{F}_{e}\left(\stvec{Q}_{n}\right)},
\end{split}
\label{eq:IBT-stability:delta-e-rot}
\end{equation}
where for the entropy flux 
\begin{equation}
\svec{F}^{\mathcal E}\cdot\svec{n}_L = \left(\frac{1}{2}\rho V_{tot}^2 + P\right)\svec{U}\cdot\svec{n}_L = \left(\frac{1}{2}\rho \left(U_n^2+V_{t1}^2 + V_{t2}^{2}\right) + P\right)U_n = F^{\mathcal E}\left(\stvec{Q}_{n}\right).
\end{equation}

Replacing the inviscid and entropy flux expressions in \eqref{eq:IBT-stability:delta-e-rot}, we obtain the condition to be satisfied by the exact Riemann problem solution,
\begin{equation}
\begin{split}
\Delta_{e} =& - \left(\jump{-\frac{1}{2}V_{tot}^2},\jump{U_n},\jump{V_{t1}},\jump{V_{t2}},M_0^2\jump{P}\right)\left(\begin{array}{c}\rho^{\star}U_n^{\star} \\
\rho^{\star}\left(U_n^{\star}\right)^2+P^{\star} \\
\rho^{\star}U_n^{\star}V_{t1}^{\star} \\
\rho^{\star}U_n^{\star}V_{t2}^{\star} \\
\frac{1}{M_0^{2}}U_n^{\star} 
\end{array}\right)- \jump{\left(\frac{1}{2}\rho V_{tot}^2 + P\right)U_n}+\jump{-\frac{1}{2}\rho V_{tot}^2U_n+\rho V_{tot}^2 U_n + 2PU_n} \\
=&\frac{1}{2}\rho^{\star}\jump{V_{tot}^2}U_n^{\star} - \rho^{\star}\left(U_n^{\star}\right)^2\jump{U_n} - P^{\star}\jump{U_n}-\rho^{\star}U_n^{\star}V_{t1}^{\star}\jump{V_{t1}}-\rho^{\star}U_n^{\star}V_{t2}^{\star}\jump{V_{t2}} - \jump{P}U_n^{\star} +\jump{P U_n}.
\end{split}
\end{equation}
Next, we use \eqref{eq:stability:average-jump-relation} to write,
\begin{equation}
\frac{1}{2}\jump{V_{tot}^2} = \aver{U_n}\jump{U_n} + \aver{V_{t1}}\jump{V_{t1}} + \aver{V_{t2}}\jump{V_{t2}},~~\jump{PU_n} = \aver{P}\jump{U_n}+\jump{P}\aver{U_n},
\end{equation}
which implies that
\begin{equation}
\begin{split}
\Delta_{e} =&\phantom{{}+{}}\left(\rho^{\star}\left(\aver{U_n}-U_n^{\star}\right)U_n^{\star}+\aver{P}-P^{\star}\right)\jump{U_n} + \left(\aver{U_n}-U_n^{\star}\right)\jump{P} \\
&+ \rho^{\star}\left(\aver{V_{t1}}-V_{t1}^{\star}\right)U_n^{\star}\jump{V_{t1}}+ \rho^{\star}\left(\aver{V_{t2}}-V_{t2}^{\star}\right)U_n^{\star}\jump{V_{t2}}.
\end{split}
\end{equation}
Now we replace the star region solution \eqref{eq:Riemann:star-region-solution}. To do so, we consider the case with $U_n^{\star}\geqslant 0$, where $\rho^{\star} = \rho^{\star}_L$ and $V_{ti}^{\star} = V_{tiL}$,
\begin{equation}
\begin{split}
\Delta_{e} =&\phantom{{}+{}}\left(\rho^{\star}_{L}\left(\aver{U_n}-U_n^{\star}\right)U_n^{\star}+\aver{P}-P^{\star}\right)\jump{U_n} + \left(\aver{U_n}-U_n^{\star}\right)\jump{P} \\
&+ \rho^{\star}_{L}\left(\aver{V_{t1}}-V_{t1L}\right)U_n^{\star}\jump{V_{t1}}+ \rho^{\star}_{L}\left(\aver{V_{t2}}-V_{t2L}\right)U_n^{\star}\jump{V_{t2}}.
\end{split}
\label{eq:IBT-stability:delta-e-ready}
\end{equation}
The last part, which involves tangential velocities, is stable since
\begin{equation}
\rho_L^{\star}\left(\frac{V_{tiL}+V_{tiR}}{2}-V_{tiL}\right)U_n^{\star}\jump{V_{ti}} = \frac{1}{2}\rho_L^{\star}U_n^{\star}\jump{V_{ti}}^{2}\geqslant 0.
\end{equation}

We write all pressures involved in \eqref{eq:IBT-stability:delta-e-ready} in terms of the velocities. First,
\begin{equation}
U_n^{\star} = - \frac{\jump{P} +\rho_R U_{nR}\lambda_{R}^{-} -\rho_L U_{nL}\lambda_L^{+}}{\rho_L\lambda_L^{+} - \rho_R\lambda_R^{-}},~~\jump{P} = -\Lambda U_n^{\star} - \rho_R U_{nR}\lambda_{R}^{-}+\rho_L U_{nL}\lambda_L^{+},
\end{equation}
where,
\begin{equation}
\Lambda = \rho_L \lambda_L^{+} - \rho_R \lambda_R^{-} > 0.
\end{equation}
Next, the averaged pressure minus the star region pressure is
\begin{equation}
\begin{split}
\aver{P}-P^{\star} =& \frac{P_L+P_R}{2} - P_L - \rho_L  \lambda_L^{+}\left(U_{nL}-U_n^{\star}\right) = \frac{1}{2}\jump{P} - \rho_L  \lambda_L^{+}\left(U_{nL}-U_n^{\star}\right)\\
=&-\frac{1}{2}\Lambda U_n^{\star} - \frac{1}{2}\rho_R U_{nR}\lambda_{R}^{-}+\frac{1}{2}\rho_L U_{nL}\lambda_L^{+}- \rho_L  \lambda_L^{+}\left(U_{nL}-U_n^{\star}\right) \\
=&-\frac{1}{2}\Lambda U_n^{\star} - \frac{1}{2}\rho_R U_{nR}\lambda_{R}^{-} - \frac{1}{2}\rho_L\lambda_L^{+}U_{nL} + \rho_L \lambda_L^{+}U_n^{\star}.
\end{split}
\end{equation}
which replaced in \eqref{eq:IBT-stability:delta-e-ready}, and replacing the star region density $\rho_L^{\star}$ from \eqref{eq:Riemann:star-region-solution} gives,
\begin{equation}
\begin{split}
\Delta_{e} =&\left(\rho^{\star}_{L}\left(\frac{U_{nL}+U_{nR}}{2}-U_n^{\star}\right)U_n^{\star}-\frac{1}{2}\Lambda U_n^{\star} - \frac{1}{2}\rho_R U_{nR}\lambda_{R}^{-} - \frac{1}{2}\rho_L\lambda_L^{+}U_{nL} + \rho_L \lambda_L^{+}U_n^{\star}\right)\left(U_{nR}-U_{nL}\right) \\
&+ \left(\frac{U_{nL}+U_{nR}}{2}-U_n^{\star}\right)\left(-\Lambda U_n^{\star} - \rho_R U_{nR}\lambda_{R}^{-}+\rho_L U_{nL}\lambda_L^{+}\right) + \frac{1}{2}\rho_L^{\star}U_n^{\star}\left(\jump{V_{t1}}^{2}+\jump{V_{t2}}^{2}\right) \\
=&\frac{1}{U_n^{\star}-\lambda_L^{-}}\biggl\{\left[\frac{1}{2}\rho_L\lambda_L^{+}\left(U_{nL}+U_{nR}-2U_n^{\star}\right)U_n^{\star}-\left(\frac{1}{2}\Lambda U_n^{\star} + \frac{1}{2}\rho_R U_{nR}\lambda_{R}^{-} + \frac{1}{2}\rho_L\lambda_L^{+}U_{nL} - \rho_L \lambda_L^{+}U_n^{\star}\right)\left(U_n^{\star}-\lambda_L^{-}\right)\right]\left(U_{nR}-U_{nL}\right)\\
&+ \frac{1}{2}\left(U_{nL}+U_{nR}-2U_n^{\star}\right)\left(-\Lambda U_n^{\star} - \rho_R U_{nR}\lambda_{R}^{-}+\rho_L U_{nL}\lambda_L^{+}\right)\left(U_n^{\star}-\lambda_L^{-}\right) + \frac{1}{2}\rho_L\lambda_L^{+}U_n^{\star}\left(\jump{V_{t1}}^{2}+\jump{V_{t2}}^{2}\right)\biggr\}.
\end{split}
\end{equation}
Next, we define $\hat{U}_L=U_{nL}-U_n^{\star}$ and $\hat{U}_R=U_{nR}-U_n^{\star}$ and use them to write $\Delta_{e}$,
\begin{equation}
\begin{split}
\Delta_{e}=&\frac{1}{U_n^{\star}-\lambda_L^{-}}\biggl\{\left[\frac{1}{2}\rho_L\lambda_L^{+}\left(\hat{U}_L+\hat{U}_{R}\right)U_n^{\star}-\left(\frac{1}{2}\Lambda U_n^{\star} + \frac{1}{2}\rho_R \left(\hat{U}_{R}+U_n^{\star}\right)\lambda_{R}^{-} + \frac{1}{2}\rho_L\lambda_L^{+}\left(\hat{U}_{L}+U_n^{\star}\right) - \rho_L \lambda_L^{+}U_n^{\star}\right)\left(U_n^{\star}-\lambda_L^{-}\right)\right]\left(\hat{U}_{R}-\hat{U}_L\right)\\
&+ \frac{1}{2}\left(\hat{U}_L+\hat{U}_{R}\right)\left(-\Lambda U_n^{\star} - \rho_R  \left(\hat{U}_{R}+U_n^{\star}\right)\lambda_{R}^{-}+\rho_L  \left(\hat{U}_{L}+U_n^{\star}\right)\lambda_L^{+}\right)\left(U_n^{\star}-\lambda_L^{-}\right) + \frac{1}{2}\rho_L\lambda_L^{+}U_n^{\star}\left(\jump{V_{t1}}^{2}+\jump{V_{t2}}^{2}\right)\biggr\}\\
=&\frac{1}{U_n^{\star}-\lambda_L^{-}}\biggl\{\left[\frac{1}{2}\rho_L\lambda_L^{+}\left(\hat{U}_L+\hat{U}_{R}\right)U_n^{\star}-\frac{1}{2}\left(\rho_R \hat{U}_{R}\lambda_{R}^{-} +\rho_L\lambda_L^{+}\hat{U}_{L} \right)\left(U_n^{\star}-\lambda_L^{-}\right)\right]\left(\hat{U}_{R}-\hat{U}_L\right)\\
&+ \frac{1}{2}\left(\hat{U}_L+\hat{U}_{R}\right)\left(\rho_L \hat{U}_{L}\lambda_L^{+} - \rho_R\hat{U}_{R}\lambda_{R}^{-}\right)\left(U_n^{\star}-\lambda_L^{-}\right) + \frac{1}{2}\rho_L\lambda_L^{+}U_n^{\star}\left(\jump{V_{t1}}^{2}+\jump{V_{t2}}^{2}\right)\biggr\} \\
=&\frac{1}{U_n^{\star}-\lambda_L^{-}}\biggl\{\frac{1}{2}\rho_L\lambda_L^{+}\left(\hat{U}_{R}^2-\hat{U}_L^2\right)U_n^{\star}+\left(\rho_L\lambda_L^{+}\hat{U}_{L}^{2} - \rho_R\lambda_R^{-}\hat{U}_R^{2}\right)\left(U_n^\star-\lambda_{L}^{-}\right)+ \frac{1}{2}\rho_L\lambda_L^{+}U_n^{\star}\left(\jump{V_{t1}}^{2}+\jump{V_{t2}}^{2}\right)\biggr\} \\
=&\frac{1}{U_n^{\star}-\lambda_L^{-}}\left(\frac{1}{2}\rho_L\lambda_L^{+}\left(U_n^{\star}-2\lambda_L^{-}\right)\hat{U}_{L}^{2}+\left(\frac{1}{2}\rho_L\lambda_L^{+}U_n^{\star}-\rho_R\lambda_R^{-}\left(U_n^{\star}-\lambda_L^{-}\right)\right)\hat{U}_{R}^{2}+ \frac{1}{2}\rho_L\lambda_L^{+}U_n^{\star}\left(\jump{V_{t1}}^{2}+\jump{V_{t2}}^{2}\right)\right)\geqslant 0.
\end{split}
\end{equation}
Therefore, 
\begin{equation}
\Delta_{e}=\frac{1}{U_n^{\star}-\lambda_L^{-}}\left(\frac{1}{2}\rho_L\lambda_L^{+}\left(U_n^{\star}-2\lambda_L^{-}\right)\hat{U}_{L}^{2}+\left(\frac{1}{2}\rho_L\lambda_L^{+}U_n^{\star}-\rho_R\lambda_R^{-}\left(U_n^{\star}-\lambda_L^{-}\right)\right)\hat{U}_{R}^{2}+ \frac{1}{2}\rho_L\lambda_L^{+}U_n^{\star}\left(\jump{V_{t1}}^{2}+\jump{V_{t2}}^{2}\right)\right)\geqslant 0,
\end{equation}
which confirms that $\Delta_{e}$ is always positive since $U_n^{\star}\geqslant 0$, $\lambda_L^{+}>0$, $\lambda_L^{-}<0$, and $\lambda_R^{-}<0$. 

The magnitude of the dissipation introduced depends on the square of the tangential speed jumps and the square of the jumps between the star solution $U_n^{\star}$ and the left and right states,
\begin{equation}
\begin{split}
\Delta_{e}=\frac{1}{U_n^{\star}-\lambda_L^{-}}\biggl(&\frac{1}{2}\rho_L\lambda_L^{+}\left(U_n^{\star}-2\lambda_L^{-}\right)\left(U_{nL}-U_n^{\star}\right)^{2}+\left(\frac{1}{2}\rho_L\lambda_L^{+}U_n^{\star}-\rho_R\lambda_R^{-}\left(U_n^{\star}-\lambda_L^{-}\right)\right)\left(U_{nR}-U_n^{\star}\right)^{2}\\
+& \frac{1}{2}\rho_L\lambda_L^{+}U_n^{\star}\left(\jump{V_{t1}}^{2}+\jump{V_{t2}}^{2}\right)\biggr)\geqslant 0.
\end{split}
\end{equation}

In the other possible case, $U_n^{\star}<0$, the exact Riemann solver is also dissipative and
\begin{equation}
\begin{split}
\Delta_{e}=\frac{1}{U_n^{\star}-\lambda_R^{+}}\biggl(&-\frac{1}{2}\rho_R\lambda_R^{-}\left(U_n^{\star}-2\lambda_R^{+}\right)\left(U_{nR}-U_n^{\star}\right)^{2}+\left(-\frac{1}{2}\rho_R\lambda_R^{-}U_n^{\star}+\rho_L\lambda_L^{+}\left(U_n^{\star}-\lambda_R^{+}\right)\right)\left(U_{nL}-U_n^{\star}\right)^{2}\\
-& \frac{1}{2}\rho_R\lambda_R^{-}U_n^{\star}\left(\jump{V_{t1}}^{2}+\jump{V_{t2}}^{2}\right)\biggr)\geqslant 0.
\end{split}
\end{equation}
Therefore, the the positivity of \eqref{eq:stability:IBT-e} in the split--form DGSEM with the exact Riemann solver derived in \cite{2017:Bassi} is 
established.

If one uses the two--point entropy flux as the Riemann solver, 
$\ssvec{F}_{e}^{\star}=\ssvec{F}^{ec}$, obtains neutral stability, 
$\Delta_{e}=0$, as a result of Tadmor's jump condition 
\eqref{eq:stability:two-point-entropy-contraction}. Neutral stability usually 
leads to undesirable, oscillating solutions and even--odd error behaviours (see 
\cite{2016:gassner}). A systematical procedure to add controlled dissipation is by 
augmenting the two--point entropy flux at the boundaries with a matrix 
dissipation term. However, we have not explored how the latter should be performed for 
the set of equations studied.

\section{Equivalence of the two--point flux to a split form equation}\label{appendix:SplitFormAppendix}
The divergence at a point is approximated by the two--point formula \eqref{eq:dg:split-form-op}, which we repeat here making the usual definition of the derivative matrix, $D_{im} \equiv l'_{m}(\xi_{i})$, etc,
\begin{equation}
\begin{split}
\mathbb D\left(\cssvec{F}_{e}\right)^{\#}_{ijk} = 2\sum_{m=0}^{N}&\phantom{{}+{}}D_{im}\ssvec{F}_{e}^{\#}\left(Q_{ijk},Q_{mjk}\right)\cdot\aver{\mathcal J\svec{a}^{1}}_{(im)jk}\\
&+D_{jm}\ssvec{F}_{e}^{\#}\left(Q_{ijk},Q_{imk}\right)\cdot\aver{\mathcal J\svec{a}^{2}}_{i(jm)k} \\
&+D_{km}\ssvec{F}_{e}^{\#}\left(Q_{ijk},Q_{ijm}\right)\cdot\aver{\mathcal 
J\svec{a}^{3}}_{ij(km)}.
\end{split}
\label{eq:dg:split-form-opB}
\end{equation}
Note that the approximation contains products of the two--point ($i$ and $m$) averages of the fluxes and metric terms. The fluxes themselves, as seen in \eqref{eq:dg:two-point-flux-iNS}, also contain products of two--point averages. 

Two--point approximations are algebraically equivalent to the approximation of the divergence written in a more recognizable split form that contains averages of conservative and nonconservative forms, C.F. \cite{2016:gassner}.

To show the equivalent PDEs approximations, we first show the equivalence for a split flux that depends only on a single average vector,
\begin{equation}
\svec{F}^{\# 1} = \left( \average{A_{1}}\;\average{A_{2}}\;\average{A_{3}} \right)
\end{equation}
for a (space) vector flux
\begin{equation}
\svec{F} = \left( {A_{1}}\;{A_{2}}\;{A_{3}} \right),
\end{equation}
whose components are polynomials of the reference space variables.
Then \eqref{eq:dg:split-form-opB} becomes
\begin{equation}
\begin{split}
\mathbb D\left(\svec{F}\right)^{\# 1}_{ijk} &= 2\sum_{m=0}^{N}\left\{ D_{im}\average{A_{1}}_{(im)jk} + D_{jm}\average{A_{2}}_{i(jm)k} +D_{km}\average{A_{3}}_{ij(km)} \right\}
\\&
=\sum_{m=0}^{N}\left\{ D_{im}(A_{1})_{mjk} + D_{jm}(A_{2})_{imk} +D_{km}(A_{3})_{ijm} \right\} 
\\& +\sum_{m=0}^{N}\left\{ D_{im}(A_{1})_{ijk} + D_{jm}(A_{2})_{ijk} +D_{km}(A_{3})_{ijk} \right\}
\end{split}
\label{eq:SimpleAverageDerivative}
\end{equation}
Since the derivative of a constant function is zero and $\sum_{m=0}^{N} D_{im}=0$, the second sum in \eqref{eq:SimpleAverageDerivative} vanishes so that
\begin{equation}
\mathbb D\left(\svec{F}\right)^{\# 1}_{ijk} = \sum_{m=0}^{N}\left\{ D_{im}(A_{1})_{mjk} + D_{jm}(A_{2})_{imk} +D_{km}(A_{3})_{ijm} \right\} = \left(\nabla_{\xi}\cdot\svec{F}\right)_{ijk}
\label{eq:SingleAverageEquiv}
\end{equation}
Thus, the divergence approximation using simple two--point average is equivalent to the divergence of the original flux.

Next in complexity is if the flux
is the product of two quantities
\begin{equation}
\svec{F} = \left( {A_{1}}B\;{A_{2}}B\;{A_{3}} B\right).
\end{equation}
 and the associated two--point flux replaces the product with the product of two averages
\begin{equation}
\svec{F}^{\# 2} = \left( \average{A_{1}}\average{B}\;\average{A_{2}}\average{B}\;\average{A_{3}}\average{B} \right).
\end{equation}
Since
\begin{equation}
\begin{split}
4&\sum_{m=0}^{N}D_{im}\average{A_{1}}\average{B} \\&= \sum_{m=0}^{N}D_{im}\left( \left(A_{1}\right)_{ijk}+\left(A_{1}\right)_{mjk}\right)\left(B_{ijk}+B_{mjk}\right) \\&
= \sum_{m=0}^{N}D_{im}\left(A_{1}\right)_{mjk}B_{mjk}+ \left(A_{1}\right)_{ijk}\sum_{m=0}^{N}D_{im}B_{mjk}+ B_{ijk}\sum_{m=0}^{N}D_{im}\left(A_{1}\right)_{mjk},
\end{split}
\end{equation}
it follows that
\begin{equation}
\mathbb D\left(\svec{F}\right)^{\# 2}_{ijk} = \frac{1}{2}\left(\nabla_{\xi}\cdot\svec{F}\right)_{ijk}+ \frac{1}{2}\left( \svec A \cdot \nabla_{\xi} B +  B  \nabla_{\xi}\cdot \svec A \right)_{ijk},
\label{eq:DoubleAverageEquiv}
\end{equation}
which is the average of the conservative form of the divergence and the product rule applied to it. 

Finally, we follow the same steps to find the equivalent approximation for a triple product flux of the form
\begin{equation}
\svec{F} = \left( C{A_{1}}B\;C{A_{2}}B\;C{A_{3}} B\right).
\end{equation}
approximated by
\begin{equation}
\svec{F}^{\# 3} = \left( \average{C}\average{A_{1}}\average{B}\;\average{C}\average{A_{2}}\average{B}\;\average{C}\average{A_{3}}\average{B} \right).
\end{equation}
The result is a specific application \cite{2008:Kennedy} of the product rule to the product of three polynomials, 
\begin{equation}
\begin{split}
\mathbb D\left(\svec{F}\right)^{\# 3}_{ijk} = \frac{1}{4}\left(\nabla_{\xi}\cdot\svec{F}\right)_{ijk} 
&+ \frac{1}{4}\left\{ (C\svec A)_{ijk}\cdot\nabla_{\xi}B_{ijk} 
+ (CB)_{ijk}\nabla_{\xi}\cdot \svec A 
+ C_{ijk}\nabla_{\xi}\cdot\mathbb I^{N}\left(\svec A B\right)_{ijk}  \right. 
\\& \left.
+ (B\svec A)_{ijk}\cdot \nabla_{\xi}C_{ijk} 
+ \svec A_{ijk}\cdot\nabla_{\xi}\mathbb I^{N}(CB)_{ijk} 
+ B_{ijk}\nabla_{\xi}\cdot\mathbb I^{N}(C\svec A)_{ijk} \right\}.
\end{split}
\label{eq:TripleAverageEquiv}
\end{equation}

The second form, \eqref{eq:DoubleAverageEquiv}, can be used to assess the influence of the metric terms in \eqref{eq:dg:split-form-opB}. Replacing $\svec A$ by the matrix $\tens M$ whose columns are $\mathcal J\svec a^{i}, i=1,2,3$, the \emph{contravariant vector} form of the divergence becomes
\begin{equation}
\mathbb D\left(\svec{\tilde F}\right)^{\# 1}_{ijk} =\frac{1}{2}\left(\svec{\nabla}_{\xi}\cdot\svec{\tilde{F}}\right)_{ijk} + \frac{1}{2}\left(\mathcal J \svec{a}^{1}\cdot\frac{\partial\svec{\tilde{F}}}{\partial\xi}+\mathcal J \svec{a}^{2}\cdot\frac{\partial\svec{\tilde{F}}}{\partial\eta}+\mathcal J \svec{a}^{3}\cdot\frac{\partial\svec{\tilde{F}}}{\partial\zeta}\right)_{ijk}
\label{eq:equiv-split:metrics-dealiasing}
\end{equation}
by virtue of the discrete metric identities
\begin{equation}
  \sum_{m=0}^{N}\left(D_{im}\mathcal J\svec{a}^{1}_{mjk}+D_{jm}\mathcal J\svec{a}^{2}_{imk}+D_{km}\mathcal 
J\svec{a}^{3}_{ijm}\right) = \left(\svec{\nabla}_{\xi}\cdot\tens{M}\right)_{ijk}=\left(\sum_{n=1}^{3}\frac{\partial \mathcal J\svec{a}^{n}}{\partial 
\xi^{n}}\right)_{ijk}=0.
\end{equation}
In light of \eqref{eq:equiv-split:metrics-dealiasing}, we simplify the discussion below and examine the approximation \eqref{eq:dg:split-form-opB} using the pointwise values of the metric terms rather than the averages.

From these three forms \eqref{eq:SingleAverageEquiv}, \eqref{eq:DoubleAverageEquiv} and \eqref{eq:TripleAverageEquiv}, we can derive the split form approximation equivalent to the two--point fluxes given by \eqref{eq:dg:two-point-flux-iNS} with the two choices for the momentum approximation, \eqref{eq:stability:momentum-skewsymmetric-options}. Notice that with the first momentum average, the entropy conserving fluxes are single or product averages. Using the second momentum average the flux is made up of double or triple averages.\\

 \noindent{\it Option 1:  $\tilde{\rho u_i}^{(1)}=\aver{\rho u_i}$}\\

When the momentum is implemented by the average, the continuity equation is in single average form and hence by \eqref{eq:SingleAverageEquiv},
\begin{equation}
\svec{\nabla}\cdot\left(\rho \svec{ u}\right)_{ijk} \approx \svec{\nabla}_\xi\cdot\mathbb I^{N}\left(\rho \svec{\tilde U}\right)_{ijk},
\end{equation}
where $\svec{\tilde u}=\tens{M}^{T}\svec{u}$ is the contravariant velocity.

The momentum equation is approximated using the product of two averages for the momentum and  the average for the pressure so that it includes forms \eqref{eq:DoubleAverageEquiv} and \eqref{eq:SingleAverageEquiv}. The two--point flux form therefore approximates

\begin{equation}
\begin{split}
\svec{\nabla}_{\xi}\cdot\left(\rho \svec{\tilde u}u_l + p \tens{M}^{T}\svec{e}_{l}\right)\approx& \frac{1}{2}\svec{\nabla}_{\xi}\cdot\mathbb I^{N}\left(\rho U_l \svec{\tilde U}\right)_{ijk} + \frac{1}{2}U_{l,ijk}\svec{\nabla}_{\xi}\cdot\mathbb I^{N}\left(\rho \svec{\tilde U}\right)_{ijk} \\
&+ \frac{1}{2}\left(\rho \svec{U}\right)_{ijk}\cdot\left(\svec{\nabla}_{\xi}\cdot\mathbb I^{N}\left(U_l\tens{M}\right)\right)_{ijk} + \svec{\nabla}_{\xi}\cdot\mathbb I^{N}\left(P\tens{M}^{T}\svec{e}_{l}\right)_{ijk}.
\end{split}
\end{equation}
Finally, the two--point flux form of the artificial compressibility equation, like the continuity equation, is equivalent to using standard DG since it is linear in the velocities.\\

\noindent{\it Option 2: $\tilde{\rho u_i}^{(2)}=\aver{\rho}\aver{u_i}$}\\

Approximating the momentum as the product of two averages (vs. the average of the product) leads to the approximation of a different form of the equations by increasing the number of products in each equation.

Under the second approximation, the continuity equation now has the product of two averages and hence is the equivalent to the approximation
\begin{equation}
\svec{\nabla}_\xi\cdot\left(\rho \svec{\tilde u}\right)_{ijk} \approx  \frac{1}{2}\svec{\nabla}_{\xi}\cdot\mathbb I^{N}\left(\rho\svec{\tilde U}\right)_{ijk} + \frac{1}{2}\rho_{ijk}\svec{\nabla}_{\xi}\cdot\left(\svec{\tilde U}\right)_{ijk}+ \frac{1}{2}\left(\svec{U}\right)_{ijk}\cdot\left(\svec{\nabla}_{\xi}\cdot\mathbb I^{N} \left(\rho\tens{M}\right)\right)_{ijk}.
\end{equation}
The momentum equation is approximated with the triple product for the momentum, and still the standard approximation in the pressure,
\begin{equation}
\begin{split}
\svec{\nabla}_{\xi}&\cdot\left(\rho \svec{\tilde u}u_l + p \tens{M}^{T}\svec{e}_{l}\right)\approx\frac{1}{4}\svec{\nabla}_{\xi}\cdot\mathbb I^{N}\left(\rho U_{l} \svec{\tilde{U}}\right)_{ijk} + \frac{1}{4}\rho_{ijk}U_{l,ijk}\left(\svec{\nabla}_{\xi}\cdot\svec{\tilde{U}}\right)_{ijk} + \frac{1}{4}\left(\rho \svec{\tilde U}\right)_{ijk}\cdot\left(\svec{\nabla}_{\xi}\cdot\mathbb I^{N}\left(U_l\tens{M}\right)\right)_{ijk} \\
&+\frac{1}{4}\rho_{ijk}\svec{\nabla}_{\xi}\cdot\mathbb I^{N}\left(U_l\svec{\tilde{U}}\right)_{ijk} + \frac{1}{4}\left(u_l \svec{\tilde U}\right)_{ijk}\cdot\left(\svec{\nabla}_{\xi}\cdot\mathbb I^{N}\left(\rho\tens{M}\right)\right)_{ijk}+\frac{1}{4}U_{l,ijk}\svec{\nabla}_{\xi}\cdot\mathbb I^{N}\left(\rho\svec{\tilde{U}}\right)_{ijk}\\
&+ \frac{1}{4}\left(\svec{\tilde U}\right)_{ijk}\cdot\left(\svec{\nabla}_{\xi}\cdot\mathbb I^{N}\left(\rho U_l\tens{M}\right)\right)_{ijk}.
\end{split}
\end{equation}

\bibliography{mybibfile}

\begin{thebibliography}{10}
\expandafter\ifx\csname url\endcsname\relax
  \def\url#1{\texttt{#1}}\fi
\expandafter\ifx\csname urlprefix\endcsname\relax\def\urlprefix{URL }\fi
\expandafter\ifx\csname href\endcsname\relax
  \def\href#1#2{#2} \def\path#1{#1}\fi

\bibitem{1981:Hirt}
C.~W. Hirt, B.~D. Nichols, Volume of fluid ({VOF}) method for the dynamics of
  free boundaries, Journal of computational physics 39~(1) (1981) 201--225.

\bibitem{1980:Nichols}
B.~Nichols, C.~Hirt, R.~Hotchkiss, Sola-vof: A solution algorithm for transient
  fluid flow with multiple free boundaries, Tech. rep., Los Alamos Scientific
  Lab. (1980).

\bibitem{1994:Sussman}
M.~Sussman, P.~Smereka, S.~Osher, A level set approach for computing solutions
  to incompressible two-phase flow, Journal of Computational physics 114~(1)
  (1994) 146--159.

\bibitem{1995:Adalsteinsson}
D.~Adalsteinsson, J.~A. Sethian, A fast level set method for propagating
  interfaces, Journal of computational physics 118~(2) (1995) 269--277.

\bibitem{1999:Jacqmin}
D.~Jacqmin, Calculation of two-phase navier--stokes flows using phase-field
  modeling, Journal of Computational Physics 155~(1) (1999) 96--127.

\bibitem{2003:Badalassi}
V.~Badalassi, H.~Ceniceros, S.~Banerjee, Computation of multiphase systems with
  phase field models, Journal of computational physics 190~(2) (2003) 371--397.

\bibitem{1997:Shen}
J.~Shen, Pseudo-compressibility methods for the unsteady incompressible
  {N}avier--{S}tokes equations, in: Proceedings of the 1994 Beijing symposium
  on nonlinear evolution equations and infinite dynamical systems, 1997, pp.
  68--78.

\bibitem{2012:Ferrer}
E.~Ferrer, R.~H. Willden, A high order discontinuous {G}alerkin--{F}ourier
  incompressible {3D} {N}avier--{S}tokes solver with rotating sliding meshes,
  Journal of Computational Physics 231~(21) (2012) 7037--7056.

\bibitem{2013:Karniadakis}
G.~Karniadakis, S.~Sherwin, Spectral/hp element methods for computational fluid
  dynamics, Oxford University Press, 2013.

\bibitem{2014:Ferrer}
E.~Ferrer, D.~Moxey, R.~Willden, S.~Sherwin, Stability of projection methods
  for incompressible flows using high order pressure--velocity pairs of same
  degree: continuous and discontinuous {G}alerkin formulations, Communications
  in Computational Physics 16~(3) (2014) 817--840.

\bibitem{2017:Ferrer}
{\relax E. Ferrer}, {\relax An interior penalty stabilised incompressible
  Discontinuous {G}alerkin--{F}ourier solver for implicit Large Eddy
  Simulations}, Journal of Computational Physics 348 (2017) 754--775.

\bibitem{2016:Cox}
C.~Cox, C.~Liang, M.~W. Plesniak, A high--order solver for unsteady
  incompressible {N}avier--{S}tokes equations using the flux reconstruction
  method on unstructured grids with implicit dual time stepping, Journal of
  Computational Physics 314 (2016) 414--435.

\bibitem{2016:gassner}
{\relax G.J. Gassner, A.R. Winters and D.A. Kopriva}, {\relax Split form nodal
  discontinuous Galerkin schemes with Summation-By-Parts property for the
  compressible Euler equations}, Journal of Computational Physics, in Press.

\bibitem{2017:Gassner}
{\relax G.J. Gassner, A. Winters, Andrew, F. Hindenlang, D.A. Kopriva}, The
  {BR1} scheme is stable for the compressible {N}avier--{S}tokes equations,
  Journal of Scientific Computing 77~(1) (2018) 154--200.

\bibitem{2017:Bassi}
F.~Bassi, F.~Massa, L.~Botti, A.~Colombo, Artificial compressibility {G}odunov
  fluxes for variable density incompressible flows, Computers \& Fluids 169.

\bibitem{1997:Bassi}
F.~Bassi, S.~Rebay, A high-order accurate discontinuous finite element method
  for the numerical solution of the compressible {N}avier--{S}tokes equations,
  Journal of computational physics 131~(2) (1997) 267--279.

\bibitem{2014:Carpenter}
M.~H. Carpenter, T.~C. Fisher, E.~J. Nielsen, S.~H. Frankel, Entropy stable
  spectral collocation schemes for the {N}avier--{S}tokes equations:
  Discontinuous interfaces, SIAM Journal on Scientific Computing 36~(5) (2014)
  B835--B867.

\bibitem{2006:Kopriva}
D.~A. Kopriva, Metric identities and the discontinuous spectral element method
  on curvilinear meshes, Journal of Scientific Computing 26~(3) (2006) 301.

\bibitem{2003:Tadmor}
E.~Tadmor, Entropy stability theory for difference approximations of nonlinear
  conservation laws and related time-dependent problems, Acta Numerica 12
  (2003) 451--512.

\bibitem{2013:Fisher}
T.~C. Fisher, M.~H. Carpenter, High-order entropy stable finite difference
  schemes for nonlinear conservation laws: Finite domains, Journal of
  Computational Physics 252 (2013) 518--557.

\bibitem{1980:Williamson}
{\relax J.H. Williamson}, Low-storage {R}unge--{K}utta schemes, Journal of
  Computational Physics.

\bibitem{2009:Toro}
{\relax E. Toro}, {\relax Riemann solvers and numerical methods for fluid
  dynamics}, Springer, 2009.

\bibitem{2019:Hinderlang}
F.~J. Hindenlang, G.~J. Gassner, D.~A. Kopriva, Stability of wall boundary
  condition procedures for discontinuous {G}alerkin spectral element
  approximations of the compressible {E}uler equations, arXiv preprint
  arXiv:1901.04924.

\bibitem{1948:Kovasznay}
L.~I.~G. Kovasznay, Laminar flow behind a two-dimensional grid, Mathematical
  Proceedings of the Cambridge Philosophical Society 44~(1) (1948) 58--62.

\bibitem{2000:Guermond}
J.-L. Guermond, L.~Quartapelle, A projection {FEM} for variable density
  incompressible flows, Journal of Computational Physics 165~(1) (2000) 167 --
  188.

\bibitem{2018:Manzanero}
{\relax J. Manzanero, A.M. Rueda--Ram\'irez, G. Rubio and E. Ferrer}, {\relax
  The Bassi Rebay 1 scheme is a special case of the Symmetric Interior Penalty
  formulation for discontinuous Galerkin discretisations with Gauss--Lobatto
  points}, Journal of Computational Physics 363 (2018) 1 -- 10.

\bibitem{2019:Manzanero}
J.~Manzanero, G.~Rubio, D.~A. Kopriva, E.~Ferrer, E.~Valero, A free--energy
  stable nodal discontinuous {G}alerkin approximation with summation-by-parts
  property for the {C}ahn--{H}illiard equation, arXiv preprint
  arXiv:1902.08089.

\bibitem{2017:Flad}
{\relax David Flad and Gregor Gassner}, On the use of kinetic energy preserving
  {DG}-schemes for large eddy simulation, Journal of Computational Physics
  350~(Supplement C) (2017) 782 -- 795.

\bibitem{2017:Fernandez}
P.~Fernandez, N.-C. Nguyen, J.~Peraire, Subgrid-scale modeling and implicit
  numerical dissipation in dg-based large-eddy simulation, in: 23rd AIAA
  Computational Fluid Dynamics Conference, 2017, p. 3951.

\bibitem{2018:Manzanero-role}
J.~Manzanero, E.~Ferrer, G.~Rubio, E.~Valero, On the role of numerical
  dissipation in stabilising under-resolved turbulent simulations using
  discontinuous {G}alerkin methods, arXiv preprint arXiv:1805.10519.

\bibitem{1958:Cahn}
J.~W. Cahn, J.~E. Hilliard, Free energy of a nonuniform system. {I}.
  interfacial free energy, The Journal of chemical physics 28~(2) (1958)
  258--267.

\bibitem{2008:Kennedy}
{\relax C.A. Kennedy and A. Gruber}, Reduced aliasing formulations of the
  convective terms within the {N}avier--{S}tokes equations for a compressible
  fluid, Journal of Computational Physics 227 (2008) 1676--1700.

\end{thebibliography}

\end{document}